\documentclass[11pt, reqno]{amsart}
\usepackage{amssymb}
\usepackage{hyperref}
\usepackage{tikz}
\usepackage{tikz-cd}
\usepackage[margin=1in]{geometry}
\usepackage{array}
\usepackage{enumitem}

\usepackage[backend=biber,style=alphabetic,maxnames=99,minnames=3]{biblatex}
\usepackage{listings}
\usepackage{xcolor}

\definecolor{codegreen}{rgb}{0,0.6,0}
\definecolor{codegray}{rgb}{0.5,0.5,0.5}
\definecolor{codepurple}{rgb}{0.58,0,0.82}
\definecolor{backcolour}{rgb}{0.95,0.95,0.92}

\lstdefinestyle{mystyle}{
    backgroundcolor=\color{backcolour},   
    commentstyle=\color{codegreen},
    keywordstyle=\color{magenta},
    numberstyle=\tiny\color{codegray},
    stringstyle=\color{codepurple},
    basicstyle=\ttfamily\lineskiplimit0pt\footnotesize,
    breakatwhitespace=false,         
    breaklines=true,                 
    captionpos=b,                    
    keepspaces=true,                 
    numbers=left,                    
    numbersep=5pt,                  
    showspaces=false,                
    showstringspaces=false,
    showtabs=false,                  
    tabsize=1,
    aboveskip=0pt,
    belowskip=0pt,
    lineskip=0pt,
}
\newtheorem{theorem}{Theorem}[section]
\newtheorem*{theorem*}{Theorem}
\newtheorem{lemma}[theorem]{Lemma}

\newtheorem{corollary}[theorem]{Corollary}
\newtheorem{proposition}[theorem]{Proposition}

\theoremstyle{definition}
\newtheorem{definition}[theorem]{Definition}
\newtheorem{notation}[theorem]{Notation}
\newtheorem*{question}{Question}
\newtheorem{examplex}[theorem]{Example}

\newenvironment{example}
  {\pushQED{\qed}\examplex}
  {\popQED\endexamplex}

\theoremstyle{remark}
\newtheorem*{remark}{Remark}

\newcommand\CC{\mathbb{C}}

\newcommand\PP{\mathbb{P}}
\newcommand\QQ{\mathbb{Q}}

\newcommand\ZZ{\mathbb{Z}}

\renewcommand\H{\mathcal{H}}

\newcommand\N{\mathcal{N}}

\renewcommand\SS{\mathcal{S}}

\newcommand{\brB}{\mathcal{N}_B}

\newcommand{\br}{\mathcal{N}}
\newcommand{\haus}{\mu_{\mathrm{haus}}}

\newcommand\wt{\widetilde}
\newcommand\wh{\widehat}

\newcommand\sep{^\mathrm{sep}}

\newcommand\inv{^{-1}}

\newcommand\aff{\operatorname{Aff}}

\newcommand{\bArb}{\operatorname{\overline{Arb}}}
\newcommand\Aut{\operatorname{Aut}}

\newcommand\gal{\operatorname{Gal}}
\newcommand\Gal{\operatorname{Gal}}

\newcommand\St{\operatorname{St}}

\newcommand{\bpimg}{\operatorname{\overline{pIMG}}}
\newcommand{\pimg}{\operatorname{pIMG}}
\newcommand{\bpfImg}{\operatorname{\overline{pIMG}}}
\newcommand{\pfImg}{\operatorname{pIMG}}

\newcommand{\arb}{\operatorname{Arb}}
\newcommand{\barb}{\operatorname{\overline{Arb}}}
\newcommand\st{\operatorname{St}}

\newcommand\cchar{\operatorname{char}}
\newcommand\sgn{\operatorname{sgn}}

\newcommand\cyc{\mathrm{cyc}}

\newcommand\ord{\mathrm{ord}}
\newcommand{\aut}{\mathrm{Aut}}

\newcommand\ab{\mathrm{ab}}
\newcommand\V{\mathcal{V}}

\newcommand\ve{\varepsilon}

\newcommand\Mper{\mathcal{A}}
\newcommand\Mpre{\mathcal{B}}

\newcommand\res{\!\!\mid\!\!}

\def\lp{(\!(}
\def\rp{)\!)}

\newcommand\s{\sigma}

\newcommand\ca[1]{\hat{a}_{#1}}
\newcommand\cb[1]{\hat{b}_{#1}}
\newcommand\cc[1]{\hat{c}_{#1}}

\newcommand\et{\text{\'et}}

\newcommand\piettame{\pi_1^{\et,\textrm{tame}}}
\newcommand\img{\mathrm{IMG}}

\newcommand{\mon}{\operatorname{Mon}}
\newcommand{\bmon}{\operatorname{\overline{Mon}}}
\newcommand{\Mon}{\operatorname{Mon}}
\newcommand{\bMon}{\operatorname{\overline{Mon}}}

\def\lprof{\langle\!\langle}
\def\rprof{\rangle\!\rangle}

\begin{document}

\title[Iterated Monodromy Groups]{Profinite Iterated Monodromy Groups of\\ Unicritical Polynomials}

\author{Ophelia Adams}
\address{Dept. of Mathematics\\University of Rochester\\Rochester, NY, 14627}
\email{ophelia.adams@rochester.edu}

\author{Trevor Hyde}
\address{Dept. of Mathematics and Statistics\\
Vassar College\\
Poughkeepsie, NY 12604\\
}
\email{thyde@vassar.edu}

\begin{abstract}
    Let $f(x) = ax^d + b \in K[x]$ be a unicritical polynomial with degree $d \geq 2$ which is coprime to $\cchar K$.
    We provide an explicit presentation for the profinite iterated monodromy group of $f$, analyze the structure of this group, and use this analysis to determine the constant field extension in $K(f^{-\infty}(t))/K(t)$.
\end{abstract}

\maketitle

\tableofcontents

\section{Introduction}

Let $K$ be a field and let $f(x) \in K(x)$ be a rational function of degree $d \geq 2$ coprime to $\cchar K$.
We write $f^n := f \circ f \circ \cdots \circ f$ to denote the $n$-fold composition of $f$.
Given an element $\beta$ in some extension of $K$ and a choice $K(\beta)\sep$ of separable closure, let $K(f^{-n}(\beta))$ denote the splitting field of $f^n(x) = \beta$ in $K(\beta)\sep$ over $K(\beta)$.

If $\beta$ is not post-critical, then $f^{-n}(\beta)$ contains $d^n$ distinct elements in $K(\beta)\sep$, and together the sets $f^{-n}(\beta)$ naturally carry the structure of a regular rooted $d$-ary tree $T_d^\infty$ with root $\beta$ where each $\alpha \in f^{-n}(\beta)$ is a child of $f(\alpha)$.
The absolute Galois group $\gal(K(\beta)\sep/K(\beta))$ fixes $K$, hence the coefficients of $f$, and thus permutes the nodes at each level of this tree while respecting the tree structure.
This gives us an \emph{arboreal Galois representation} 
\[
    \rho : \gal(K(\beta)\sep/K(\beta)) \to \aut(T_d^\infty).
\]
These Galois representations encode interesting arithmetic dynamical information, much of which is yet to be understood. 
Arboreal representations were first introduced and studied by Odoni~\cite{odoni1985galois,odoni1985prime,odoni1988realising} who used the Chebotarev density theorem to link these representations to the density of prime divisors in orbits.
See \cite{jones2013galois} for a survey of arboreal representations up to 2013 and \cite[Sec. 5]{benedetto2019current} for updates through 2018.

Of particular interest is when the function $f$ is \emph{post-critically finite} or PCF: when every critical point has a finite forward orbit.
In the PCF case, iterated pre-image extensions $K(f^{-n}(\beta))/K(\beta)$ have ramification uniformly constrained to a finite set of places, and the image of the arboreal Galois representation is topologically finitely generated.
There is a growing literature on arboreal representations for PCF maps, including but not limited to \cite{pilgrim2000dessins,nekrashevych2024self,bartholdi_nekrashevych, aitken2005finitely,benedetto2017large, bouw2021dynamical,adams2023dynamical, ejder2024iterated, benedetto2025specializations}.

In this paper we focus on unicritical polynomials $f(x) = ax^d + b$ with a generic base point $\beta = t$ where $t$ is transcendental over $K$.
Note that in this case, with $d$ coprime to $\cchar K$, the extensions $K(f^{-n}(t))/K(t)$ are separable, hence Galois.
Let $\arb f$ denote the image of the arboreal Galois representation of $f$ with transcendental base point and let $\barb f \subseteq \arb f$ denote the arboreal representation of $\gal(K(t)\sep/K\sep(t))$.
The groups $\arb f$ and $\barb f$ are (isomorphic to) the \emph{arithmetic profinite iterated monodromy group} and \emph{geometric profinite iterated monodromy group} of $f$, respectively.
There is a short exact sequence
\[
    1 \to \barb f \to \arb f \to \gal(\wh K_f/K) \to 1
\]
where $\wh K_f$ is the algebraic closure of $K$ in $K(f^{-\infty}(t))$, the extension formed by adjoining all iterated $f$-preimages of $t$.
We refer to $\wh K_f$ as the \emph{constant field extension} of $f$.
As with any short exact sequence, there is an outer action of $\gal(\wh K_f/K)$ on $\barb f$ given by lifting to $\arb f$ and conjugating.
To understand $\arb f$, we study the group $\barb f$, the constant field extension $\wh K_f/K$, and the outer action of $\gal(\wh K_f/K)$ on $\barb f$.

\subsection{Results}

Our first main result provides an explicit recursive topological presentation of $\barb f$ for any unicritical PCF polynomial.
First, some notation.
If $G \subseteq S_d$ is a subgroup of the symmetric group, then we write $[G]^\infty$ to denote the iterated wreath product of $G$ with itself
\[
    [G]^\infty = G \ltimes ([G]^\infty)^d,
\]
where $G$ acts on $d$-tuples $(g_1, g_2,\ldots, g_d)$ by permuting indices (see Section~\ref{section: iterated wreath products}).
We write elements of $[G]^\infty$ as $g(g_1, g_2,\ldots, g_d)$ where $g\in G$ and $g_i\in [G]^\infty$.
The element $g(g_1, g_2,\ldots, g_d)$ acts on the tree $T_d^\infty$ by $g$ on the first level and $g_i$ on the $i$th subtree.

Let $\s := (123\cdots d) \in S_d$ be a $d$-cycle and $C_d = \langle \s \rangle$. In Section~\ref{section: pimg}, we show that the tree $T_d^\infty$ may be labeled so that $\barb f \subseteq [C_d]^\infty$ when $f(x) = ax^d + b$ is unicritical of degree $d$ prime to $\cchar K$.
In Section~\ref{section: model groups} we show that the structure of the group $\barb f$ is determined by the orbit of the critical point $0$ under $f$, naturally splitting into three cases: post-critically infinite (PCI), periodic, or (strictly) preperiodic.
In the first two cases, the combinatorial structure of the orbit alone entirely determines $\barb f$, but in the preperiodic case a small arithmetic input is also required.

\begin{theorem}
\label{theorem: intro presentation}
    Let $f(x) = ax^d + b \in K[x]$ where $d$ is coprime to $\cchar K$.
    There exists a labeling of $T_d^\infty$ such that $\barb f \subseteq [C_d]^\infty$, and
    \begin{enumerate}
        \item (PCI) If $0$ has an infinite orbit under $f$, then
        \[
            \barb f = [C_d]^\infty.
        \]

        \item (Periodic) If $0$ is periodic with period $n$ under $f$, then $\barb f = \lprof a_1, a_2, \ldots, a_n \rprof$ where
        \[
            a_i =
            \begin{cases}
                \s (1,\ldots,1,a_n) & \text{if }i = 1,\\
                (1,\ldots,1,a_{i-1}) & \text{if }i \neq 1.
            \end{cases}
        \]

        \item (Preperiodic) If $0$ is strictly preperiodic under $f$, let $m < n$ be the smallest integers such that $f^{m + 1}(0) = f^{n + 1}(0)$.
        Let $\zeta_d$ be a choice of primitive $d$th root of unity in $K\sep$.
        Since $f^m(0)$ and $f^n(0)$ have the same image under $f(x) = ax^d +b$, there exists some $1 \leq \omega < d$ such that $f^n(0) = \zeta_d^\omega f^m(0)$.
        Then $\barb f = \lprof b_1, b_2,\ldots, b_n \rprof$ where
        \[
            b_i =
            \begin{cases}
                \s & \text{if }i = 1,\\
                (1,\ldots,1,b_n,1,\ldots,1,b_m) &\text{if }i = m + 1, \text{ where $b_n$ is in the $\omega$th component,}\\
                (1,\ldots,1,b_{i-1}) & \text{if }i \neq 1, m + 1.
            \end{cases}
        \]
    \end{enumerate}
\end{theorem}

The key technical result which allows us to derive these explicit recursive presentations for $\barb f$ in the periodic and preperiodic cases is the following \emph{semirigidity} result, which implies that $\barb f$ is determined by the $[C_d]^\infty$ conjugacy classes of its inertia generators.

\begin{theorem}
    \label{theorem: intro semirigidity}
    Let $f(x) = ax^d + b \in K[x]$ where $d$ is coprime to $\cchar K$.
    Suppose $f$ is PCF and that the strict forward orbit of $0$ has $n$ elements.
    Then $\barb f = \lprof c_1, c_2,\ldots, c_n\rprof \subseteq [C_d]^\infty$ where each $c_i$ is the image of an inertia generator over $f^i(0)$.
    If $c_i' \in [C_d]^\infty$ are elements such that $c_i$ is conjugate to $c_i'$ in $[C_d]^\infty$ for each $i$, then there exists an element $w \in [C_d]^\infty$ and elements $u_i \in \barb f$ such that
    \[
        wc_i'w\inv = u_i c_i u_i\inv
    \]
    for each $i$ and
    \[
        w\lprof c_1', c_2',\ldots, c_n'\rprof w\inv = \barb f.
    \]
\end{theorem}

Theorem~\ref{theorem: intro presentation} allows us to say a lot about the structure of $\barb f$.
For example, we determine the abelianization $\barb f$.
Let $\ZZ_d$ denote the additive group of $d$-adic integers.

\begin{proposition}
\label{prop: intro abelianization}
    Let $f(x) = ax^d + b \in K[x]$ where $d$ is coprime to $\cchar K$ and let
    $\barb f^\ab$ denote the abelianization of $\barb f$.
    If $f$ is PCF, let $n$ denote the length of the strict forward orbit of 0 under $f$.
    Then
    \[
        \barb f^\ab \cong
        \begin{cases}
            (\ZZ/d\ZZ)^\infty & \text{(PCI)}\\
            \ZZ_d^n & \text{(Periodic)}\\
            (\ZZ/d\ZZ)^n & \text{(Preperiodic)}
        \end{cases}
    \]
\end{proposition}

Other properties of $\barb f$ require us to further split the preperiodic case into several subcases:
\begin{enumerate}[leftmargin=5em]
    \item[$(A_1)$] $\omega \neq d/2$,
    \item[$(A_2)$] $m > 1$ and $(d,n) \neq (2,m+1)$,
    \item[$(A_3)$]  $d = 2$, $m > 2$, and $n = m + 1$,
    \item[$(B_1)$] $d > 2$, $\omega = d/2$, and $m = 1$,
    \item[$(B_2)$] $d = 2$, $m = 1$, and $n > 2$,
    \item[$(C)$] $d = 2$, $m = 2$, and $n = 3$.
    \item[$(D)$] $d = 2$, $m = 1$, and $n = 2$.
\end{enumerate}
We let $(A)$ refer to the assumption $(A_1)$,$(A_2)$, or $(A_3)$, and let $(B)$ refer to $(B_1)$ or $(B_2)$.
Note that these hypotheses exhaust all possible cases with $d \geq 2$, $1 \leq m < n$, and $1 \leq \omega < d$.
When $d = 2$, we have $\omega = d/2 = 1$ by default.
Case $(D)$ is exceptional throughout; it corresponds, up to conjugacy, to the quadratic Chebyshev polynomial $x^2 - 2$.

We also determine the Hausdorff dimension and orders of finite level truncations of $\barb f$.
Given $\ell \geq 0$, let $[C_d]^\ell$ denote the $\ell$-fold iterated wreath product of $C_d$ and let $\rho_\ell : [C_d]^\infty \to [C_d]^\ell$ denote the truncation map.
Let $[\ell]_d := 1 + \ell + \ell^2 + \ldots + \ell^{d-1}$.
Given a subgroup $H \subseteq [C_d]^\infty$ we define $\ord_\ell(H) := \ord(\rho_\ell(H))$ to be the order of the level $\ell$ truncation of $H$ and we define the \emph{Hausdorff dimension} of $H$ to be
\[
    \haus(H) := \lim_{\ell\to\infty}\frac{\log_d\ord_\ell(H)}{\log_d\ord_\ell([C_d]^\infty)} = \lim_{\ell\to\infty}\frac{\log_d\ord_\ell(H)}{[\ell]_d},
\]
provided the limit exists.

\begin{proposition}
\label{prop: intro finite level}
    Let $f(x) = ax^d + b \in K[x]$ where $d$ is coprime to $\cchar K$, let $n$ be as in Theorem~\ref{theorem: intro presentation}, and let $q_\ell$ and $r_\ell$ be the unique integers such that $\ell = q_\ell n + r_\ell$ and $0 \leq r_\ell < n$.
    Then the values of $\ord_\ell(\barb f)$ and $\haus(\barb f)$ are as listed in the table below.

    \begin{center}
    \vspace{1em}
    \begin{tabular}{|c|l|l|}
    \hline
            \rule[-.5em]{0pt}{1.5em} & $\log_d\ord_\ell(\barb f)$ & $\haus(\barb f)$ \\
    \hline
        PCI \rule[-.5em]{0pt}{1.5em} & $[\ell]_d$ & $1$\\
    \hline
        Periodic \rule[-.5em]{0pt}{1.5em} & $[\ell]_d - d^{r_\ell}[q_\ell]_{d^n} + q_\ell$ & $1 - \frac{d - 1}{d^{n} - 1}$\\
    \hline
       $(A)$\rule[-.5em]{0pt}{1.5em} & $[\ell]_d + d[\ell-n]_d - 2[\ell-n+1]_d + 2$ & $1 - \frac{1}{d^{n-1}}$ \\
    \hline
       $(B)$\rule[-.5em]{0pt}{1.5em} & $[\ell]_d + \frac{3d}{2}[\ell-n]_d - 3[\ell-n+1]_d + 3$ & $1 - \frac{3}{2d^{n-1}}$ \\
    \hline
       $(C)$\rule[-.5em]{0pt}{1.5em} & $11\cdot 2^{\ell-1} + 2$ & $\frac{11}{16}$ \\
    \hline
       $(D)$\rule[-.5em]{0pt}{1.5em} & $\ell + 1$ &  $0$\\
    \hline
    \end{tabular}
    \vspace{1em}
\end{center}
\end{proposition}

If $\ell \geq 0$, let $\wh K_{f,\ell}$ denote the algebraic closure of $K$ in $K(f^{-\ell}(t))$.
We establish a general bound on $\wh K_{f,\ell}$ which holds for all polynomials with degree coprime to $\cchar K$ and show that $\wh K_{f,\ell}$ is completely encoded within the structure of $\barb f$.
Given $g, h \in [C_d]^\infty$ and $\ell \geq 1$ we say that $g \sim_\ell h$ if $\rho_\ell(g)$ is conjugate to $\rho_\ell(h)$.

\begin{theorem}
\label{theorem: intro general constant field bound}
    Let $K$ be a field and let $f(x) \in K[x]$ be a polynomial with degree $d$ coprime to $\cchar K$.
    Let $g_\infty \in \barb f$ denote the image of an inertia generator over $\infty$ in $\barb f$ and let $\chi_\cyc : \gal(K\sep/K) \to \wh\ZZ^\times$ denote the cyclotomic character of $K$.
    Then for $1\leq \ell \leq \infty$ we have $\wh K_{f,\ell} \subseteq K(\zeta_{d^\infty})$ and $\tau \in \gal(K\sep/K)$ fixes $\wh K_{f,\ell}$ if and only if $g_\infty \sim_\ell g_\infty^{\chi_\cyc(\tau)}$ in $\barb f$.
\end{theorem}

Theorem~\ref{theorem: intro general constant field bound} reduces the analysis of $\wh K_{f,\ell}$ for polynomials to a purely group theoretic problem about $\barb f$.
We leverage our explicit recursive presentation of $\barb f$ provided by Theorem~\ref{theorem: intro presentation} to precisely determine the $\wh K_{f,\ell}$.

\begin{theorem}
\label{theorem: intro constant field}
    Let $f(x) = ax^d + b \in K[x]$ where $d$ is coprime to $\cchar K$.
    Let $m, n, \omega$ be as in Theorem~\ref{theorem: intro presentation} and cases $(A)$ through $(D)$ as described above.
    Let $1 \leq \ell \leq \infty$.
    \begin{enumerate}
        \item (PCI) $\wh K_{f,\ell} = K(\zeta_d)$.
        \item (Periodic) $\wh K_{f,\ell} = K(\zeta_{d^{\lfloor (\ell-1)/n\rfloor + 1}})$. Hence $\wh K_f = K(\zeta_{d^\infty})$.
        \item (Preperiodic) If $\ell \leq n$, then $\wh K_{f,\ell} = K(\zeta_d)$ and if $\ell > n$, then
        \[
            \wh K_{f,\ell} =
            \begin{cases}
                K(\zeta_{d^2/\gcd(d,\omega)}) & \text{if $(A)$ and either $m > 1$ or $d$ odd,}\\
                K(\zeta_{d^2/\gcd(d,\omega+d/2)}) & \text{if $(A)$ and $m = 1$ and $d$ even,}\\
                K(\zeta_{4d}) & \text{if $(B)$ or $(C)$,}\\
                K(\zeta_{2^\ell} + \zeta_{2^\ell}\inv) & \text{if $(D)$.}
            \end{cases}
        \]
        In particular, if not $(D)$, then
        \[
            K(\zeta_d) \subseteq \wh K_f \subseteq K(\zeta_{2d^2}).
        \]
    \end{enumerate}
\end{theorem}

\subsection{Related work}

Nekrashevych \cite[Sec. 6.4.2]{nekrashevych2024self} attributes the introduction of profinite iterated monodromy groups $\arb f$ to private communication with Richard Pink from the year 2000.
The first calculation of finite level truncations of $\barb f$ were carried out by Pilgrim \cite[Thm. 4.2]{pilgrim2000dessins} for a certain subfamily of dynamical Belyi polynomials.

If $K$ has characteristic 0, then one may use topological methods to analyze the group $\barb f$.
In particular, $\barb f$ is the profinite completion of the \emph{discrete iterated monodromy group} of $f$, denoted $\img f$, which is the arboreal representation of the fundamental group of $\PP^1(\CC)$ punctured at each point of the post-critical set of $f$.
Bartholdi and Nekrashevych \cite{bartholdi_nekrashevych} determined $\img f$ for all PCF quadratic polynomials, showing that these groups are determined by the kneading sequence of the polynomial.
Nekrashevych \cite[Thm. 5.5.3]{nekrashevych2024self} showed that the Julia set of $f$ may be recovered from $\img f$.

In 2013, Pink posted a series of preprints \cite{pinklift, pink2013profinite, PinkPCF} analyzing the groups $\arb f$ and $\barb f$ for quadratic polynomials $f(x) \in K[x]$ where $K$ has odd characteristic.
The first paper \cite{pinklift} takes an algebro-geometric approach, arguing that quadratic PCF polynomials can be lifted to characteristic 0 maps with the same post-critical combinatorics.
This combined with Grothendieck's comparison theorem for the tame \'etale fundamental group allowed Pink to show that $\barb f$ must be the same in both cases.

Our work is inspired by Pink's subsequent papers \cite{pink2013profinite, PinkPCF} in which he takes a purely group theoretic approach to analyzing $\barb f$ for quadratic polynomials.
One insight gleaned from this perspective is that while the groups $\img f$ depend on the kneading sequence of $f$, their closures $\barb f$ only depend on the combinatorial structure of the post-critical orbit; hence $\barb f$ is a coarser invariant of $f$ than $\img f$.
At the heart of his strategy is a semirigidity result \cite[Thm. 0.3]{PinkPCF} which we have generalized to all unicritical polynomials in Theorem \ref{theorem: intro semirigidity}.
Pink uses semirigidity to deduce the degree 2 cases of our main results: Compare Theorem \ref{theorem: intro presentation} with \cite[Thm. 2.4.1, Thm. 3.4.1]{PinkPCF}; Proposition \ref{prop: intro abelianization} with \cite[Thm. 2.2.7, Thm. 3.1.6]{PinkPCF}; Proposition \ref{prop: intro finite level} with \cite[Prop. 2.3.1, Prop. 3.3.3]{PinkPCF}; Theorem \ref{theorem: intro constant field} with \cite[Thm. 2.8.4, Thm. 3.10.5, Cor. 3.10.6]{PinkPCF}.

A number of challenges arise while generalizing Pink's semirigidity results from degree 2 to all $d\geq 2$.
For example, Pink leverages the fact that $\barb f$ is a pro-$p$ group in the degree 2 case (with $p = 2$) in a crucial way.
We circumvent this issue with an intermediate rigidity result which works uniformly for all degrees (Lemma \ref{lemma: generator rigidity}).
Even the identification of $\barb f$ and $\arb f$ with subgroups of $\Aut T_d^\infty$ requires more care; important ``character maps'', natural analogues of the sign homomorphisms, are not defined on $\Aut T_d^\infty$, nor preserved by conjugation from $\Aut T_d^\infty$. When $d>2$ we are required to make a careful choice of ``algebraic paths'' to embed $\barb f$ into $[C_d]^\infty$ (Proposition \ref{prop: odometer labeling}). This did not appear in the degree $d=2$ case, where $S_2 = C_2$ and the character maps are precisely the usual sign homomorphisms.

Our strategy for analyzing the constant field extensions $\wh K_f/K$ differs from the one taken by Pink.
Pink first determines the quotient $N(\barb f)/\barb f$, where $N(\barb f)$ is the normalizer of $\barb f$ in $[S_2]^\infty$, and the natural action of this group on the abelianization of $\barb f$; then he determines the image of $\gal(\wh K_f/K)$ in this action.
Instead we use the \emph{branch cycle lemma} (attributed to Fried \cite{FriedBCL}; see Lemma \ref{lemma: BCL}) to reduce the analysis of $\wh K_f/K$ to the group theoretic problem of determining which powers of an inertia generator at infinity $c_\infty$ are conjugate to $c_\infty$ in $\barb f$ (Theorem \ref{theorem: outer action determined by infinity}).
Our study of the normalizer $N(\barb f)$ limited to showing it acts transitively on odometers in $\barb f$ (Proposition \ref{prop: N(B) odometer criteria}).
Circumventing the analysis of the normalizer quotient provides a more efficient route to the constant field calculations.

Unfortunately, Pink did not publish his results; they are only available as arXiv preprints. 
While our work draws significant inspiration from Pink, there is no logical dependence on these unpublished papers.
We include all the details covering degree 2 for completeness.

The containment $K(\zeta_{d^\infty}) \subseteq \wh K_f$ in the periodic case of Theorem \ref{theorem: intro constant field} was anticipated by a result of Hamblen and Jones \cite[Thm. 2.1]{hamblen_jones}, building on a construction from Benedetto et al. \cite[Lem. 1.4]{ahmad2022arithmetic}.
We refine this containment to an exact determination of $\wh K_{f,\ell}$ using our group theoretic techniques in Proposition \ref{prop: periodic case constant field}.

\subsection{Directions for future work}

Much of the work on arboreal representations has focused on the extensions $K(f^{-\infty}(\beta))/K(\beta)$ where $\beta$ is algebraic over $K$.
These may be viewed as specializations of the extensions $K(f^{-\infty}(t))/K(t)$ that we study.
We have not considered specializations in this paper, but this is a natural next direction to pursue.
In \cite{benedetto2025specializations}, the authors develop a technique for analyzing these specializations via the Frattini subgroup of $\barb f$.
They focus on settings where $\barb f$ is a pro-nilpotent group and topologically finitely generated, in which case, the Frattini subgroup is especially well-behaved. Their results apply to unicritical polynomials with prime power degree.
When $d$ is divisible by at least two primes, the groups $\barb f$ will no longer be pro-nilpotent.
It is unclear what the Frattini subgroups look like in those case and how their arguments may be adapted.

Of course we would like to understand the groups $\arb f$ for more general PCF rational functions, but it is hard to predict how far our techniques will extend.
There are two main obstacles: semirigidity and the branch cycle lemma reduction of the constant field extension problem.
Semirigidity seems too strong to generalize in the same form much further beyond the unicritical case.
For example, one interesting new phenomenon that arises for $d > 2$ is the parameter $\omega$ in the preperiodic case.
When $d = 2$, we have $\omega = 1$ by default.
The structure of $\barb f$ and $\wh K_f$ depends on $\omega$ in subtle and interesting ways.
Perhaps semirigidity could be extended by accounting for the right invariants?
As for the branch cycle lemma reduction, the fact that polynomials have a totally tamely ramified fixed point, and hence a self-centralizing inertia subgroup, significantly constrains the outer action.
For example, we use it to show that the outer action is faithful. 
This is almost certainly true for rational functions, but would seem to require a different approach.
Beyond that, this distinguished self-centralizing subgroup is pro-cyclic, which further simplifies the group-theoretic analysis. 
For rational functions, we cannot expect such straightforward behavior; our single conjugacy problem likely expands to a family of entangled conjugacy problems.

\subsection{Overview}

In Section \ref{section: preliminaries} we review conjugacy in wreath products and systems of recurrences in iterated wreath products.
Here we define the notion of a \emph{system of cyclic conjugate recurrences} which are essential for semirigidity.

Section \ref{section: pimg} establishes the foundations of iterated preimage extensions and profinite iterated monodromy groups.
We formally define self-similar embeddings of iterated monodromy groups using collections of algebraic paths and construct an especially nice collection of paths for any PCF polynomial in Proposition \ref{prop: odometer labeling}.
The main result of this section is Proposition \ref{prop: geometric generators}, which provides a nice family of conjugate recurrences for topological generators of $\barb f$.

Sections~\ref{section: model groups}~and~\ref{section: further properties} are the technical heart of the paper.
In Section~\ref{section: model groups}, we define \emph{model groups} with recursive topological presentations reflecting the conjugate recurrences deduced in the previous section.
We analyze the structure of these groups in detail, ultimately leading to the semirigidity results.
This allows us to identify $\barb f$ with one of these model groups and thereby translate all the results about model groups into results about $\barb f$.

Then, in Section~\ref{section: further properties}, we study the normalizer of the model groups, the odometers, and certain distinguished elements called power conjugators. Building on these group-theoretic facts, in Section~\ref{section: constant fields} we prove our results on constant field extensions.
Here we obtain the cyclotomic bounds on $\wh K_f$ for any polynomial $K$, review the branch cycle lemma, and characterize $\wh K_{f,\ell}/K$ in terms of inertia at $\infty$ for polynomials.
We finish with a translation of structural results about $\barb f$ into a precise determination of the $\wh K_{f,\ell}/K$ for all unicritical polynomials.

\subsection{Acknowledgments}

We thank Rob Benedetto, Rafe Jones, Jamie Juul, Joe Silverman, and Tom Tucker for helpful conversations and encouragement throughout this project. This work was also supported in part by the AIM Workshop, ``Groups of Dynamical Origin''.

\section{Preliminaries}
\label{section: preliminaries}

In this section we review the construction of wreath products, characterize their conjugacy classes, and discuss systems of recursion in iterated wreath products.
Throughout this paper we repeatedly leverage the recursive structure of iterated wreath products to make inductive arguments; Proposition~\ref{prop: inductive criteria} is the essential tool that makes this work.

\subsection{Conjugacy in wreath products}
Let $d\geq 2$ be an integer.
Let $G$ and $H$ be groups and suppose that $G$ acts on the set $\{1,2,\ldots, d\}$ from the left.
The wreath product $G \wr H$ of $G$ and $H$ is the semidirect product $G\ltimes H^d$ associated to the action of $G$ on $H^d$ given by permuting coordinates.

The underlying set is $G \times H^d$, so we write elements of $G \wr H$ in the form $g(h_1,h_2,\ldots,h_d)$ where $g\in G$ and $h_i \in H$ for each $i$. 
The natural inclusion of $G$ into $G \wr H$ is given by $g\mapsto g(1,\ldots,1)$ and conjugation permutes coordinates:
\[
    (h_1,\dots,h_d)^g = g\inv (h_1,\ldots, h_d)g = (h_{g(1)}, \ldots, h_{g(d)}).
\]
The product of two general elements $g(h_1,\ldots,h_d)$ and $g'(h_1',\ldots,h_d')$ is given by
\[
    g(h_1, \ldots, h_d) \cdot g'(h_1',\ldots, h_d')
    = gg'(h_{g'(1)}h_1', \ldots, h_{g'(d)}h_d').
\]

\begin{definition}
    Given $u = g(h_1,\ldots, h_d)$ and $1\leq i \leq d$ with a $g$-orbit of length $n$, let $\pi_{u,i}$ denote the product
    \[
        \pi_{u,i} := h_{g^{n-1}(i)}\cdots h_{g(i)}h_i \in H.
    \]
    Replacing $i$ by $g^j(i)$ for any $j$  cyclically permutes the factors in this product, hence $\pi_{u,i} \sim \pi_{u,g^j(i)}$.
    Therefore the conjugacy class of $\pi_{u,i}$ only depends on the $g$-orbit of $i$.
\end{definition}

\begin{proposition}
\label{prop: general wreath conj crit}
    If $u = g(h_1,\ldots,h_d)$ and $u' = g'(h_1',\ldots,h_d')$ are elements of $G \wr H$, then $u \sim u'$ if and only if there exists an $a \in G$ and $b_i \in H$ for $1\leq i \leq d$ such that both $g' = g^a$ and $\pi_{u',i} = \pi_{u, a(i)}^{b_i}$ for all $1 \leq i \leq d$.
\end{proposition}

\begin{proof}
    Well, $u\sim u'$ if and only if there is some $v := a(b_1, \ldots, b_d) \in G \wr H$, such that $u' = u^v$, equivalently
    \begin{align*}
          v\inv u v 
          &= (b_1\inv, \ldots, b_d\inv)a\inv g(h_1,\ldots,h_d) a(b_1,\ldots,b_d)\\
          &= g^a (b_{g^a(1)}\inv h_{a(1)} b_1,\ldots, b_{g^a(d)}\inv h_{a(d)}b_d).
    \end{align*}
    Thus $u \sim u'$ if and only if there is some $v$ such that $g' = g^a$ and $h_i' = b_{g^a(i)}\inv h_{a(i)} b_i$ for each $i$.
    
    If $1 \leq i \leq d$ has a $g$-orbit of length $n$, then the $g' = g^a$ orbit of $a(i)$ also has length $n$, and conversely. Then we calculate
    \begin{align*}
        \pi_{u', i} &= h'_{g'^{n-1}(i)}\cdots h'_i\\
            &= (b_i\inv h_{a(g^a)^{n-1}(i)}b_{(g^a)^{n-1}(i)})\cdots(b_{g^a(i)}\inv h_{a(i)}b_i)\\
            &= b_i\inv h_{g^{n-1}a(i)}\cdots h_{ga(i)}h_{a(i)}b_i\\
            &= \pi_{u,a(i)}^{b_i}.
    \end{align*}
    
    Therefore $u \sim u'$ if and only if there exists some $a \in G$ and $b_i \in H$ for $1 \leq i \leq d$ such that $g' = g^a$ and $\pi_{u',i} = \pi_{u,a(i)}^{b_i}$ for each $i$.
\end{proof}

\subsection{Iterated wreath products}\label{section: iterated wreath products}

Here, we define iterated wreath products of finite groups and give topological/inductive criteria equality and conjugacy in these groups.

Let $G$ denote a group acting faithfully on the set $\{1,2,\ldots,d\}$ on the left.
In particular, $G$ is finite.

\begin{definition}
    For $\ell\geq 0$, the \emph{$\ell$th iterated wreath product} of $G$, denoted $[G]^\ell$, is defined inductively by $[G]^0 := 1$ and $[G]^\ell := G \wr [G]^{\ell-1}$ for $\ell \geq 1$.
\end{definition}

Let $T_d^\ell$ denote the regular rooted $d$-ary tree of height $\ell$, where the children of each node are labeled by the elements of $\{1,2,\ldots,d\}$.
There is a correspondence between the leaves of $T_d^\ell$ and words of length $\ell$ in the alphabet $\{1,2,\ldots,d\}$ given by reading the labels of the nodes along the unique path from the root to a given leaf.
Suppose $w = iw'$ is a word where $i \in \{1,\ldots, d\}$ and $w'$ is a word of length $\ell-1$.
If $u := g(g_1, g_2, \ldots, g_d) \in [G]^\ell$ then $u$ acts on $w$ by
\[
    u(w) = g(i)g_i(w')
\]
This action respects the tree structure of $T_d^\ell$ under the correspondence between words of length $\ell$ and leaves of $T_d^\ell$.
Since $G$ acts faithfully on $\{1,2,\ldots,d\}$, the group $[G]^\ell$ acts faithfully on $T_d^\ell$.
In particular, an element of $[G]^\ell$ is uniquely determined by its action on $T_d^\ell$.

Restricting an element of $[G]^\ell$ to words of length $\ell-1$ gives an element of $[G]^{\ell-1}$ and defines a surjective homomorphism $ [G]^\ell \twoheadrightarrow [G]^{\ell-1}$.
Hence the groups $[G]^\ell$ form an inverse system
\begin{equation}
\label{eqn inverse system}
    1 = [G]^0 \twoheadleftarrow [G]^1 \twoheadleftarrow [G]^2 \twoheadleftarrow [G]^3 \twoheadleftarrow\ldots.
\end{equation}
\begin{definition}
    The \emph{iterated wreath product} of $G$, denoted $[G]^\infty$, is the inverse limit of \eqref{eqn inverse system}. 
\end{definition}

The trees $T_d^\ell$ likewise form an inverse system whose limit we denote $T_d^\infty$. 
Since $G$ is finite, so too are the groups $[G]^\ell$ and the trees $T_d^\ell$.
Hence $[G]^\infty$ and $T_d^\infty$ are both naturally endowed with a profinite topology and the group $[G]^\infty$ acts continuously on $T_d^\infty$, or equivalently on right-infinite words $i_1i_2i_3\ldots$ in the alphabet $\{1,2,\ldots,d\}$.

We require a little more language to discuss iterated wreath products and their subgroups effectively:
\begin{definition}
    Let $\rho_\ell: [G]^\infty \twoheadrightarrow [G]^\ell$ denote the restriction to words of length $\ell$.
    The kernel of $\rho_\ell$, called the \emph{level \(\ell\) stabilizer} and denoted $\st_\ell [G]^\infty$, is the subgroup of all elements of $[G]^\infty$ which stabilize the first $\ell$ levels of the tree.
    
    Given two elements $u, v \in [G]^\infty$ and an integer $\ell\geq 0$ we write $u =_\ell v$ as a shorthand for $\rho_\ell(u) = \rho_\ell(v)$.
    We write $u \sim_\ell v$ as a shorthand for $\rho_\ell(u) \sim \rho_\ell(v)$ in $[G]^\ell$.
    If $U, V \subseteq [G]^\infty$ are subgroups, then we define $U =_\ell V$ and $U \sim_\ell V$ analogously.
    Given an integer $\ell\geq 0$ and an element $u \in [G]^\infty$, let $\ord_\ell(u)$ denote the order of $\rho_\ell(a)$ in the group $[G]^\ell$.
    If $K \subseteq H \subseteq [G]^\infty$, then we define
    \[
        [H : K]_\ell := [\rho_\ell(H) : \rho_\ell(K)] = [H\st_\ell [G]^\infty : K\st_\ell [G]^\infty].
    \]
\end{definition}

The recursive structure of the iterated wreath product makes it amenable to inductive arguments. Some subgroups have a similar property, including those we will study in subsequent sections,
\begin{definition}
A subgroup $H \subseteq [G]^\infty$ is said to be \emph{self-similar} if $\st_1 H \subseteq H^d$.
Note that self-similarity is not generally stable under conjugation, which is to say that self-similarity depends on the labeling of the tree.
\end{definition}

By construction, the stabilizers $\St_\ell[G]^\infty$ form a neighborhood basis of the identity. With the associated truncatiosn $\rho_\ell$, they furnish a convenient criterion for equality of closed sets:
\begin{lemma}
\label{lemma: pre-inductive criteria}
    Suppose that $U, V \subseteq [G]^\infty$ are closed subsets.
    Then
    \begin{enumerate}
        \item $U = V$ if and only if $U =_\ell V$ for all $\ell \geq 0$,
        \item $U \sim V$ if and only if $U \sim_\ell V$ for all $\ell \geq 0$.
    \end{enumerate}
\end{lemma}

\begin{proof}
    The forward implications of (1) and (2) follow immediately by taking quotients.
    We now consider the reverse implications.

    (1)
    If $\rho_\ell(U) = \rho_\ell(V)$ for all $\ell\geq 0$, then the definition of the profinite topology implies that $\overline{U} = \overline{V}$.
    Since $U$ and $V$ are closed by assumption, we conclude that $U = V$.
    
    (2)
    Let $\wt{H}_\ell \subseteq [G]^\ell$ be the stabilizer of $\rho_\ell(U)$ under the conjugation action, and let 
    \[
        H_\ell := \rho_\ell^{-1}(\wt{H}_\ell) \subseteq [G]^\infty.
    \]
    Note that $[G]^\ell$ finite implies that $\wt H_\ell$ is closed, hence compact, and of finite index.
    
    Let $M_\ell \subseteq [G]^\infty$ be the set
    \[
        M_\ell := \{m \in [G]^\infty : m\inv Um =_\ell V\}.
    \]
    Each $M_\ell$ is a union of cosets of $H_\ell$, necessarily finite, and hence compact. The $M_\ell$ are nested
    \[
        M_0 \supseteq M_1 \supseteq M_2 \supseteq \ldots,
    \]
    and the assumption $U \sim_\ell V$ implies that $M_\ell$ is nonempty for each $\ell \geq 0$.
    
    Therefore, the intersection $\bigcap_{\ell\geq 0} M_\ell$ is nonempty and for any $m \in \bigcap_{\ell\geq 0} M_\ell$, we have $m\inv Um =_\ell V$ for all $\ell\geq 0$. Thus by (1) we have $m\inv Um = V$.
\end{proof}

The following proposition provides inductive criteria for checking equality and conjugacy in iterated wreath products. 
We make frequent use of these criteria throughout the paper.

\begin{proposition}
\label{prop: inductive criteria}
    If $u, v \in [G]^\infty$ are elements and $U, V \subseteq [G]^\infty$ are closed subgroups, then
    \begin{enumerate}
        \item $u = v$ if and only if $u =_\ell v$ for all $\ell\geq 0$,
        \item $U = V$ if and only if $U =_\ell V$ for all $\ell\geq 0$,
        \item $u \sim v$ if and only if $u \sim_\ell v$ for all $\ell\geq 0$,
        \item $U \sim V$ if and only if $U \sim_\ell V$ for all $\ell\geq 0$.
    \end{enumerate}
\end{proposition}

\begin{proof}
    This is an immediate consequence of Lemma~\ref{lemma: pre-inductive criteria}; note singleton subsets are closed in the Hausdorff group $[G]^\infty$.
\end{proof}

\begin{remark}
    In Proposition~\ref{prop: inductive criteria}, the assumption that the subgroups $U$ and $V$ are closed is essential.
    For example, in Bartholdi and Nekrashevych's resolution of the twisted rabbit problem~\cite{twistedrabbit}, the authors identify three distinct subgroups of $[S_2]^\infty$ which nevertheless coincide at every finite level truncation. These are the iterated monodromy groups of the rabbit, co-rabbit, and airplane, and
    Proposition~\ref{prop: inductive criteria} implies that they have the same closure.
\end{remark}

\subsection{Systems of recurrences}
One simple way to construct elements of $[G]^\infty$ is via recursive relations. 
Suppose $x_1, x_2, \ldots, x_n$ is a list of $n$ indeterminates, where $n$ may be infinite, and let $[G]^\infty(x_1, x_2, \ldots, x_n)$ denote the group formed by freely adjoining the $x_i$ to $[G]^\infty$.

\begin{definition}
    A \emph{system of recurrences} in $[G]^\infty$ is a list of $n$ equations
    \[
        x_i = g_i(h_{i,1},\ldots, h_{i,d}),
    \]
    where $g_i \in G$ and $h_{i,j} \in [G]^\infty(x_1,x_2,\ldots, x_n)$.
\end{definition}

\begin{example}
\label{ex sys recurs}
    Let $G = S_3$ and $n = 2$.
    Let $\s := (123)$ and $\tau := (23)$ be elements of $S_3$.
    Then
    \begin{align*}
        x_1 &= \s(1,x_1,x_2)\\
        x_2 &= \tau(x_1x_2, 1,1)
    \end{align*}
    is an example of a system of recurrences with two equations.
    The system of recurrences completely determines how a solution $(x_1,x_2) = (a_1,a_2)$ acts on words.
    For example, we calculate
    \begin{align*}
        a_1(21131)
        &= 3a_2(1131)\\
        &= 31 a_1a_2(131)\\
        &= 31 a_1(1a_1a_2(31))\\
        &= 312 a_1a_1a_2(31)\\
        &= 312 a_1a_1(21)\\
        &= 312 a_1(3a_21)\\
        &= 3121 a_2(1)\\
        &= 31211.
    \end{align*}
    Since elements of $[S_3]^\infty$ are completely determined by how they act on right-infinite words, it follows that the system of congruences has a unique solution.
    Lemma~\ref{lemma: recursive solution} generalizes this observation.
\end{example}

\begin{lemma}
    \label{lemma: recursive solution}
    Any system of recurrences $x_i = g_i(h_{i,1}, h_{i,2},\ldots, h_{i,d})$ in $[G]^\infty$ has a unique solution $x_i = a_i \in [G]^\infty$ for $1 \leq i \leq n$.
\end{lemma}

\begin{proof}
    For each $\ell \geq 0$, let $A_\ell \subseteq ([G]^\infty)^n$ denote the set of all $n$-tuples $(a_1, \ldots, a_n)$ such that $x_i =_\ell a_i$ satisfies the system of recursions in $[G]^\ell$.
    The $A_\ell$ are compact and nested
    \[
        A_0 \supseteq A_1 \supseteq A_2 \supseteq \ldots.
    \]
    We prove by induction on $\ell$ that $A_\ell$ is non-empty and that if $(a_1, a_2, \ldots, a_n), (a_1', a_2',\ldots, a_n') \in A_\ell$ are two elements, then $a_i =_\ell a_i'$ for $1 \leq i \leq n$.

    Since $g =_0 1$ for all $g \in [G]^\infty$, it follows that $(1,\ldots,1) \in A_0$ and that any $(a_1,\ldots, a_n) \in A_0$ satisfies $a_i =_0 1$.
    Now suppose $\ell \geq 1$ and that our assertion holds for $\ell - 1$.
    Let $(\tilde{a}_1, \ldots, \tilde{a}_n) \in A_{\ell-1}$ and define $(a_1, \ldots, a_n)$ by
    \[
        a_i := g_i(h_{i,1}(\tilde{a}), \ldots, h_{i,d}(\tilde{a}),
    \]
    where $h_{i,j}(\tilde{a}) \in [G]^\infty$ is the element we get by substituting $x_k \mapsto \tilde{a}_k$ into $h_{i,j}$ for each $k$.
    Define $h_{i,j}(a)$ analogously via the substitution $x_k \mapsto a_k$.
    Since $(\tilde{a}_1,\ldots, \tilde{a}_n) \in A_{\ell-1}$, we have
    \[
        \tilde{a}_i =_{\ell-1} a_i = g_i(h_{i,1}(\tilde{a}), \ldots, h_{i,d}(\tilde{a})).
    \]
    
    Thus, $h_{i,j}(a) =_{\ell-1} h_{i,j}(\tilde{a})$ for all $i$ and $j$, implying that
    \[
        a_i = g_i(h_{i,1}(\tilde{a}), \ldots, h_{i,d}(\tilde{a})) =_{\ell} g_i(h_{i,1}(a), \ldots, h_{i,d}(a)).
    \]
    Hence, $(a_1,\ldots, a_n) \in A_{\ell}$.
    If $(a_1', \ldots, a_n') \in A_{\ell}$ is another element, then our inductive hypothesis implies that $a_i =_{\ell-1} a_i'$.
    Thus $h_{i,j}(a) =_{\ell-1} h_{i,j}(a')$ for each $i$ and $j$, and
    \[
        a_i =_{\ell} g_i(h_{i,1}(a), \ldots, h_{i,d}(a))
        =_{\ell} g_i(h_{i,1}(a'), \ldots, h_{i,d}(a'))
        =_{\ell} a_i'.
    \]
    This completes our inductive step and hence our induction.
    Therefore, $A := \bigcap_{\ell\geq 0}A_\ell$ is nonempty and Proposition~\ref{prop: inductive criteria} implies that $A$ contains a unique solution of the system of recursions.
\end{proof}

In the application we are working towards, we encounter elements of an iterated wreath product which satisfy a system of recurrences up to an unknown conjugacy.

\begin{definition}
    A \emph{system of $n$ conjugate recurrences} in $[G]^\infty$, where $0 \leq n \leq \infty$, is a list of conjugation identities
    \[
        x_i \sim g_i(h_{i,1}, \ldots, h_{i,d}),
    \]
    one for each $1\leq i \leq n$, where $g_i \in G$ and $h_{i,j} \in [G]^\infty(x_1, \ldots, x_n)$. 
    We say a system of conjugate recurrences is \emph{cyclic} if for each $1 \leq i \leq n$ and each $1\leq j \leq d$ there exists integers $k, m_k$ with $1 \leq k \leq n$ such that $\pi_{x_i,j} \sim x_k^{m_k}$ in $[G]^\infty(x_1,\ldots,x_n)$.
\end{definition}

A solution to a system of conjugate recurrences is a solution to any particular system of recursions obtained from a choice of conjugating elements for each $i$.
Once conjugating elements are chosen, we get a system of recurrences which has a unique solution by Lemma~\ref{lemma: recursive solution}.
When the system of conjugate recurrences is cyclic, the following proposition shows that all solutions of the system are themselves conjugate in $[G]^\infty$.

\begin{proposition}
\label{prop: weak recursion}
    Let $x_i \sim g_i(h_{i,1},\ldots,h_{i,d})$ be a system of $n$ conjugate cyclic recurrences in $[G]^\infty$.
    If $x_i = a_i$ for $1\leq i \leq n$ is one solution of this system, then $x_i = a_i'$ is another solution if and only if $a_i \sim a_i'$ in $[G]^\infty$ for $1 \leq i \leq n$.
    In other words, conjugate cyclic recurrences have unique solutions up to conjugacy.
\end{proposition}

\begin{remark}
    It is not hard to see that if $(a_1, \ldots, a_n)$ is a solution of a system of conjugate recurrences, then so is $(ua_1u^{-1}, \ldots, ua_nu^{-1})$ for any $u \in [G]^\infty$.
    Proposition~\ref{prop: weak recursion} shows something stronger: 
    If the system of conjugate recurrences is cyclic, then we can conjugate the $a_i$ independently for each $i$ and still get a solution of the system. 
\end{remark}

\begin{proof}
First suppose that $x_i = a_i$ for $1\leq i \leq n$ is a solution of the system of conjugate cyclic recurrences and that $a_i'$ are elements such that $a_i \sim a_i'$ in $[G]^\infty$ for each $i$.
Cyclicity implies that for each $1 \leq i \leq n$ and $1\leq j \leq d$ there exists integers $k, m_k$ such that $1 \leq k \leq n$ and $\pi_{a_i,j} \sim a_k^{m_k}$ in $[G]^\infty$.
Then $a_k \sim a_k'$ for all $k$ implies that $a_k^{m_k} \sim a_k'^{m_k}$.
Thus Proposition~\ref{prop: general wreath conj crit} implies
\[
    a_i' \sim a_i \sim g_i(h_{i,1}(a),\ldots,h_{i,d}(a)) \sim g_i(h_{i,1}(a'),\ldots,h_{i,d}(a')).
\]
Hence $x_i = a_i'$ for $1 \leq i \leq n$ is a solution of the system of conjugate recurrences.

Next suppose that $x_i = a_i$ and $x_i = a_i'$ are two solutions of the system of cyclic conjugate recurrences.
We wish to show that $a_i \sim a_i'$ in $[G]^\infty$ for all $1 \leq i \leq n$.
By Proposition~\ref{prop: inductive criteria} it suffices to prove by induction on $\ell$ that $a_i \sim_\ell a_i'$ for each $1\leq i \leq n$ and for all $\ell \geq 0$.
First note that $a_i =_0 a_i' =_0 1$ for each $i$, which establishes the base case.
Next suppose that $\ell \geq 1$ and that for each $i$ we have that $a_i \sim_{\ell-1} a_i'$.
Cyclicity implies that for each $1\leq i \leq n$ and $1 \leq j \leq d$,
\[
    \pi_{a_i,j} \sim a_k^{m_k} \sim_{\ell-1} a_k'^{m_k} \sim \pi_{a_i',j}.
\]
Thus Proposition~\ref{prop: general wreath conj crit} implies that 
\[
    a_i \sim g_i(h_{i,1}(a),\ldots, h_{i,d}(a)) \sim_\ell g_i(h_{i,1}(a'),\ldots, h_{i,d}(a')) \sim a_i',
\]
which completes our induction.
\end{proof}

\subsection{Iterated wreath products of cyclic groups}
For the remainder of this paper, we let $\s = (123\cdots d) \in S_d$ denote the standard $d$-cycle and let $C_d := \langle \s \rangle$ denote the cyclic subgroup generated by $\s$.
We construct the iterated wreath product $[C_d]^\infty$ with respect to the natural action of $C_d$ on $\{1,2,\ldots, d\}$.
Note that for each $\ell \geq 0$ the order of the group $[C_d]^\ell$ is $d^{[\ell]_d}$, where $[\ell]_d := \frac{d^\ell-1}{d-1}$.
Hence $\ord_\ell(g)$ divides a power of $d$ for every $g \in [C_d]^\infty$ and every $\ell \geq 0$, which implies that $g^m$ is well-defined for every $d$-adic integer $m \in \ZZ_d$.
The groups we study in this paper are all closed subgroups of $[C_d]^\infty$, hence will contain these powers.

Lemma~\ref{lemma: sigma conjugate} highlights a special case of Proposition~\ref{prop: general wreath conj crit} which is particularly useful for checking conjugacy in $[C_d]^\infty$.

\begin{lemma}
\label{lemma: sigma conjugate}
   If $\s (h_1, \ldots, h_d) \in C_d \wr H$, then
    \[
        \s (h_1, \ldots, h_d) \sim \s (1,\ldots,1, h_dh_{d-1}\cdots h_1).
    \]
\end{lemma}

The following lemma is used in the proof of Proposition~\ref{prop: geometric generators}.

\begin{lemma}
\label{lemma: d-adic unit conjugate}
    If $g \in [C_d]^\infty$ and $\ve \in \ZZ_d^\times$ satisfies $\ve \equiv 1 \bmod d$, then $g \sim g^\ve$ in $[C_d]^\infty$.
\end{lemma}

\begin{proof}
    By Proposition~\ref{prop: inductive criteria}, it suffices to prove that $g \sim_\ell g^\ve$ for all $g \in [C_d]^\infty$ and all $\ell \geq 0$.
    We proceed by induction on $\ell$.
    The base case is trivial since $g =_0 g^\ve =_0 1$.
    Now suppose $\ell\geq 1$ and $g \sim_{\ell-1} g^\ve$.
    
    If $g =_1 \s^j$, let $k$ be the number of $\s^j$-orbits in $\{1,2,\ldots,d\}$.
    Then Proposition~\ref{prop: general wreath conj crit} implies that $g$ is conjugate to
    \[
        g \sim g' := \s^j(g_1,\ldots,g_k,1,\ldots,1)
    \]
    for some $g_i \in [C_d]^\infty$.
    Note that
    \[
        g'^d = (g_1,\ldots,g_k,g_1,\ldots,g_k,\ldots,g_1,\ldots,g_k).
    \]
    By assumption we may write $\ve = 1 + d\ve'$ for some $\ve'\in \ZZ_d$.
    Then
    \begin{align*}
        g^\ve &\sim g'^\ve\\
        &= g'(g'^d)^{\ve'}\\
        &= \s^j(g_1^{1 + \ve'}, \ldots, g_k^{1 + \ve'}, g_1^{\ve'},\ldots, g_k^{\ve'},\ldots,g_1^{\ve'}, \ldots, g_k^{\ve'}, g_1^{\ve'}, \ldots, g_k^{\ve'})\\
        &\sim \s^j(g_1^{1 + d\ve'}, \ldots, g_k^{1 + d\ve'},1,\ldots,1)\\
        &= \s^j(g_1^\ve, \ldots, g_k^\ve,1,\ldots,1)
    \end{align*}
    where the second conjugacy follows from Proposition~\ref{prop: general wreath conj crit}.
    Our inductive hypothesis implies that $g_i^\ve \sim_{\ell-1} g_i$ for each $1\leq i \leq k$.
    Appealing again to Proposition~\ref{prop: general wreath conj crit} we have
    \[
        g^\ve 
        \sim \s^j(g_1^\ve, g_2^\ve,\ldots, g_k^\ve,1,\ldots,1)
        \sim_\ell \s^j(g_1,\ldots, g_k,1,\ldots,1)\
        = g'
        \sim g.
    \]
    This completes the induction.
\end{proof}

The abelianization of a wreath product $G \wr H$ is isomorphic to $G^\ab \times H^\ab$ with $g(h_1,\ldots,h_d)$ mapping to $(g, h_1\cdots h_d)$.
It follows that the abelianization of $[C_d]^\infty$ is isomorphic to the infinite direct product $C_d^\infty$.

\begin{definition}\label{definition: characters}
    Given an integer $\ell \geq 1$, the (additive) \emph{level $\ell$ character}, denoted $\chi_\ell$, is the function $\chi_\ell : [C_d]^\infty \to \ZZ/d\ZZ$ given by the composition of the abelianization map, projection onto the $\ell$th coordinate, and finally the isomorphism sending $\s^j \mapsto j \bmod d$.
\end{definition}

The functions $\chi_\ell$ may be calculated recursively by $\chi_1(\s^j(g_1,\ldots, g_d)) = j$ and for $\ell \geq 1$,
\[
    \chi_\ell(\s^j(g_1,\ldots, g_d)) = \sum_{i=1}^d \chi_{\ell-1}(g_i).
\]

\begin{remark}
    When $d = 2$, we have $S_2 = C_2$, hence $[S_2]^\infty = [C_2]^\infty$.
    In that case the functions $\chi_\ell$ are equivalent to the $\ell$th level sign functions $\sgn_\ell$, differing only by an isomorphism of their codomains.
    The sign functions play an important role in analyzing $[S_d]^\infty$ and its subgroups.
    When $d > 2$, there are also sign functions defined on $[C_d]^\infty$ by restriction from $[S_d]^\infty$, but they are less useful.
    For example, when $d$ is odd, all the sign functions are trivial on $[C_d]^\infty$, and when $d$ is even, they are $d/2$ powers of our sign functions.
    The characters $\chi_\ell$ are a finer invariant, and the appropriate generalization to this family of groups.
\end{remark}

\begin{definition}
    The \emph{standard odometer} is the element $c_\infty \in [C_d]^\infty$ defined recursively by
    \[
        c_\infty = \s(1,\ldots, 1, c_\infty).
    \]
    An \emph{odometer} is any element of $[S_d]^\infty$ that is conjugate to $c_\infty$ in $[S_d]^\infty$.
    A \emph{strict odometer} is any element of $[C_d]^\infty$ that is conjugate to $c_\infty$ in $[C_d]^\infty$.
\end{definition}

Since $\chi_1(c_\infty) = 1$ and $\chi_\ell(c_\infty) = \chi_{\ell-1}(c_\infty)$ for $\ell \geq 1$, it follows that $\chi_\ell(c_\infty) = 1$ for all $\ell\geq 1$.
Hence if $c \in [C_d]^\infty$ is any strict odometer, then we also have $\chi_\ell(c) = 1$ for all $\ell\geq 1$.
The following lemma shows that this is a sufficient condition to be an odometer.

\begin{lemma}
\label{lemma: chi odometer}
    If $c \in [C_d]^\infty$, then $c$ is a strict odometer if and only if $\chi_\ell(c) = 1$ for all $\ell \geq 1$.
\end{lemma}

\begin{proof}
    We prove by induction on $\ell$ that if $c \in [C_d]^\infty$ is an element such that $\chi_k(c) = 1$ for all $1\leq k \leq \ell$, then $c \sim_\ell c_\infty$ in $[C_d]^\infty$.
    Note that if this holds for all $\ell \geq 0$, then Proposition~\ref{prop: inductive criteria} implies that $c \sim c_\infty$ in $[C_d]^\infty$.    
    If $\ell = 1$, then $\chi_1(c) = 1$ implies that $c =_1 \s =_1 c_\infty$, and the conclusion is immediate.
    
    Now suppose that $\ell > 1$ and that our hypothesis holds for $\ell - 1$.
    Thus $c =_1 \s$ and Lemma~\ref{lemma: sigma conjugate} implies that $c$ is conjugate to $\s(1,\ldots,1,c')$ for some $c' \in [C_d]^\infty$.
    If $1 < k \leq \ell$, then $1 = \chi_k(c) = \chi_{k-1}(c')$.
    Our inductive hypothesis implies that $c' \sim_{\ell-1} c_\infty$.
    Thus by Lemma~\ref{lemma: sigma conjugate},
    \[
        c \sim_\ell \s(1,\ldots,1,c') \sim_\ell \s(1,\ldots,1,c_\infty) = c_\infty,
    \]
    where both the conjugacies are in $[C_d]^\infty$.
    This completes our induction.
\end{proof}

\begin{lemma}
\label{lemma: odometer order at finite level}
    If $c \in [S_d]^\infty$ is an odometer, then $c$ acts via a $d^\ell$-cycle on words of length $\ell$, and in particular $\ord_\ell(c) = d^\ell$.
\end{lemma}

\begin{proof}
    Identifying our alphabet $\{1,2,\ldots,d\}$ with $\{0,1,\ldots,d-1\}$ via $i \mapsto i - 1$, there is a natural bijection between right-infinite words and elements of $\ZZ_d$.
    Let $\tau : \ZZ_d \to \ZZ_d$ be the translation by 1 function, $\tau(x) := x + 1$.
    Then $\tau$ satisfies the recursion
    \[
        \tau = \s (1,\ldots,1,\tau),
    \]
    where the nontrivial restriction comes from carrying.
    Clearly $\tau$ acts via a $d^\ell$-cycle on words of length $\ell$.
    Since $c_\infty$ and $\tau$ only differ by a relabeling of the alphabet, the same must be true for $c_\infty$.
    Since all odometers are conjugate to $c_\infty$ we conclude that every odometer $c$ acts via a $d^\ell$-cycle on words of length $\ell$ and thus, $\ord_\ell(c) = d^\ell$.
\end{proof}

The essential property of odometers that we use in our analysis is that they are self-centralizing within the full tree automorphism group $[S_d]^\infty$.

\begin{notation}
    If $G$ is a topological group and $g_1, \ldots, g_n \in G$, then we write $\lprof g_1,\ldots, g_n\rprof$ to denote the subgroup of $G$ topologically generated by the $g_i$.
\end{notation}

\begin{proposition}
\label{prop: odometers are self-centralizing}
    If $c \in [S_d]^\infty$ is an odometer, $0 \leq \ell \leq \infty$, and $g \in [S_d]^\infty$ is an element such that $[g,c] =_\ell 1$, then $g \in_\ell \lprof c \rprof$.
\end{proposition}

\begin{proof}
    It suffices to prove this for $c = c_\infty$.
    Suppose $g \in [S_d]^\infty$ and $[g,c] =_\ell 1$.
    Lemma~\ref{lemma: odometer order at finite level} implies that $c$ acts on words of length $\ell$ via a $d^\ell$-cycle.
    Since $d^\ell$-cycles are self-centralizing in $S_{d^\ell}$, we conclude that $g \in_\ell \lprof c \rprof$.
    The $\ell = \infty$ case then follows from Proposition~\ref{prop: inductive criteria}.
\end{proof}

\subsection{Heisenberg group}
\label{section: heisenberg}

The \emph{Heisenberg group} $\H_d$ arises naturally as a quotient of the iterated monodromy groups we study in Section \ref{section: model groups}.

\begin{definition}
    Let $\H_d := \langle \s \rangle \ltimes (\ZZ/d\ZZ))^2$ where $\s$ acts on $(\ZZ/d\ZZ)^2$ via $\s\inv (i,j)\s = (i + j, j)$.
\end{definition}

Lemma~\ref{lemma: H_d presentation} provides a convenient presentation for $\H_d$.

\begin{lemma}
\label{lemma: H_d presentation}
    For all $d \geq 1$,
    \[
        \langle g_1, g_2 : g_1^d = g_2^d = [g_1,[g_1,g_2]] = [g_2,[g_1,g_2]] = 1\rangle \cong \H_d
    \]
    where the isomorphism sends $g_1 \mapsto \s$ and $g_2 \mapsto (0,1)$.
\end{lemma}
\begin{proof}
    Let $G := \langle g_1, g_2 : g_1^d = g_2^d = [g_1,[g_1,g_2]] = [g_2,[g_1,g_2]] = 1\rangle$.
    Consider the map from the free group generated by $g_1, g_2$ to $\H_d$ which sends $g_1 \mapsto (0,1)$ and $g_2 \mapsto \s$.
    Observe that $[g_1, g_2] \mapsto (1,0)$, which by construction commutes with $g_1$ and $g_2$.
    Since $(0,1)$ and $\s$ both have order $d$, this map factors through $G$.
    Thus it suffices to show that $G$ has order at most $d^3$.

    Since $g_1$ and $g_2$ generate $G$, the relations imply $[g_1, g_2]$ is in the center. 
    Therefore $g\in G$ may be written in the form
    \[
        g = g_1^i g_2^j [g_1,g_2]^k
    \]
    by writing $g$ as a word in $g_1$ and $g_2$ and iteratively applying $g_2g_1 = g_1g_2[g_1,g_2]\inv$. 
    It suffices to show that $[g_1,g_2]$ has order at most $d$.

    Every commutator in $G$ may be expressed as a word in conjugates of $[g_1, g_2]$, which lies in the center of $G$, hence belongs to the center of $G$.
    Therefore if $h_1, h_2, h_3 \in G$, then
    \[
        [h_1,h_2h_3] = [h_1,h_3]h_3\inv [h_1,h_2] h_3 = [h_1,h_2][h_1, h_3].
    \]
    It follows that
    \[
        [g_1,g_2]^d = [g_1,g_2^d] = [g_1,1] = 1.\qedhere
    \]
\end{proof}

The group $\H_d$ has a distinguished involution which plays an important role later.

\begin{lemma}
\label{lemma: H_d involution}
    With respect to the presentation in Lemma~\ref{lemma: H_d presentation}, there is a unique involution $\tau$ of $\H_d$ such that $\tau(g_1) = g_2$. Explicitly, in terms of $\langle \sigma \rangle \ltimes (\ZZ/d\ZZ)^2$, it is given by
    \[
        \tau(\s^r (s,t)) = \s^t (rt -s, r).
    \]
\end{lemma}
\begin{proof}
    Exchanging $g_1$ and $g_2$ in the presentation for $\H_d$ provided by Lemma~\ref{lemma: H_d presentation} yields an isomorphic group.
    Hence there is an involution $\tau: \H_d \to \H_d$ which satisfies $\tau(\s) = (0,1)$ and $\tau(0,1) = \s$; this involution is unique since $\s$ and $(0,1)$ generate $\H_d$.
    Note that $\s^r (s,t) = \s^r (0,1)^t [(0,1),\s]^s$, hence
    \[
        \tau(\s^r(s,t) = (0,1)^r \s^t [(0,1),\s]^{-s} = (0,r)\s^t(-s,0) = \s^t(rt -s, r).\qedhere
    \]
\end{proof}

\section{Profinite iterated monodromy groups}
\label{section: pimg}

Let $K$ be a field and $f$ a rational function over $K$. In this section, we construct the geometric and arithmetic profinite iterated monodromy groups of $f$ over $K$, relate them to the constant field extension, and produce embeddings into iterated wreath products which realize them as self-similar groups in a Galois-equivariant fashion. For polynomials, we produce embeddings that are especially well-behaved with respect to the ramification at infinity.

\subsection{Monodromy groups}

Let $K$ be a field, $f$ rational function over $K$, and $t$ some transcendental over $K$. Then $f$ determines a branched self-cover of $\PP^1_K$, which, in terms of function fields, corresponds to the extension extension $K(x)/K(t)$ defined by clearing denominators in $f(x)=t$. The resulting polynomial has degree $d$, and is irreducible by Gauss's lemma. We will further assume that this extension is separable; for a polynomial, $d$ coprime to $\cchar K$ is sufficient.

Fix a separable closure $(K(t))\sep$ and let $K\sep$ be the separable closure of $K$ within $(K(t))\sep$.

Three natural field extensions and their associated Galois groups appear:

\begin{definition}    
    The \emph{arithmetic monodromy group of $f$}, denoted $\mon f$ is
    \[
        \mon f := \gal(K(f^{-1}(t))/K(t)).
    \]
    The \emph{geometric monodromy group of $f$}, denoted $\bmon f$ is the subgroup
    \[
        \bmon f := \gal(K\sep(f^{-1}(t))/K\sep(t)).
    \]
    Let $\wh K_f$ be the algebraic closure of $K$ inside $K(f\inv(t))$, or equivalently $K\sep \cap K(f^{-1}(t))$. We call $\wh{K}$ the \emph{constant field extension} associated to $f$. The last group is simply $\gal(\wh K_f/K)$, the \emph{constant field Galois group}.
\end{definition}

Restriction induces a natural inclusion of $\bMon f$ into $\Mon f$. There is a further restriction from $\Mon f$ to $\Gal(\wh K_f/K)$. It is clear that $\bMon f$ is in the kernel of the latter restriction. In fact, it can be shown that the following sequence of restrictions is exact:
\[
    0\rightarrow \bmon f \rightarrow \mon f \rightarrow \gal(\wh K_f/K) \rightarrow 0.
\]

Recall that finite extensions of $\CC(t)$ correspond to finite branched covers of $\PP^1_\CC$.
In this case, the arithmetic and geometric monodromy groups of $f$ coincide and are isomorphic to the group of deck transformations of the branched cover associated to the Galois closure of $f$.
For more general fields $K$, the arithmetic monodromy group of $f$ incorporates more delicate information about the interaction between the arithmetic Galois theory of $K\sep/K$ and the geometric deck transformations of the branched cover associated to $f$.

\subsection{Iterated preimage extensions}

Let $f^n$ denote the $n$-fold composition of $f$ with itself.
The iterated preimage extensions of $K(t)$ naturally form a tower
\[
    K(t) \subseteq K(f^{-1}(t)) \subseteq K(f^{-2}(t)) \subseteq \ldots \subseteq K(f^{-\infty}(t)) := \bigcup_{\ell\geq 0} K(f^{-\ell}(t)),
\]
and similarly with $K$ replaced by $K\sep$.
With $f$ understood, let $\wh K_\ell$ be a shorthand for $\wh{K}_{f^\ell}$, and let $\wh K_\infty$ be the field such that $K(f^{-\infty}(t)) \cap K\sep(t) = \wh K_\infty (t)$.
Equivalently, $\wh K_\infty = \bigcup_{\ell\geq 0} \wh K_\ell$. The arithmetic and geometric monodromy groups for each iterate $f^n$ and the associated constant field extensions are all compatible, and passing to the limit gives rise to the profinite iterated monodromy groups we study:

\begin{definition}
    The \emph{arithmetic profinite iterated monodromy group of $f$} is the Galois group
    \[
        \pimg(f) := \gal(K(f^{-\infty}(t))/K(t)) \cong \varprojlim \mon f^n.
    \]
    The \emph{geometric profinite iterated monodromy group of $f$} is the Galois group 
    \[
        \bpimg(f) :=\gal(K(f^{-\infty}(t)/\wh{K}_f(t)) \cong \gal(K\sep(f^{-\infty}(t))/K\sep(t)) \cong \varprojlim \bmon f^n.
    \]
\end{definition}

The ``profinite'' in the names of these groups
refers to the fact that both groups, being Galois groups, are profinite. This is in contrast with the \emph{discrete} iterated monodromy group associated to a rational function defined over the complex numbers.
The discrete iterated monodromy group of $f(x) \in \CC(x)$ is constructed topologically and the geometric profinite iterated monodromy group is isomorphic to its profinite completion \cite[Ch. 5]{nekrashevych2024self}.

\subsection{Self-Similarity}

The iterated $f$-preimages of $t$ naturally form a regular rooted $d$-ary tree.
Since $f$ has coefficients in $K$, the group $\pimg f$ acts via tree automorphisms on these iterated preimages.
This is an example of an \emph{arboreal representation} of the Galois group $\pimg f$.
We will construct labelings of this tree which give this arboreal representation several useful properties.

\begin{definition}
    If $t'$ is transcendental over $K$, then a \emph{path} from $t$ to $t'$ is a $K\sep$-isomorphism $\lambda : K(t)\sep \to K(t')\sep$ such that $\lambda(t) = t'$.
\end{definition}

Note that a path from $t$ to $t'$ induces a $K\sep$-isomorphism between $K(f^{-\infty}(t))$ and $K(f^{-\infty}(t'))$.
Paths exist between any two elements transcendental over $K$, but there is typically no natural way to choose or construct an explicit path.
If $K = \CC$ and $a, b \in \CC$, then a path from $a$ to $b$ in the sense of topology induces a path from $t - a$ to $t - b$ in the algebraic sense defined above, but not all algebraic paths from $t - a$ to $t - b$ arise in this way: there are far more algebraic paths, in the same way that profinite iterated monodromy groups have far more elements than discrete iterated monodromy groups.

Since $t$ is transcendental over $K$, each $t' \in f^{-1}(t)$ is also transcendental over $K$.
If $\lambda$ is a path from $t$ to $t'$, then $f(\lambda(t)) = t$
and it follows that for all $\ell\geq 0$ and $t'' \in f^{-\ell}(t)$ we have $f^\ell(\lambda(t'')) = t''$.
Also, because $t' \in K(t)\sep$ we may suppose that $K(t')\sep \subseteq K(t)\sep$.
Hence $\lambda$ is a $K\sep$-endomorphism of $K(t)\sep$.

\begin{definition}
    A \emph{preimage labeling or ordering} (for $f$ and $t$) is a bijection $i \mapsto t_i$, where the $t_i$ are the distinct solutions of $f(x) = t$ in $K(t)\sep$.
    Given a preimage labeling, a \emph{choice of paths} $\Lambda$ for $f$ and $t$ consists of a collection of paths $\lambda_i$ from $t$ to $t_i$ for each $i$.
\end{definition}

\begin{remark}
    When $K\sep = \CC$, one may use topology to construct a convenient choices of paths for analyzing $\bpimg f$. See \cite{bartholdi_nekrashevych} for an example in the quadratic case.
    Since we are working with an arbitrary field $K$, we cannot as directly rely on topology to help us choose paths.
    We instead take a purely algebraic approach.
\end{remark}

A preimage labeling determines an identification of $\Mon f$ and $\bMon f$ with subgroups of $S_d$ with its usual permutation action. We use the choice of paths $\Lambda$ to propagate this upward, and label the tree of iterated $f$-preimages of $t$ by words in the alphabet $\{1,2,\ldots,d\}$ as follows:
If $w = i_1i_2\ldots i_\ell$ is a word in the alphabet $\{1,2,\ldots, d\}$, let $\lambda_w := \lambda_{i_1}  \lambda_{i_2} \cdots  \lambda_{i_\ell}$.
We then define $t_w := \lambda_w(t)$.
Observe that $f^\ell(t_w) = t$, hence $t_w \in f^{-\ell}(t)$. 

Suppose we have fixed a preimage labeling $i \mapsto t_i$.
Let $\gamma \mapsto \tilde\gamma$ denote the composition of the natural restriction map $\pimg f \to \mon f$ with the induced inclusion into $S_d$.
If $1 \leq i, j \leq d$, then $\tilde\gamma(i) = j$ if and only if $\gamma(t_i) = t_j$.
In this case, $\gamma$ induces a $K(t)$-isomorphism 
\[
    \gamma \res_i : K(f^{-\infty}(t_i)) \to K(f^{-\infty}(t_j)).
\]
If $\delta \in \pimg f$, then $(\gamma \delta)\res_i = \gamma\res_{\tilde\delta(i)}\delta\res_i$.

\begin{lemma}
    Let $\Lambda$ be a choice of paths.
    If $\gamma \in \pimg f$ and $1 \leq i \leq d$, let $\gamma_{\Lambda,i} = \lambda_{\tilde\gamma(i)}^{-1}\gamma \res_{i} \lambda_i$.
    Then the function $\rho_\Lambda' : \pimg f \to \mon f \ltimes (\pimg f)^d$ defined by
    \[
        \rho_\Lambda'(\gamma) := \tilde\gamma(\gamma_{\Lambda,1},\ldots, \gamma_{\Lambda,d}),
    \]
    is an injective homomorphism.
\end{lemma}

\begin{proof}
    Observe that
    \begin{align*}
        \rho_\Lambda'(\gamma)\rho_\Lambda'(\delta) &= \tilde\gamma(\gamma_{\Lambda,1},\ldots, \gamma_{\Lambda,d})\tilde\delta(\delta_{\Lambda,1},\ldots,\delta_{\Lambda,d})\\
        &= \tilde\gamma\tilde\delta(\gamma_{\Lambda,\tilde\delta(1)}\delta_{\Lambda,1},\ldots, \gamma_{\Lambda,\tilde\delta(d)}\delta_{\Lambda,d}).
    \end{align*}
    On the other hand, we have $\widetilde{\gamma\delta} = \tilde{\gamma}\tilde\delta$ and for each $i$,
    \[
        \gamma_{\Lambda,\tilde\delta(i)}\delta_{\Lambda,i} = (\lambda_{\tilde{\gamma}\tilde\delta(i)}^{-1}\gamma\res_{\tilde\delta(i)}\lambda_{\tilde\delta(i)})
        (\lambda_{\tilde\delta(i)}^{-1}\delta\res_i \lambda_{i})
        = \lambda_{\tilde{\gamma}\tilde\delta(i)}^{-1}\gamma\res_{\tilde\delta(i)}
        \delta\res_i \lambda_{i}
        = \lambda_{\wt{\gamma\delta}(i)}^{-1}(\gamma\delta)\res_i \lambda_{i}
        = (\gamma\delta)_{\Lambda,i}.
    \]
    Therefore $\rho_\Lambda'(\gamma\delta) = \rho_\Lambda'(\gamma)\rho_\Lambda'(\delta)$.
    If $\gamma \in \ker \rho_\Lambda'$, then $\tilde\gamma(i) = i$ for each $i$ and $\gamma\res_i$ is the identity on $K(f^{-\infty}(t_i))$. 
    Hence $\gamma$ acts trivially on the set $\bigcup_{\ell\geq 0} f^{-\ell}(t)$, implying $\gamma = 1$ in $\pimg f$.
    Thus $\rho_\Lambda'$ is injective.
\end{proof}

Iterating $\rho_\Lambda'$ gives a homomorphism $\rho_\Lambda : \pimg f \to [\mon f]^\infty$ defined recursively by
\[
    \rho_\Lambda(\gamma) := \tilde\gamma(\rho_\Lambda(\gamma_{\Lambda,1}),\ldots,\rho_\Lambda(\gamma_{\Lambda,d})).
\]
The following proposition shows how $\rho_\Lambda(\gamma)$ acts with respect to the labeling of the iterated pre-image tree corresponding to $\Lambda$.

\begin{lemma}
\label{lemma: iterated monodromy embedding}
    If $w$ is a word in the alphabet $\{1,2,\ldots,d\}$ and $\gamma \in \pimg f$, then
    \[
        \gamma(t_w) = t_{\rho_{\Lambda}(\gamma)(w)}.
    \]
\end{lemma}

\begin{proof}
    We proceed by induction on the length $\ell$ of $w$.
    The claim is trivially true for words of length 0.
    Now suppose that $\ell \geq 1$ and that the assertion holds for all words of length $\ell - 1$ and all $\gamma \in \pimg f$.

    Every word of length $\ell$ may be expressed as $iw$ for some $1 \leq i \leq d$ and some word $w$ of length $\ell - 1$.
    If $\gamma \in \pimg f$, then there is some word $w'$ with the same length as $w$ such that $\gamma(t_{iw}) = t_{\tilde\gamma(i)w'}$.
    Since $t_{iw} = \lambda_i(t_w)$ and $t_{\tilde\gamma(i)w'} = \lambda_{\tilde\gamma(i)}(t_{w'})$, it follows that
    \[
        \gamma_{\Lambda,i}(t_w) = \lambda_{\tilde\gamma(i)}^{-1}\gamma\res_i\lambda_i(t_w) = t_{w'}.
    \]
    On the other hand, our inductive hypothesis implies that
    $
        \gamma_{\Lambda,i}(t_w) = t_{\rho_\Lambda(\gamma_{\Lambda,i})}.
    $
    Therefore $w' = \rho_\Lambda(\gamma_{\Lambda,i})$.
    Hence
    \[
        \gamma(t_{iw}) = t_{\tilde\gamma(i)\rho_\Lambda(\gamma_{\Lambda,i})(w)} = t_{\rho_\Lambda(\gamma)(iw)},
    \]
    which completes our induction.
\end{proof}

\begin{lemma}
\label{lemma: choice of paths}
    Fix a preimage labeling and suppose $\Lambda_1$ and $\Lambda_2$ are two choices of paths.
    There exists an element $w \in \st_1 [\mon f]^\infty = ([\Mon f]^\infty)^d$ such that $w\inv\rho_{\Lambda_1}w = \rho_{\Lambda_2}$.
\end{lemma}

\begin{proof}
    Proposition~\ref{prop: inductive criteria} implies that it suffices to prove that for all $\ell \geq 0$ there exists a $w \in \st_1 [\mon f]^\infty$ such that $w\inv\rho_{\Lambda_1}w =_\ell \rho_{\Lambda_2}$.
    We proceed by induction on $\ell$; the $\ell = 0$ case is immediate.

    Suppose that $\ell \geq 1$ and that the claim is true for $\ell-1$.
    Let $w \in \st_1[\mon f]^\infty$ be an element such that $w\inv\rho_{\Lambda_1}w =_{\ell-1} \rho_{\Lambda_2}$.
    Let $w_1 := (w, \ldots, w) \in \st_1[\mon f]^\infty$.
    Then for $\gamma \in \pimg f$ we have
    \begin{align*}
        w_1\inv \rho_{\Lambda_1}(\gamma) w_1 
        &= \tilde\gamma(w\inv\rho_{\Lambda_1}w(\gamma_{\Lambda_1,1}),\ldots,w\inv\rho_{\Lambda_1}w(\gamma_{\Lambda_1,d}))\\
        &=_\ell \tilde\gamma(\rho_{\Lambda_2}(\gamma_{\Lambda_1,1}),\ldots,\rho_{\Lambda_2}(\gamma_{\Lambda_1,d})).
    \end{align*}
    Let $\lambda_{1,i}$ and $\lambda_{2,i}$ denote the paths associated to $\Lambda_1$ and $\Lambda_2$ respectively.
    Define
    \[
        w_2 := (\rho_{\Lambda_2}(\lambda_{1,1}\lambda_{2,1}\inv),\ldots,\rho_{\Lambda_2}(\lambda_{1,d}\lambda_{2,d}\inv)) \in \st_1 [\mon f]^\infty,
    \]
    and set $w' := w_1w_2 \in \st_1 [\mon f]^\infty$.
    Then
    \begin{align*}
        {w'}\inv\rho_{\Lambda_1}(\gamma)w' 
        &=_\ell w_2\inv \tilde\gamma(\rho_{\Lambda_2}(\gamma_{\Lambda_1,1}),\ldots,\rho_{\Lambda_2}(\gamma_{\Lambda_1,d})) w_2\\
        &= \tilde\gamma
        (\rho_{\Lambda_2}(\lambda_{2,\tilde\gamma(1)}\lambda_{1,\tilde\gamma(1)}^{-1}\gamma_{\Lambda_1,d}\lambda_{1,1}\lambda_{2,1}^{-1}),
        \ldots,
        \rho_{\Lambda_2}(\lambda_{2,\tilde\gamma(d)}\lambda_{1,\tilde\gamma(d)}^{-1}\gamma_{\Lambda_1,d}\lambda_{1,d}\lambda_{2,d}^{-1})
        )\\
        &= \tilde\gamma(\rho_{\Lambda_2}(\gamma_{\Lambda_2,1}),
        \ldots,
        \rho_{\Lambda_2}(\gamma_{\Lambda_2,d}))\\
        &= \rho_{\Lambda_2}(\gamma).
    \end{align*}
    This completes our inductive step, hence our proof.
\end{proof}

A choice of paths $\Lambda$ for $f$ and $t$ also gives an embedding
\[
    \bar{\rho}_{\Lambda} : \bpimg f \to [\bmon f]^\infty \subseteq [\mon f]^\infty.
\]

\begin{definition}
    Let $\Lambda$ be a choice of paths for $f$ and $t$.
    Then we define
    \begin{align*}
        \arb f &:= \rho_{\Lambda}(\pimg f) \subseteq [\mon f]^\infty,\\
        \barb f &:= \bar\rho_{\Lambda}(\bpimg f) \subseteq [\bmon f]^\infty.
    \end{align*}
    The dependence on the paths is suppressed in the notation.
\end{definition}

The homomorphisms $\rho_\Lambda$ and $\bar\rho_\Lambda$ are examples of \emph{arboreal Galois representations}, hence the notation $\arb f$ and $\barb f$ for their images.
With a fixed preimage labeling, Lemma~\ref{lemma: choice of paths} implies that $\arb f$ and $\barb f$ are well-defined up to conjugation by an element of $\st_1 [\bmon f]^\infty$.

\begin{lemma}
    The groups $\arb f$ and $\barb f$ are self-similar with respect to any choice of paths $\Lambda$.
\end{lemma}

\begin{proof}
    This follows immediately from the definition of $\rho_\Lambda$ and the fact that $\gamma_{\Lambda,i} := \lambda_{\tilde\gamma i}^{-1}\gamma\res_i \lambda_i \in \pimg f$ for every $\gamma \in \pimg f$.
\end{proof}

\subsection{Post-critically finite rational functions}

Let $C_f \subseteq \PP^1_{K\sep}$ denote the set of critical points of $f(x)$.
The \emph{post-critical set} of $f$, denoted $P_f$ is the strict forward $f$-orbit of $C_f$,
\[
    P_f := \bigcup_{n\geq 1} f^n(C_f).
\]

\begin{definition}
    We say $f$ is \emph{post-critically finite} or \emph{PCF} if $P_f$ is finite.
\end{definition}

Recall that there is a natural correspondence between places in $K\sep(t)$ and points in $\PP^1_{K\sep}$.
Hence when we talk about ramification or inertia groups over a point $p \in \PP^1_{K\sep}$, we mean the ramification or inertia groups over the place corresponding to $p$.
The points which ramify in $K\sep(f^{-\infty}(t))/K\sep(t)$ correspond to the critical values of iterates of $f$, often called branch points.
The chain rule implies that $P_f$ contains the critical values of all iterates of $f$.
Thus the extension $K\sep(f^{-\infty}(t))/K\sep(t)$ is only ramified at finitely many points when $f$ is PCF. If $\cchar K$ does not divide the ramification index of any critical point of $f$, then $K\sep(f^{-\infty}(t))/K\sep(t)$ is at most tamely ramified over these points, hence the inertia groups over all points are topologically cyclic. 

\begin{lemma}
\label{lemma: inertia generators}
    Suppose that $f$ is PCF and that $P_f$ contains $n + 1$ points $p_1, p_2,\ldots, p_{n+1}$. 
    Further assume that $\cchar K$ does not divide the ramification index of any critical point of $f$.
    Then
    \[
        \bpimg f = \lprof \gamma_1, \gamma_2,\ldots, \gamma_n, \gamma_{n+1}\rprof,
    \]
    where each $\gamma_i$ is a topological generator of an inertia subgroup over the point $p_i$, and $\gamma_{n+1} = \gamma_1\cdots \gamma_n$.
\end{lemma}

\begin{proof}
    Let $K\sep(t)_{P_f}/K\sep(t)$ denote the maximal tamely ramified extension of $K\sep(t)$ which is only ramified over the points in $P_f$.
    Grothendieck proved that
    \begin{equation}
    \label{eqn: grothendieck}
        \gal(K\sep(t)_{P_f}/K\sep(t)) = \lprof \tau_1, \tau_2, \ldots, \tau_n, \tau_{n+1} : \tau_{n+1} = \tau_1\cdots \tau_n\rprof,
    \end{equation}
    where each $\tau_i$ is a topological generator of an inertia subgroup over the point $p_i$ (see \cite[Expos\'e X Corollaire 3.9 and Expos\'e XII, Corollaire 5.2]{SGA}, or \cite[Thm. 4.9.1]{szamuely} for a precise statement in English).
    Since $K\sep(f^{-\infty}(t))/K\sep(t)$ is only tamely ramified over $P_f$, it follows that $\bpimg f$ is a quotient of $\gal(K\sep(t)_{P_f}/K\sep(t))$, hence that
    \[
        \bpimg f = \lprof \gamma_1, \gamma_2,\ldots, \gamma_n, \gamma_{n+1}\rprof,
    \]
    where each $\gamma_i$ topologically generates an inertia subgroup over $p_i$ and $\gamma_{n+1} = \gamma_1\cdots \gamma_n$.
\end{proof}

\begin{remark}
    The Galois group $\gal(K\sep(t)_{P_f}/K\sep(t))$, the maximal extension of \(K(t)\) unramified outside of $P_f$ may be interpreted as $\piettame(\PP^1_{K\sep} \setminus P_f, p)$, the tame \'etale fundamental group of $\PP^1_{K\sep}\setminus P_f$ with respect to the geometric point $p$ corresponding to our choice of separable closure of $K\sep(t)$ in which $K\sep(t)_{P_f}$ lives.
    This presentation of $\piettame(\PP^1_{K\sep} \setminus P_f, p)$ follows from the fact that tamely ramified extensions can be lifted to characteristic 0 where the \'etale fundamental group is known to be the profinite completion of the topological fundamental group of $\PP^1_\CC \setminus \tilde P_f$, where $\tilde P_f$ is a lift of the set $P_f$.
    This topological fundamental group is classically known to have a presentation matching \eqref{eqn: grothendieck}, where the generators correspond to loops winding once around each of the corresponding punctures.
\end{remark}

\begin{remark}
    Note that Lemma~\ref{lemma: inertia generators} only tells us that $\bpimg f$ is topologically finitely generated by inertia generators. This choice of inertia generators is not unique, and not all choices necessarily generate the group. There will typically be many intricate relations among these generators which depend on dynamical properties of $f$.
\end{remark}

The following Lemma is useful for determining the conjugacy class of $\rho_\Lambda(\gamma)$ when $\gamma$ is an inertia generator.

\begin{lemma}
\label{lemma: inertia recursion}
    Fix a preimage labeling and a choice of paths $\Lambda$ for $f$ and $t$.
    Let $(t - p)$ be a prime in $K\sep(t)$ and let $P$ be a prime in $K\sep(f^{-\infty}(t))$ over $(t - p)$.
    Let $\gamma \in \bpimg f$ be a topological generator for the inertia group of $P$ over $(t - p)$.
    
    For each $1 \leq i \leq d$, we have $P \cap K\sep(t_i) = (t_i - q_i)$ for some $q_i \in f^{-1}(p)$.
    Let $e_i$ denote the ramification index of $f$ at $q_i$; note that $e_i$ is also the length of the orbit of $i$ under $\tilde \gamma$.
    Then
    \begin{equation}
    \label{eqn: gamma product}
        \gamma_{\Lambda,\tilde\gamma^{e_i-1}(i)}\cdots \gamma_{\Lambda,i} = \lambda_i\inv \gamma^{e_i}\res_i\lambda_i
    \end{equation}
    is a topological generator for the inertia group of $P_i := \lambda_i^{-1}(P \cap K\sep(f^{-\infty}(t_i)))$ over $(t - q_i)$.
\end{lemma}

\begin{proof}
    The identity \eqref{eqn: gamma product} follows from the definition of $\gamma_{\Lambda,i}$ and the assumption that $e_i$ is the length of the $\tilde\gamma$ orbit of $i$.
    Since $e_i$ is the ramification index of $f$ at $q_i$, it follows that $\gamma^{e_i}$ topologically generates the inertia group for $\lambda_i(P_i) := P\cap K\sep(f^{-\infty}(t_i))$ over $(t_i - q_i)$.
    Therefore $\lambda_i\inv \gamma^{e_i}\res_i\lambda_i$ topologically generates the inertia group of $P_i$ over $\lambda_i^{-1}(t_i - q_i) = (t - q_i)$.
\end{proof}

\subsection{Choosing paths for polynomials}

Suppose now that $f(x) \in K[x]$ is a polynomial of degree $d$ prime to $\cchar K$.
It follows that $\infty$ is a totally tamely ramified fixed point for $f$, and hence of any extension $K(x)/K(t)$ defined by $f^n(x) = t)$, and so an inertia subgroup $\lprof\gamma_\infty\rprof$ over infinity in the Galois closure $K\sep(f^{-\infty}(t))/K\sep(t)$ is isomorphic to $\ZZ_d$.

Recall the standard odometer $c_\infty$ in $[C_d]^\infty$ defined recursively by $c_\infty = \s (1,\ldots,1,c_\infty)$. If one views the tree labeling as encoding digits of $d$-adic integers, the standard odometer represents addition by $1$ in $\ZZ_d$, so $\lprof c_\infty\rprof \cong \ZZ_d$.

Abstractly, then, inertia subgroups are isomorphic to subgroups generated by an odometer. It turns out this isomorphisms can be realized Galois-theoretically: Proposition~\ref{prop: odometer labeling} shows that for each choice of topological generator $\gamma_\infty \in \bpimg f$ of an inertia subgroup at $\infty$, there exists a choice of paths with respect to which $\gamma_\infty$ acts via the standard odometer.

\begin{proposition}
\label{prop: odometer labeling}
    For each topological generator $\gamma_\infty \in \bpimg f$ of an inertia group over $\infty$, there exists a preimage labeling $i \to t_i$ and a path $\lambda$ from $t$ to $t_1$ such that the collection of paths $\Lambda$ defined by $\lambda_i := \gamma_\infty^{i-1}\res_1\lambda$ from $t$ to $t_i$ satisfies $\rho_\Lambda(\gamma_\infty) = c_\infty$.
\end{proposition}

\begin{proof}
    Up to a linear change of coordinate defined over $K$, we may write $f(x) = ax^d + b$.
    Let $K\sep\lp1/t\rp \supseteq K\sep(t)$ be the field of formal Laurent series in $1/t$, which may be interpreted as the completion of $K\sep(t)$ at $\infty$.
    Let $K\lp1/t\rp\sep$ be a choice of separable closure; note that $K\lp1/t\rp\sep$ contains a separable closure of $K(t)$.
    Since $d$ is coprime to $\cchar K$, the extension $K\sep\lp1/t\rp(f^{-\infty}(t))/K\sep\lp1/t\rp$ is a pro-$d$ cyclic extension with Galois group canonically isomorphic to the inertia group over $\infty$ of $K\sep(f^{-\infty}(t))/K\sep(t)$.
    Hence we may extend $\gamma_\infty$ to a topological generator of $\gal(K\sep\lp1/t\rp(f^{-\infty}(t))/K\sep\lp1/t\rp)$.
    Note that $\lprof \gamma_\infty \rprof$ is isomorphic to $\ZZ_d$, the additive group of $d$-adic integers.

    Let $t_1 \in f^{-1}(t) \subseteq K\lp1/t\rp\sep$ and let $\lambda : K\lp 1/t\rp\sep \to K\lp 1/t_{1}\rp\sep$ be a $K\sep$-isomorphism such that $\lambda(t) = t_{1}$; note that $\lambda$ induces a path from $t$ to $t_{1}$.
    Furthermore, $|g| = |\lambda(g)|_1$ for all $g \in K\lp1/t\rp\sep$ where $|\cdot|$ and $|\cdot|_1$ are the standard normalized absolute values on $K\lp1/t\rp\sep$ and $K\lp1/t_1\rp\sep$ respectively.
    Let $t_{0,\ell} := \lambda^\ell(t)$.
    Then $(t_{0,\ell})$ is a tower of iterated $f$-preimages of $t$.

    Kummer theory implies that for each $\ell \geq 0$,
    \[
        K\sep\lp1/t\rp(f^{-\ell}(t)) = K\sep\lp1/t^{1/d^\ell}\rp.
    \]
    
    Considering the Newton polygon of $f^\ell(x) - t$ with respect to the $1/t$-adic valuation, it follows that each element of $f^{-\ell}(t)$ has valuation $-1/d^{\ell}$. 
    This implies that for each element $t' \in f^{-\ell}(t)$ there is a unique root of $g_\ell(x) :=a^{d^{\ell-1}}x^{d^\ell} - t$ which is closest to $t'$ with respect to the $1/t$-adic metric.
    Let $\tau_\ell$ be the root of $g_\ell(x)$ closest to $t_{0,\ell}$.
    Note that the these roots are compatible in the sense that \(a\tau_{\ell+1}^d = \tau_\ell\).

    Let $(\zeta_{d^\ell})$ be the sequence of primitive roots of unity in $K\sep$ such that 
    \[
        \gamma_\infty(\tau_\ell) = \zeta_{d^\ell}\tau_\ell.
    \]
    Compatibility of the $\tau_\ell$ implies that $\zeta_{d^{\ell+1}}^d = \zeta_{d^\ell}$.
    By construction of $\tau_\ell$, we may express $t_{0,\ell}$ as a Laurent series in $\tau_\ell$ with coefficients in $K\sep$ of the form
    \[
        t_{0,\ell} = \tau_{\ell} + \sum_{i=1-d}^\infty a_{i,\ell} \tau_{\ell}^{-i}.
    \]
    For $\ve \in \ZZ/d^\ell\ZZ$, define $t_{\ve,\ell} := \gamma_\infty^\ve t_{0,\ell}$.
    Thus the leading term of $t_{\ve,\ell}$ expanded as a Laurent series in $\tau_\ell$ is $\zeta_{d^\ell}^\ve \tau_\ell$.
    Since the leading term of 
    \[
        \lambda(t_{\ve,\ell}) = \zeta_{d^{\ell}}^\ve\lambda(t_{0,\ell}) = \zeta_{d^\ell}^\ve t_{0,\ell+1}
    \]
    is $\zeta_{d^{\ell}}^\ve\tau_{\ell+1} = \zeta_{d^{\ell+1}}^{d\ve}\tau_{\ell+1}$, it follows that
    \[
        \lambda(t_{\ve,\ell}) = t_{d\ve,\ell+1}.
    \]
    
    Define $\Lambda$ to be the collection of paths $\lambda_i$ from $t$ to $t_i$ defined by $\lambda_i := \gamma_\infty^{i-1}\res_1\lambda$.
    If $1 \leq i < d$, then
    \[
        (\gamma_\infty)_{\Lambda,i} = \lambda_{i+1}^{-1}\gamma_{\infty}\res_i\lambda_i
        = \lambda^{-1}(\gamma_\infty^i\res_1)^{-1}\gamma_\infty\res_i\gamma_\infty^{i-1}\res_1\lambda
        = 1,
    \]
    and if $i = d$, then
    \[
        (\gamma_\infty)_{\Lambda,d} = \lambda_{1}^{-1}\gamma_{\infty}\res_d\lambda_d
        = \lambda^{-1}\gamma_\infty\res_d\gamma_\infty^{d-1}\res_1\lambda
        = \lambda^{-1}\gamma_\infty^d\res_1\lambda.
    \]
    To complete the proof it suffices to show that $(\gamma_\infty)_{\Lambda,d}(t_{\ve,\ell}) = \gamma_\infty(t_{\ve,\ell})$ for all $\ell \geq 0$ and $\ve \in \ZZ/d^\ell\ZZ$.
    By construction we have $\gamma_\infty(t_{\ve,\ell}) = t_{\ve+1,\ell}$.
    On the other hand,
    \[
        \lambda^{-1}\gamma_\infty^d\res_1\lambda(t_{\ve,\ell})
        = \lambda^{-1}\gamma_\infty^d(t_{d\ve,\ell+1})
        = \lambda^{-1}(t_{d\ve + d,\ell+1})
        = t_{\ve + 1,\ell}.\qedhere
    \]    
\end{proof}

\subsection{Unicritical polynomials}
\label{section: unicritical PCF}

A polynomial $f(x) \in K[x]$ is said to be \emph{unicritical} if $f$ has a unique finite critical point.
Since $\gal(K\sep/K)$ permutes critical points of $f$, the unique finite critical point must belong to $K$.
Changing coordinates over $K$ we may assume that $f(x) = ax^d + b$ for some $a, b \in K$.

The post-critical set $P_f$ consists of $\infty$ and the forward orbit of $0$.
Let $p_i := f^i(0)$.
We say $f$ is \emph{post-critically infinite} if $0$ has an infinite orbit under $f$; otherwise $f$ is PCF.
Let $\gamma_i$ for $1 \leq i \leq n$ and $\gamma_\infty$ be topological generators for inertia groups over each $p_i$.
Since $f$ is unicritical, we have identifications $\bmon f = C_d$ and $\mon f \subseteq \aff_{1,d}$, where 
\[
    \aff_{1,d} \cong \{ rx + s : r \in \ZZ/d\ZZ^\times, s \in \ZZ/d\ZZ\}
\]
is the one dimensional affine group modulo $d$.
Note that $p_1$ is the only finite critical value of $f$, hence $\gamma_i =_1 1$ for $i > 1$.
The group $\barb f$ is as large as possible in the post-critically infinite case.

\begin{proposition}
\label{prop: PCI arb}
    Suppose $f(x) = ax^d + b$ is post-critically infinite.
    Then $\barb f = [C_d]^\infty$.
\end{proposition}

\begin{proof}
    The element $\gamma_i$ acts trivially on all levels $j < i$, acts (without loss of generality) by $\sigma$ on a single branch at level $i$, and acts trivially on each of the sub-trees above that branch.
    Such elements generate the group $[C_d]^\infty$.
    Therefore $\barb f = [C_d]^\infty$.
\end{proof}

We now turn to the case of a PCF unicritical polynomial.
Suppose that $f$ has exactly $n$ finite post-critical points.
Proposition~\ref{prop: odometer labeling} provides us with a preimage labeling and a choice of paths $\Lambda$ for $f$ and $t$ such that $\rho_\Lambda(\gamma_\infty) = c_\infty$ is the standard odometer.
Let $c_i := \rho_\Lambda(\gamma_i) \in [C_d]^\infty$. 
Thus $\barb f = \lprof c_1, c_2, \ldots, c_n\rprof$ and $c_\infty = c_1\cdots c_n$.

The identity $c_\infty = c_1 \cdots c_n$ implies that $c_1 =_1 c_\infty =_1 \s$.
If $t' = \sqrt[d]{\frac{t - b}{a}} \in f^{-1}(t)$, then there exists a primitive $d$th root of unity $\zeta_d \in K\sep$ such that $\gamma_1(t') = \gamma_\infty(t') = \zeta_dt'$.

The qualitative structure of $\barb f$ depends fundamentally on whether the finite critical point of $f$ is periodic or (strictly) preperiodic.
We refer to these as the \emph{periodic} and \emph{preperiodic} cases, respectively.
In the periodic case, $p_1$ is periodic with period $n$.
In the preperiodic case, there is an integer $1 \leq m < n$ such that $f^m(p_1) = f^n(p_1)$.
This is equivalent to saying that $f(p_m) = f(p_n) = p_{m+1}$.
Since $f$ is unicritical, there must be some integer $1\leq \omega < d$ such that $p_n = \zeta_d^\omega p_m$.

Proposition~\ref{prop: geometric generators} shows that the $c_i$ satisfy a system of cyclic conjugate recurrences in $[C_d]^\infty$ which depend on $d$ and $n$ in the periodic case, and $d$, $m$, $n$, and $\omega$ in the preperiodic case.

\begin{proposition}
\label{prop: geometric generators}
    Let $f(x) = ax^d + b$ be PCF.
    The topological generators $c_1, c_2, \ldots, c_n$ of $\barb f$ satisfy the following system of conjugate recurrences in $[C_d]^\infty$,
    
    \begin{enumerate}
        \item (Periodic Case)
        \[
            c_i \sim 
            \begin{cases}
                \s(1,\ldots,1,c_n) & \text{if }i = 1,\\
                (1,\ldots, 1, c_{i-1}) & \text{if }i \neq 1
            \end{cases}
        \]

        \item (Preperiodic Case)
        \[
            c_i \sim 
            \begin{cases}
                \s & \text{if }i = 1\\
                (1,\ldots,1,c_n,1,\ldots,1,c_m) & \text{if }i = m+1, \text{where $c_n$ is in the $\omega$th component}\\
                (1,\ldots, 1, c_{i-1}) & \text{if }i \neq 1, m + 1.
            \end{cases}
        \]
    \end{enumerate}
\end{proposition}

\begin{proof}
    (1)
    If $i \neq 1$, then $p_i$ is not a critical value of $f$.
    Hence $c_i =_1 1$.
    Thus 
    \[
        \rho_\Lambda(\gamma_i) = (\rho_\Lambda(\gamma_{\Lambda,1}),\ldots,\rho_\Lambda(\gamma_{\Lambda,d}))
    \]
    where Lemma~\ref{lemma: inertia recursion} implies that each $\rho_\Lambda(\gamma_{\Lambda,j})$ is the image of an inertia generator over $(t - q_{i,j})$ where $q_{i,j} \in f^{-1}(p_i)$.
    Our assumption that the unique finite critical point of $f$ is periodic implies that there is a unique $j$ such that $q_{i,j} = p_{i-1}$ and for all other $j'$ we have $q_{i,j'} \notin P_f$. 
    Hence $\rho_{\Lambda}(\gamma_{\Lambda,j}) \sim c_{i-1}^{\ve_{i-1}}$ in $\barb f$ for some $d$-adic unit $\ve_{i-1}$ and $\rho_\Lambda(\gamma_{\Lambda,j'}) = 1$.
    Therefore, conjugating by a power of $\s$ we have the following conjugacy in $[C_d]^\infty$,
    \[
        c_i \sim (1,\ldots,1, c_{i-1}^{\ve_{i-1}}).
    \]

    Next suppose $i = 1$.
    Then $p_1$ is a critical value of $f$ with unique preimage $p_n  = 0$.
    Thus $c_1 =_1 \s$.
    To ease notation, let us write $\gamma := \gamma_1$.
    Lemma~\ref{lemma: inertia recursion} implies that
    $
        \rho_\Lambda(\gamma_{\Lambda,d})\cdots \rho_\Lambda(\gamma_{\Lambda,1})
    $
    is a topological generator of an inertia group over $p_n$.
    Therefore Proposition~\ref{prop: general wreath conj crit} implies that
    \[
        c_1 \sim \s (1,\ldots,1,c_n^{\ve_n})
    \]
    in $[C_d]^\infty$ for some $d$-adic unit $\ve_n$.

    Since $c_\infty$ is the standard odometer, we have $\chi_\ell(c_\infty) = 1$ for all $\ell \geq 1$.
    Since $c_i =_\ell 1$ for all $\ell < i$, we have $\chi_\ell(c_i) = 0$ for $\ell < i$.
    The conjugation relations imply that $\chi_1(c_1) = 1$ and that $\chi_\ell(c_i) = \ve_{i-1}\chi_{\ell-1}(c_{i-1})$ for $\ell > 1$ with the $i-1$ subscripts interpreted modulo $n$.
    These recursive relations combined with the identity $c_\infty = c_n\cdots c_1$ imply that for $1 \leq \ell \leq n$,
    \[
        \chi_\ell(c_i) = 
        \begin{cases}
            1 & i = \ell\\
            0 & i \neq \ell.
        \end{cases}
    \]
    Thus for $1 \leq i < n$,
    \[
        1 = \chi_i(c_i) = \ve_{i-1}\chi_{i-1}(c_{i-1}) = \ve_{i-1}.
    \]
    If $1 < i \leq n$, then $\chi_{n+1}(c_i) = \chi_n(c_{i-1}) = 0$.
    Hence
    \[
        1 = \chi_{n+1}(c_\infty) = \chi_{n+1}(c_1) = \ve_n\chi_n(c_n) = \ve_n.
    \]
    Therefore $\ve_i \equiv 1 \bmod d$ for all $1 \leq i \leq n$.
    Lemma~\ref{lemma: d-adic unit conjugate} implies that $c_i^{\ve_i} \sim c_i$ in $[C_d]^\infty$, and we conclude that the following conjugacies hold in $[C_d]^\infty$,
    \[
        c_i \sim 
            \begin{cases}
                \s(1,\ldots,1,c_n) & \text{if }i = 1,\\
                (1,\ldots, 1, c_{i-1}) & \text{if }i \neq 1.
            \end{cases}
    \]

    (2)
    If $i \neq 1, m + 1$, then we argue as in the periodic case that there exists $d$-adic units $\ve_{i-1}$ such that
    \[
        c_i \sim (1,\ldots,1,c_{i-1}^{\ve_{i-1}})
    \]
    in $[C_d]^\infty$.
    If $i = 1$, then $c_1 =_1 \s$ and Lemma~\ref{lemma: inertia recursion} implies that $\rho_\Lambda(\gamma_{\Lambda,d})\cdots \rho_\Lambda(\gamma_{\Lambda,1})$ is a topological generator of an inertia group over 0.
    In the preperiodic case, 0 is not an element of $P_f$, hence
    \[
        \rho_\Lambda(\gamma_{\Lambda,d})\cdots \rho_\Lambda(\gamma_{\Lambda,1}) = 1.
    \]
    Therefore Proposition~\ref{prop: general wreath conj crit} implies that $c_1 \sim \s$ in $[C_d]^\infty$.
    If $i = m + 1$, then $c_{m+1} =_1 1$.
    Suppose $P$ is a prime in $K\sep(f^{-\infty}(t))$ over $(t - p_{m+1})$ such that $\gamma_{m+1}$ topologically generates the inertia group for $P$. 
    Let $q_j \in f^{-1}(p_{m+1})$ be such that $P \cap K\sep(t_j) = (t_j - q_j)$.
    Replacing $\gamma_{m+1}$ with a conjugate by a power of $\gamma_1$, we may suppose that $q_d = p_m$.
    Since $f$ is unicritical, we have $K\sep(f^{-1}(t)) = K\sep(t_j)$ and $t_j = \zeta_d^j t_d$ for each $j$.
    Thus for each $j$,
    \[
        (t_d - p_m) = (t_d - q_d) = (t_j - q_j) = (\zeta_d^jt_d - q_j) = (t_d - \zeta_d^{-j}q_j).
    \]
    Hence $q_j = \zeta_d^jp_m$.
    Therefore $q_\omega = \zeta_d^\omega p_m = p_n$.
    If $j \neq \omega, d$, then $q_j \notin P_f$.
    Thus there exists $d$-adic units $\ve_m$ and $\ve_n$ such that
    \[
        c_{m+1} \sim (1,\ldots,1,c_n^{\ve_n},1,\ldots,c_m^{\ve_m}),
    \]
    in $[C_d]^\infty$ where the $c_n^{\ve_n}$ is in the $\omega$th component.

    The conjugation identities imply that
    \[
        \chi_\ell(c_1) =   
        \begin{cases}
            1 & \ell = 1\\
            0 & \ell \neq 1,
        \end{cases}
    \]
    and $\chi_1(c_i) = 0$ for $i > 1$;
    if $i \neq 1, m + 1$, then $\chi_\ell(c_i) = \ve_{i-1}\chi_{\ell-1}(c_{i-1})$ for all $\ell > 1$; and $\chi_\ell(c_{m+1}) = \ve_m\chi_{\ell-1}(c_m) + \ve_n\chi_{\ell-1}(c_n)$.
    These recursive relations combined with $\chi_\ell(c_i) = 0$ for $\ell < i$ and $c_\infty = c_1\cdots c_n$ imply that for $1 \leq \ell \leq n$,
    \[
        \chi_\ell(c_i) = 
        \begin{cases}
            1 & i = \ell\\
            0 & i \neq \ell.
        \end{cases}
    \]
    Hence for $i \neq 1, m +1$,
    \[
        1 = \chi_i(c_i) = \ve_{i-1}\chi_{i-1}(c_{i-1}) = \ve_{i-1},
    \]
    and
    \[
        1 = \chi_{m+1}(c_{m+1}) = \ve_m\chi_m(c_m) + \ve_n\chi_m(c_n) = \ve_m.
    \]
    If $i \neq 1, m + 1$, then
    \[
        \chi_{n+1}(c_i) = \chi_n(c_{i-1}) = 0
    \]
    Hence, since \(\chi_{n+1}(c_1) = 0\), 
    \[
        1 = \chi_{n+1}(c_\infty) = \chi_{n+1}(c_{m+1}) = \chi_n(c_m) + \ve_n\chi_n(c_n) = \ve_n.
    \]
    Therefore $\ve_i \equiv 1 \bmod d$ for all $i$.

    Since the unique finite critical point is strictly preperiodic and $f$ is unicritical, it follows that each $c_i$ has order $d$ for $1 \leq i \leq n$.
    Thus
    \[
        c_i \sim 
        \begin{cases}
            \s & \text{if }i = 1\\
            (1,\ldots,1,c_n,1,\ldots,1,c_m) & \text{if }i = m+1, \text{where $c_n$ is in the $\omega$th component}\\
            (1,\ldots, 1, c_{i-1}) & \text{if }i \neq 1, m + 1.
        \end{cases}
        \qedhere
    \]
\end{proof}

\section{Model groups and semirigidity}
\label{section: model groups}

In Proposition~\ref{prop: geometric generators} we showed that for a unicritical PCF polynomial $f(x) \in K[x]$, the group $\barb f$ is topologically generated by a collection of elements in $[C_d]^\infty$ satisfying a system of \textit{conjugate} recurrences which is determined entirely by the combinatorial structure of the post-critical orbit in the periodic case, and in the preperiodic case further requires only a small piece of arithmetic information, the parameter $\omega$.
In this section we study two families of \emph{model groups} $\Mper(d,n)$ and $\Mpre(d,m,n,\omega) \subseteq [C_d]^\infty$ which are topologically generated by elements satisfying the same systems of recursion as $\barb f$ up to equality rather than conjugacy.
Ultimately we show that these generators have a certain \emph{semirigidity property} which allows us to deduce that $\Mper$ and $\Mpre$ are conjugate within $[C_d]^\infty$ to $\barb f$ in the periodic and preperiodic cases, respectively.

\begin{definition}
    Let $d \geq 2$.
    \begin{enumerate}
        \item Given an integer $n \geq 1$, let $\Mper = \Mper(d,n) \subseteq [C_d]^\infty$ be the closed subgroup defined by 
        \[
            \Mper := \lprof a_1, \ldots, a_n\rprof
        \]
        where
        \[
            a_i := 
            \begin{cases} 
                \s(1,\ldots,1,a_n) & \text{if }i = 1,\\ 
                (1\ldots,1,a_{i-1}) & \text{if }i\neq 1.
            \end{cases}
        \]
        Let $a_\infty := a_1a_2\cdots a_n \in \Mper$.

        \item Given integers $n > m \geq 0$ and $1 \leq \omega < d$, let $\Mpre = \Mpre(d,m,n,\omega) \subseteq [C_d]^\infty$ be the closed subgroup defined by
        \[
            \Mpre := \lprof b_1, \ldots, b_n\rprof
        \]
        where
        \[
            b_i :=
            \begin{cases}
                \s & \text{if }i = 1\\
                (1,\ldots,1,b_n,1,\ldots,1,b_m) & \text{if }i = m+1, \text{where $b_n$ is in the $\omega$th component}\\
                (1,\ldots, 1, b_{i-1}) & \text{if }i \neq 1, m + 1.
            \end{cases}
        \]
        Let $b_\infty := b_1b_2\cdots b_n \in \Mpre$.
    \end{enumerate}
\end{definition}

\begin{notation}
    Given an integer $1 \leq i \leq n$ we let $\wh{\Mper}_i$ denote the normal subgroup of $\Mper$ topologically generated by all conjugates of $a_j$ with $j \neq i$, and  define $\wh{\Mpre}_i$ analogously.
    
    Similarly, if $1 \leq i, j \leq n$, then  $\wh{\Mpre}_{i,j}$ denotes the normal subgroup of $\Mpre$ topologically generated by all conjugates of $b_k$ with $k\neq i, j$.
    In $\Mpre$, we write $b_{i,j} := \s^{j} b_i \s^{-j}$.
    Since $\s = b_1$, it follows that $b_{i,j} \in \Mpre$.
\end{notation}

The recursive descriptions of the model groups imply that $\st_1\Mper \subseteq \Mper^d$ and $\st_1\Mpre \subseteq \Mpre^d$, which is to say that $\Mper$ and $\Mpre$ are both self-similar.
Furthermore, each coordinate projection is surjective.

\begin{lemma}
\label{lemma: stabilizer projections are surjective}
    Both $\Mper$ and $\Mpre$ are self-similar. Moreover, the coordinate projections $\pi_j : \st_1\Mper \to \Mper$ and $\pi_j : \st_1 \Mpre \to \Mpre$ are surjective for all $1\leq j \leq d$.
\end{lemma}

\begin{proof}
    First consider $\Mper$.
    If $2 \leq i \leq d$, then conjugating $a_i$ by powers of $a_1$ gives us
    \[
        (1,\ldots,1,a_{i-1},1,\ldots,1) \in \st_1\Mper
    \]
    where the non-trivial component can be in any coordinate.
    We also have
    \[
        a_1^d = (a_n,\ldots,a_n) \in \st_1\Mper.
    \]
    Therefore $\pi_j(\st_1\Mper)$ contains the group topologically generated by all of the $a_i$, which is $\Mper$.
    Hence each $\pi_j : \st_1\Mper \to \Mper$ is surjective.

    Next consider $\Mpre$.
    Observe that conjugating $b_i$ with $i \neq 1, m + 1$ by powers of $b_1 = \s$ gives
    \[
        (1,\ldots,1,b_{i-1},1,\ldots,1) \in \st_1\Mpre
    \]
    with the non-trivial component in any coordinate.
    Conjugating $b_{m+1}$ by powers of $b_1$ gives elements of $\st_1\Mpre$ with either $b_m$ or $b_n$ in any prescribed coordinate.
    Thus $\pi_j(\st_1\Mpre)$ contains the group topologically generated by the $b_i$, namely $\Mpre$.
    Hence each $\pi_j : \st_1\Mpre \to \Mpre$ is surjective.
\end{proof}

\subsection{Orders of generators}

The defining systems of recurrence for the model groups along with the recursive formula for the characters $\chi_\ell$ (Definition~\ref{definition: characters}) allow us to calculate the values of $\chi_\ell$ at each generator.
We make frequent use of these calculations, including in the determination of the orders of the generators.

\begin{lemma}
\label{lemma: chi of generators}
Let $1 \leq i \leq n$ and $\ell \geq 1$, then
\begin{enumerate}
    \item $\displaystyle{\chi_\ell(a_i) =
    \begin{cases}
        1 & \text{if }\ell \equiv i \bmod n,\\
        0 & \text{otherwise.}
    \end{cases}
    }$

    \item $\displaystyle{\chi_\ell(b_i) =
    \begin{cases}
        1 & \ell = i \leq m\\
        1 & \ell \geq i > m \text{ and } \ell \equiv i \bmod n - m\\
        0 & \text{otherwise.}
    \end{cases}
    }$
\end{enumerate}
\end{lemma}

\begin{proof}
    (1)
    If $\ell = 1$, then by definition of the $a_i$ we have
    \[
        \chi_1(a_i) =
        \begin{cases}
            1 & \text{if } i = 1,\\
            0 & \text{otherwise}.
        \end{cases}
    \]
    If $\ell > 1$, then the definition of $a_i$ and the recursive formula for $\chi_\ell$ imply that
    \[
        \chi_\ell(a_i) = \chi_{\ell-1}(a_{i-1})
    \]
    where we interpret the subscript of $a_{i-1}$ modulo $n$.
    Hence we may conclude by induction that
    \[
        \chi_{\ell}(a_i) = \chi_{1}(a_{i - \ell + 1})
        =
        \begin{cases}
            1 & \text{if } \ell \equiv i \bmod n,\\
            0 & \text{otherwise}.
        \end{cases}
    \]

    (2)
    We proceed similarly for the $b_i$.
    Starting with $\ell = 1$, the definition of the $b_i$ imply that
    \[
        \chi_1(b_i) =
        \begin{cases}
            1 & \text{if } i = 1\\
            0 & \text{otherwise}.
        \end{cases}
    \]
    If $\ell > 1$, then
    \[
        \chi_\ell(b_i) = 
        \begin{cases}
            \chi_{\ell-1}(b_{i-1}) & \text{if }i \neq 1, m + 1,\\
            0 & \text{if }i = 1,\\
            \chi_{\ell-1}(b_m) + \chi_{\ell-1}(b_n) & \text{if }i = m + 1.
        \end{cases}
    \]
    These relations imply that $\chi_\ell(b_i) = 0$ whenever $\ell < i$.
    Furthermore, if $i \leq m$, then
    \begin{equation}
    \label{eqn chi small i}
        \chi_\ell(b_i) = 
        \begin{cases}
            1 & \text{if }\ell = i,\\
            0 & \text{otherwise.}
        \end{cases}
    \end{equation}
    Next we calculate $\chi_\ell(b_n)$ for all $\ell \geq 1$.
    If $\ell < n - m$, then $n - \ell + 1 > m + 1$ and
    \[
        \chi_\ell(b_n) = \chi_1(b_{n-\ell+1}) = 0.
    \]
    If $\ell \geq n - m$, then
    \[
        \chi_\ell(b_n) = \chi_{\ell - (n-m)}(b_m) + \chi_{\ell-(n-m)}(b_n).
    \]
    Hence if $\ell = r + q(n-m)$ with $0 \leq r < n-m$, then a simple induction implies that
    \[
        \chi_\ell(b_n) = \chi_r(b_n) + \sum_{s=0}^{q-1}\chi_{r + s(n-m)}(b_m)
        =\sum_{s=0}^{q-1}\chi_{r + s(n-m)}(b_m).
    \]
    Observe that \eqref{eqn chi small i} implies that $\chi_{r + s(n-m)}(b_m) = 1$ if and only if $r + s(n-m) = m$ for some $0 \leq s < q$, which is equivalent to $\ell > m$ and $\ell \equiv m \bmod n - m$.
    Therefore
    \[
        \chi_\ell(b_n) =
        \begin{cases}
            1 & \text{if } \ell > m \text{ and } \ell\equiv m \bmod n-m,\\
            0 & \text{otherwise.}
        \end{cases}
    \]
    Now suppose that $m < i \leq n$ and $\ell \geq i$.
    Then
    \[
        \chi_\ell(b_i) = \chi_{\ell - i + m}(b_m) + \chi_{\ell-i+m}(b_n).
    \]
    Our calculations imply that
    \[
        \chi_{\ell-i+m}(b_m) =
        \begin{cases}
            1 &\text{if }\ell = i,\\
            0 &\text{otherwise,}
        \end{cases}
    \]
    and
    \[
        \chi_{\ell-i+m}(b_n) =
        \begin{cases}
            1 &\text{if }\ell > i\text{ and }\ell\equiv i \bmod n - m,\\
            0 &\text{otherwise.}
        \end{cases}
    \]
    Therefore if $m < i \leq n$ and $\ell \geq i$, then
    \[
        \chi_\ell(b_i) =
        \begin{cases}
            1 &\text{if }\ell \equiv i \bmod n - m,\\
            0 &\text{otherwise.}
        \end{cases}
    \]
    
    Putting this all together we have
    \[
        \chi_\ell(b_i) = 
        \begin{cases}
            1 & \ell = i \leq m.\\
            1 & \ell \geq i > m \text{ and } \ell \equiv i \bmod n - m,\\
            0 & \text{otherwise.}
        \end{cases}
        \qedhere
    \]
\end{proof}

Recall that $a_\infty = a_1\cdots a_n$ and $b_\infty = b_1\cdots b_n$.
Thus Lemma~\ref{lemma: chi of generators} implies that
\[
    \chi_\ell(a_\infty) = \chi_\ell(b_\infty) = 1
\]
for all $\ell \geq 1$.
Hence Lemma~\ref{lemma: chi odometer} implies that $a_\infty$ and $b_\infty$ are both strict odometers.

\begin{proposition}
\label{prop: order of generators}
    Let $1 \leq i \leq n$ and $\ell \geq 0$, then
    \begin{enumerate}
        \item 
        $
            \ord_\ell(a_i) = d^{\big\lfloor \frac{\ell - i}{n}\big\rfloor + 1}.
        $

        \item $\displaystyle{\ord_\ell(b_i) = 
        \begin{cases}
            d & \text{if }\ell \geq i,\\
            1 & \text{otherwise.}
        \end{cases}
        }$
    \end{enumerate}
    In particular, the map $\ZZ_d \to \lprof a_i \rprof$ defined by $m \mapsto a_i^m$ is an isomorphism and each $b_i$ has order $d$.
\end{proposition}

\begin{proof}
    (1)
    It suffices to prove that for all $1 \leq i \leq n$ and all $\ell\geq 0$,
    \begin{equation}
    \label{eqn order}
        \log_d(\ord_\ell(a_i)) = \Big\lfloor \frac{\ell - i}{n}\Big\rfloor + 1.
    \end{equation}
    Note that $a_i =_0 1$ implies $\log_d(\ord_0(a_i)) = 0$ for all $1 \leq i \leq n$.
    If $1 < i \leq n$, then $a_i = (1,\ldots, 1, a_{i-1})$ implies that
    \[
        \log_d(\ord_\ell(a_i)) = \log_d(\ord_{\ell-1}(a_{i-1})).
    \]
    Since $a_1 = \s(1,\ldots, 1, a_n)$ and $\s$ has order $d$, we have
    \begin{align*}
        \log_d(\ord_\ell(a_1)) &= \log_d(\ord_\ell(a_1^d)) + 1\\
        &= \log_d(\ord_\ell((a_n,\ldots,a_n))) + 1\\
        &= \log_d(\ord_{\ell-1}(a_n)) + 1.
    \end{align*}
    Thus for all $\ell> 1$,
    \[
        \log_d(\ord_\ell(a_i)) = \begin{cases} \log_d(\ord_{\ell-1}(a_n)) + 1 & \text{if }i = 1,\\ \log_d(\ord_{\ell-1}(a_{i-1})) & \text{otherwise}.\end{cases}
    \]
    Furthermore, these recursive identities completely determine $\ord_\ell(a_i)$ for all $\ell\geq 0$ and all $1 \leq i \leq n$.
    On the other hand, let 
    \[
        \alpha_{\ell,i} := \Big\lfloor \frac{\ell - i}{n}\Big\rfloor + 1.
    \]
    We will show that $\alpha_{\ell,i}$ satisfies the same recursive identities and initial values.
    First note that for $1 \leq i \leq n$,
    \[
        \alpha_{0,1} = \Big\lfloor \frac{-i}{n}\Big\rfloor + 1 = 0.
    \]
    If $2 \leq i \leq n$, then
    \[
        \alpha_{\ell,i} = \Big\lfloor \frac{\ell - i}{n}\Big\rfloor + 1
        = \Big\lfloor \frac{(\ell-1) - (i-1)}{n}\Big\rfloor + 1
        =\alpha_{\ell-1,i-1}.
    \]
    Finally,
    \[
        \alpha_{\ell,1} = \Big\lfloor \frac{\ell - 1}{n}\Big\rfloor + 1
        = \Big(\Big\lfloor \frac{(\ell - 1) - n}{n}\Big\rfloor + 1\Big) + 1
        = \alpha_{\ell-1,n} + 1.
    \]
    Thus \eqref{eqn order} holds for all $\ell\geq 0$ and $1 \leq i \leq n$.
    Therefore
    \[
        \log_d(\ord_\ell(a_i)) = \alpha_{\ell,i} = \Big\lfloor \frac{\ell-i}{n}\Big\rfloor + 1.
    \]

    (2)
    Let $1\leq i \leq n$. 
    If $\ell < i$, then the recursive formulas for $b_i$ tell us that $b_i$ acts trivially on $T_d^\ell$, hence that $\ord_\ell(b_i) = 1$.
    Lemma~\ref{lemma: chi of generators} implies that $\chi_i(b_i) = 1 \in \ZZ/d\ZZ$ for all $i$.
    Hence the order of each $b_i$ is a multiple of $d$.
    On the other hand, note that the elements $b_i^d$ satisfy the system of recursions
    \begin{align*}
        x_1 &= 1\\
        x_{m+1} &= (1,\ldots,1,x_n,1,\ldots,1,x_m)\\
        x_i &= (1,\ldots,1, x_{i-1}) \text{ if }i\neq 1, m+1,
    \end{align*}
    which are also satisfied by $x_i = 1$.
    Hence Lemma~\ref{lemma: recursive solution} implies that $b_i^d = 1$ for each $1\leq i \leq n$.
    Therefore each $b_i$ has order exactly $d$.
    Hence $\ord_\ell(b_i) = d$ for all $\ell \geq i$.
\end{proof}

\subsection{Branching}
\label{subsection: branching}

A closed subgroup $G \subseteq [S_d]^\infty$ is called a \emph{weakly branch group} if $G$ acts transitively on every level of the tree $T_d^\infty$ and $G$ contains an infinite closed normal subgroup $K$ such that $K^d \subseteq G$, in which case we say that $G$ is \emph{weakly branched over $K$}.
If $K$ has finite index in $G$, then we say that $G$ is a \emph{regular branch group} and that $G$ \emph{is branched over $K$}.
Weakly branch and regular branch groups are an important, natural class of just infinite groups; see Grigorchuk \cite{grigorchuk2000just}, Bartholdi, Grigorchuk, and \v{S}uni\'{k} \cite{bartholdi2003branch}, or Nekrashevych \cite{nekrashevych2024self}.

We show $\Mper$ is weakly branch in Lemma \ref{lemma: A subgroup properties}, and in in Proposition~\ref{prop: branch} we prove that $\Mpre$ is branch for all but one choice of the defining parameters.
This has interesting arithmetic implications; for example, in Proposition~\ref{prop: preper finite const} we show that the constant field extension of $K(f^{-\infty}(t))/K(t)$ is finite whenever $\Mpre$ is branch.

\begin{lemma}
\label{lemma: A subgroup properties}
    The group $\Mper$ is weakly branch over $\wh\Mper_n$.
    Furthermore, $\wh{\Mper}_i \cap \lprof a_i\rprof = 1$ for $1 \leq i \leq n$ and $\st_1\Mper = \lprof a_1^d\rprof\wh{\Mper}_n^d$.
\end{lemma}

\begin{proof}
The group $\wh{\Mper}_1$ is topologically generated by all conjugates of $a_i$ with $2 \leq i \leq n$.
If $2 \leq i \leq n$, then every $\Mper$-conjugate of $a_i = (1, \ldots, 1, a_{i-1})$ belongs to $\wh{\Mper}_n^d$, hence $\wh{\Mper}_1 \subseteq \wh{\Mper}_n^d$.
On the other hand, $\wh{\Mper}_n^d$ is generated by all the $\Mper$-conjugates of elements of the form $(1,\ldots,1,a_{i-1},1,\ldots, 1)$ with $2\leq i \leq n$ where the $a_{i-1}$ can be in any component.
These elements are conjugates of $a_i$ by powers of $a_1$, hence belong to $\wh{\Mper}_1$.
Therefore $\wh{\Mper}_1 = \wh{\Mper}_n^d$.
Hence $\Mper$ is weakly branch over $\wh\Mper_n$.

By Proposition~\ref{prop: inductive criteria}, it suffices to show that $\wh{\Mper}_i \cap \lprof a_i \rprof =_\ell 1$ for each $1 \leq i \leq n$ and for all $\ell \geq 0$.
We proceed by induction on $\ell$.
The $\ell = 0$ case is trivial, so suppose that $\ell \geq 1$ and that $\wh{\Mper}_i \cap \lprof a_i \rprof =_{\ell-1} 1$ for each $1\leq i \leq n$.
Then $\wh{\Mper}_n \cap \lprof a_n \rprof =_{\ell-1} 1$ and
\[
    \st_1\Mper \cap \lprof a_1 \rprof = \lprof a_1^d \rprof = \lprof (a_n,\ldots, a_n) \rprof
\]
imply that
\[
    \wh{\Mper}_1 \cap \lprof a_1 \rprof  \subseteq \wh{\Mper}_n^d \cap \lprof (a_n,\ldots, a_n) \rprof \subseteq (\wh{\Mper}_n \cap \lprof a_n \rprof)^d =_\ell 1.
\]
Next suppose that $2 \leq i \leq n$.
Observe that $a_i = (1, \ldots, 1,a_{i-1}) \in 1^{d-1} \times \lprof a_{i-1}\rprof$ and
\[
    \wh{\Mper}_i \cap (1^{d-1}\times \Mper)\subseteq 1^{d-1} \times \wh{\Mper}_{i-1}.
\]
Thus
\begin{align*}
    \wh{\Mper}_i \cap \lprof a_i \rprof &\subseteq
    \wh{\Mper}_i \cap (1^{d-1} \times \Mper) \cap \lprof a_i \rprof\\
    &\subseteq (1^{d-1} \times \wh{\Mper}_{i-1}) \cap \lprof (1,\ldots,1,a_{i-1})\rprof\\
    &\subseteq 1^{d-1}\times (\wh{\Mper}_{i-1}\cap \lprof a_{i-1}\rprof)\\
    &=_\ell 1.
\end{align*}
This completes our induction.

The definitions of $\Mper$ and $\wh{\Mper}_1$ imply that
$\Mper = \lprof a_1 \rprof\wh{\Mper}_1$.
Then by $\wh{\Mper}_1 \subseteq \st_1\Mper$, $\st_1\Mper \cap \lprof a_1\rprof = \lprof a_1^d\rprof$, and (1) we have
\[
    \st_1\Mper = \lprof a_1^d\rprof\wh{\Mper}_1 = \lprof a_1^d\rprof\wh{\Mper}_n^d.\qedhere
\]
\end{proof}

Lemma \ref{lemma: A subgroup properties} allows us to construct a useful family of abelian quotients of $\Mper$.

\begin{proposition}
\label{prop: eta construction}
    For each $1 \leq i \leq n$, there is a continuous surjection $\eta_i: \Mper \to \ZZ_d$ such that $\ker(\eta_i) = \wh{\Mper}_i$ and
    \[
        \eta_i(a_j) = 
        \begin{cases} 
            1 & \text{if }i = j,\\ 
            0 & \text{otherwise.}
        \end{cases}
    \]
\end{proposition}

\begin{proof}
    The quotient $\Mper/\wh{\Mper}_i$ is topologically generated by the image of $a_i$, hence Lemma~\ref{lemma: A subgroup properties}(2) and Proposition~\ref{prop: order of generators} together imply that $\Mper/\wh{\Mper}_i \cong \lprof a_i \rprof \cong \ZZ_d$.
    We define $\eta_i : \Mper \to \ZZ_d$ as the composition
    \[
        \eta_i: \Mper \longrightarrow \Mper/\wh{\Mper}_i \cong \ZZ_d.
    \]
    Therefore $\ker(\eta_i) = \wh{\Mper}_i$.
    Since $a_j \in \wh{\Mper}_i$ for $i \neq j$ and the isomorphism $\lprof a_i \rprof \cong \ZZ_d$ maps $a_i \mapsto 1$ by construction, we conclude that
    \[
        \eta_i(a_j) = 
        \begin{cases} 
            1 & \text{if }i = j,\\ 
            0 & \text{otherwise.}
        \end{cases}
        \qedhere
    \]
\end{proof}

The maps $\eta_i$ and $\chi_i$ combine to give all the abelian quotients of $\Mper$ and $\Mpre$, respectively.

\begin{proposition}[Abelianizations]\,
\label{prop: abelianization}
    \begin{enumerate}
        \item The map $\eta : \Mper \to \ZZ_d^n$ defined by
        \[
            \eta(a) := (\eta_1(a), \ldots, \eta_n(a))
        \]
        is a surjective homomorphism which induces an isomorphism $\Mper^{\ab} \cong \ZZ_d^n$.

        \item The map $\chi : \Mpre \to (\ZZ/d\ZZ)^n$ defined by
        \[
            \chi(b) := (\chi_1(b), \ldots, \chi_n(b))
        \]
        is a surjective homomorphism which induces an isomorphism $\Mpre^\ab \cong (\ZZ/d\ZZ)^n$.
    \end{enumerate}
\end{proposition}

\begin{proof}
    (1)
    Proposition~\ref{prop: eta construction} implies that the generators $a_i$ map under $\eta$ to the standard basis in $\ZZ_d^n$, hence $\eta$ is a surjective homomorphism
    Since $\ZZ_d^n$ is abelian, $\eta$ factors through the abelianization, and since the $a_i$ topologically generate $\Mper$ and $\lprof a_i \rprof \cong \ZZ_d$, the abelianization factors through $\ZZ_d^n$.
    Hence $\eta$ induces an isomorphism $\Mper^\ab \cong \ZZ_d^n$.

    (2)
    If $1 \leq i \leq n$, then Lemma~\ref{lemma: chi of generators} implies that the generators $b_i$ map under $\chi$ to the standard basis in $(\ZZ/d\ZZ)^n$.
    Thus $\chi$ is a surjective homomorphism.
    On the other hand, since $\Mpre$ is generated by the $b_i$ which each have order $d$ by Proposition~\ref{prop: order of generators}(2), it follows that $\Mpre^\ab$ has order at most $d^n$, hence $\chi$ induces an isomorphism $\Mpre^\ab \cong (\ZZ/d\ZZ)^n$.
\end{proof}

Next we turn to the group $\Mpre$.
Our analysis naturally splits along several cases.
Consider the hypotheses,

\begin{enumerate}[leftmargin=5em]
    \item[$(A_1)$] $\omega \neq d/2$,
    \item[$(A_2)$] $m > 1$ and $(d,n) \neq (2,m+1)$,
    \item[$(A_3)$]  $d = 2$, $m > 2$, and $n = m + 1$,
    \item[$(B_1)$] $d > 2$, $m = 1$, and $\omega = d/2$,
    \item[$(B_2)$] $d = 2$, $m = 1$, and $n > 2$,
    \item[$(C)$] $d = 2$, $m = 2$, and $n = 3$.
    \item[$(D)$] $d = 2$, $m = 1$, and $n = 2$.
\end{enumerate}

We let $(A)$ refer to the assumption $(A_1)$,$(A_2)$, or $(A_3)$, and let $(B)$ refer to $(B_1)$ or $(B_2)$.
Note that these hypotheses exhaust all possible cases with $d \geq 2$, $1 \leq m < n$, and $1 \leq \omega < d$.
When $d = 2$, we have $\omega = d/2 = 1$ by default.
Case $(D)$ is exceptional; it corresponds, up to conjugacy, to the Chebyshev polynomial of degree $2$. 

According to these cases, we define a subgroup $\br$ over which $\Mpre$ is regular branch in all cases except $(D)$. We start by observing that a natural subgroup branches.
\begin{lemma}
\label{lemma: pre-branch}
    The subgroup $\wh \Mpre_{m,n}$ branches, meaning \((\wh{\Mpre}_{m,n})^d = \wh{\Mpre}_{1,m+1} \subseteq \Mpre\).
\end{lemma}

\begin{proof}
    The group $\wh{\Mpre}_{1,m+1}$ is generated by all conjugates of $b_i = (1,\ldots, 1, b_{i-1})$ with $i \neq 1, m + 1$.
    Since $\s = b_1 \in \Mpre$, it follows that $(\wh{\Mpre}_{m,n})^d = \wh{\Mpre}_{1,m+1} \subseteq \Mpre$.
\end{proof}

However, the subgroup $\wh\Mpre_{m,n}$ is not always finite index, nor is it the maximal normal subgroup with this property.
We define $\br$ as an extension of $\wh\Mpre_{m,n}$ by cases.

\begin{definition}
    Let $\br \subseteq \Mpre$ be the closed normal subgroup defined in cases as follows:
    \begin{enumerate}
        \item If $(A)$, then $\br$ is the closed normal subgroup generated by $[b_m, b_n]$ and $\wh{\Mpre}_{m,n}$.
        \item If $(B)$ or $(C)$, then $\br$ is the closed normal subgroup generated by $\wh{\Mpre}_{m,n}$, $[b_m, [b_m,b_n]]$, and $[b_n, [b_m,b_n]]$.
        \item If $(D)$, then $\br$ is the trivial subgroup.
    \end{enumerate}
\end{definition}

We determine the quotients $\Mpre/\br$ in Lemma~\ref{lemma: N quotients}; in particular showing that $\br$ has finite index.
In cases $(A)$ and $(B)$, the quotient $\Mpre/\br$ is induced by a quotient of $[C_d]^\infty$.
The failure of this in case $(C)$ is what distinguishes it from $(B)$.
We show that $\br^d \subseteq \Mpre$ in Proposition~\ref{prop: branch}.
    
\begin{definition}
    Let $\psi_A : [C_d]^\infty \to (\ZZ/d\ZZ)^2$ be the map defined by
    \[
        \psi_A(g) = (\chi_m(g), \chi_n(g)).
    \]
    If $g = \s^i(g_1,\ldots,g_d) \in [C_d]^\infty$, then let $\chi' : [C_d]^\infty \to \ZZ/d\ZZ$ be the function defined by
    \[
        \chi'(g) := \sum_{i=1}^d i\chi_{n-1}(g_i) \in \ZZ/d\ZZ.
    \]
    Recall that $\H_d$ is the $d$th Heisenberg group (see Section \ref{section: heisenberg}).
    Let $\psi_B: [C_d]^\infty \to \H_d$ be the map defined by
    \[
        \psi_B(g) := \s^{-\chi_1(g)}(\chi'(g),\chi_n(g)).
    \]
\end{definition}

\begin{lemma}
\label{lemma: psi properties}
    The maps $\psi_A$ and $\psi_B$ are surjective homomorphisms, and
    \[
        \br = 
        \begin{cases}
            \ker \psi_A \cap \Mpre & \textrm{ if } (A),\\
            \ker \psi_B \cap \Mpre & \textrm{ if } (B).
        \end{cases}
    \]
\end{lemma}

\begin{proof}
    First suppose $(A)$.
    As the direct product of homomorphisms, $\psi_A$ is clearly a homomorphism.
    The elements $b_m$ and $b_n$ are mapped by $\psi_A$ to generators of $(\ZZ/d\ZZ)^2$, while the other $b_i$ are in the kernel, hence $\wh\Mpre_{m,n}\subseteq \br$. 
    Since $(\ZZ/d\ZZ)^2$ is abelian, it follows that $[b_m,b_n] \in \ker \psi_A$, and therefore $\br \subseteq \ker\psi_A$.
    On the other hand, $\br$ clearly has index at most $d^2 = [\Mpre: \ker\psi_A \cap \Mpre]$, thus $\br = \ker\psi_A \cap \Mpre$.

    Next suppose $(B)$.
    Although $\chi'$ is not a homomorphism on $[C_d]^\infty$, it restricts to one on $\st_1[C_d]^\infty$.
    Therefore $\psi_B$ restricts to a homomorphism on $C_d = \langle \s \rangle$ and on $\st_1 [C_d]^\infty$.
    Let $g := (g_1,\ldots,g_d) \in \st_1[C_d]^\infty$.
    Since $[C_d]^\infty  = C_d \ltimes \st_1[C_d]^\infty$, to show that $\psi_B$ is a homomorphism on $[C_d]^\infty$, it suffices to show that
    \[
        \psi_B(\s)\inv \psi_B(g)\psi_B(\s) = \psi_B(\s\inv g\s).
    \]
    The left hand side simplifies to
    \[
        \psi_B(\s)\inv \psi_B(g)\psi_B(\s) = \s(\chi'(g),\chi_n(g))\s\inv
        = (\chi'(g) - \chi_n(g), \chi_n(g)).
    \]
    While the right hand side is
    \[
        \psi_B(\s\inv g\s) = \psi_B(g_2,\ldots,g_d,g_1)
        = \Big(\sum_{i=1}^d i\chi_{n-1}(g_{i+1}), \chi_n(\s\inv g \s)\Big)
        = \Big(\sum_{i=1}^d i\chi_{n-1}(g_{i+1}), \chi_n(g)\Big),
    \]
    where
    \begin{align*}
        \sum_{i=1}^d i \chi_{n-1}(g_{i+1})
        &= \sum_{i=1}^d (i+1) \chi_{n-1}(g_{i+1}) - \chi_{n-1}(g_{i+1})\\
        &= \Big(\sum_{i=1}^d (i+1) \chi_{n-1}(g_{i+1})\Big) - \Big(\sum_{i=1}^d\chi_{n-1}(g_{i+1})\Big)\\
        &= \chi'(g) - \chi_n(g).
    \end{align*}
    Therefore $\psi_B : [C_d]^\infty \to \H_d$ is a homomorphism.

    In case $(B)$, we have $m = 1$, so $b_m = b_1 = \sigma$, and $b_n = (1, \ldots, 1, b_{n-1})$. 
    The definition of $\psi_B$ implies that $\psi_B(b_n) = (0,1)$, which, in combination with $\psi_B(\s\inv) = \s$, generates $\H_d$, so $\psi_B$ is surjective.
    The definition of $\psi_B$ implies that $b_i \in \ker\psi_B$ for all $i \neq m, n$, hence $\wh\Mpre_{m,n} \subseteq \ker\psi_B$.
    Observe that
    \[
        \psi_B([b_m,b_n]) = [\s\inv,(0,1)] = (-1,0),
    \]
    which commutes with $\psi_B(b_m) = \s\inv$ and $\psi_B(b_n) = (0,1)$.
    Therefore $[b_m,[b_m,b_n]]$ and $[b_n,[b_m,b_n]]$ both belong to $\ker\psi_B$ as well.
    Therefore $\br \subseteq \ker\psi_B \cap \Mpre$, which is to say there is a surjective homomorphism $\Mpre/\br \to \Mpre/(\ker\psi_B \cap \Mpre) \cong \H_d$.
    
    On the other hand, the definition of $\br$ and Lemma~\ref{lemma: H_d presentation} imply there is a surjective homomorphism $\H_d \to \Mpre/\br$.
    Hence $\Mpre/\br \cong \H_d$ and $\br = \ker\psi_B \cap \Mpre$.
\end{proof}

\begin{remark}
    We use this description of $\br$ as the intersection of a normal subgroup of $[C_d]^\infty$ with $\Mpre$ provided by Lemma~\ref{lemma: psi properties} when analyzing the normalizer of $\Mpre$ in Section~\ref{section: normalizer of B}.
    However, this does not hold in case $(C)$:
    if $\br'$ is the normal closure of $\br$ in $[C_2]^\infty$, then a calculation shows that $\br$ has index 2 in $\br' \cap \Mpre$, whereas \(\br\) has index \(8\) in \(\Mpre\) (see Appendix \ref{calculation: gap code}).
\end{remark}

With these, we are now able to show that $\br$ has finite index in $\Mpre$,
\begin{lemma}
The quotients $\Mpre/\br$ are as follows:
\label{lemma: N quotients}
    \begin{enumerate}
        \item If $(A)$, then $\Mpre/\br \cong (\ZZ/d\ZZ)^2$,
        \item If $(B)$ or $(C)$, then $\Mpre/\br \cong \H_d$,
        \item If $(D)$, then $\Mpre \cong \Mpre/\N \cong \wh{D}_\infty$, the pro-2 dihedral group.
    \end{enumerate}    
\end{lemma}

\begin{proof}
    The $(A)$ and $(B)$ cases follow immediately from Lemma~\ref{lemma: psi properties}.
    In case $(C)$, the definition of $\br$ and Lemma~\ref{lemma: H_d presentation} imply there is a surjective homomorphism $\H_2 \to \Mpre/\br$.
    On the other hand, $[\Mpre : \br] \geq [\Mpre : \br]_4$ and a computer calculation (see Appendix~\ref{calculation: gap code})
    shows that $[\Mpre : \br]_4 = 8 = |\H_2|$.
    Thus this surjection is an isomorphism.

    Finally suppose $(D)$, in which case $n=2$ so $\Mpre = \lprof b_1, b_2\rprof$ and $\br$ is trivial. Both $b_1$ and $b_2$ have order 2 by Proposition~\ref{prop: order of generators}, while the identity $(b_1b_2)^2 = (b_1b_2, b_2b_1)$ implies that $b_1b_2$ is an odometer, hence \(\lprof b_1b_2\rprof \cong \ZZ_2\).
    Therefore $\Mpre \cong \Mpre/\br$ is isomorphic to the pro-2 dihedral group $\wh{D}_\infty$.
\end{proof}

\begin{proposition}
\label{prop: branch}
In cases $(A)$, $(B)$, and $(C)$, $\Mpre$ is branch over $\br$.

\end{proposition}
\begin{proof}
    We proved in Lemma \ref{lemma: N quotients} that $\br$ has finite index in cases $(A)$, $(B)$, and $(C)$.
    Hence it suffices to prove that $\br^d \subseteq \Mpre$.
    
    Recall that $b_{i,j} := \s^j b_i \s^{-j}$.
    If $(A_1)$, then $\omega \neq d/2$ which implies $g := [b_{m+1,\omega},b_{m+1}] \in \Mpre$ is supported in only the $\omega$th component.
    Since $\pi_\omega(g) = [b_m,b_n]$ and $\br= [b_m,b_n]\wh{\Mpre}_{m,n}$, we conclude that $\br^d \subseteq \Mpre$.
    
    If $(A_2)$, then $m > 1$ and $\s = b_1 \in \wh{\Mpre}_{m,n}$.
    Our assumption that $(d,n) \neq (2,m+1)$ implies that there exists some $j$ such that $b_{m,j}$ commutes with $b_n$.
    Therefore
    \[
        [b_m,b_n] \equiv [b_{m,j},b_n] \equiv1 \bmod \wh{\Mpre}_{m,n},
    \]
    or equivalently $[b_m,b_n] \in \wh{\Mpre}_{m,n}$.
    This implies that $\br= \wh{\Mpre}_{m,n}$ in case $(A_2)$.
    Therefore (1) implies that $\br^d \subseteq \Mpre$.

    If $(A_3)$, then $b_2 = (1,\s) \in \wh{\Mpre}_{m,n}$ and both $b_{m-1}$ and $b_{n-1}$ are supported in exactly one coordinate.
    Thus
    \[
        [b_m, b_n] = (1, [b_{m-1},b_{n-1}]) \equiv (1,[b_{m-1,1},b_{n-1}]) \equiv 1 \bmod \wh{\Mpre}_{m,n}.
    \]
    Hence in this case we also have $\br= \wh{\Mpre}_{m,n}$ and therefore $\br^d \subseteq \Mpre$.

    Next suppose $(B_1)$: $d > 2$, $\omega = d/2$, and $m = 1$.
    Observe that
    \[
        [b_n, [b_m, b_n]] = [b_n, b_{n,-1}\inv b_n] = [b_n,b_{n,-1}\inv]^{b_n}.
    \]
    Since $\omega \neq d - 1$, we have $[b_n,b_{n,-1}\inv] = 1$.
    Therefore, $[b_n,[b_m,b_n]] = 1$.
    Consider the element $g = [b_{m+1,\omega},[b_{m+1,\omega},b_{m+1}]] \in \Mpre$.
    Note that $g$ is supported in at most the $\omega$th and $d$th coordinates.
    We have
    $
        \pi_\omega(g) = [b_m,[b_m,b_n]],
    $
    and
    $
        \pi_d(g) = [b_n,[b_m,b_n]] = 1.
    $
    Therefore $g$ is an element of $\Mpre$ supported in a single coordinate with support $[b_m,[b_m,b_n]]$.
    Since $\br$ is the normal subgroup generated by $\wh{\Mpre}_{m,n}$, $[b_m,[b_m,b_n]]$ and $[b_n,[b_m,b_n]]$, (1) implies that $\br^d \subseteq \Mpre$.

    If $(B_2)$, then $d = 2$, $m = 1$, and $n > 2$.
    Thus $b_n = (1,b_{n-1})$ and
    \[
        [b_m,b_n] = [\s, (1,b_{n-1})] = (b_{n-1},b_{n-1}),
    \]
    which clearly commutes with both $\s$ and $b_n$.
    Hence $\brB = \wh{\Mpre}_{m,n}$.
    
    Similarly in case $(C)$ where $d = 2$, $m = 2$, and $n = 3$, we have
    \[
        [b_2,b_3] = (1,[\s,b_2]) = (1,(\s,\s)),
    \]
    which commutes with $b_2 = (1,\s)$ and $b_3 = (1,(1,\s))$.
    Therefore $\br = \wh{\Mpre}_{m,n}$ in this case as well.
    Thus (1) implies that $\br^d \subseteq \Mpre$ in cases $(B_2)$ and $(C)$.

    In cases $(A)$, $(B)$, and $(C)$ the subgroup $\br$ has finite index by Lemma~\ref{lemma: N quotients}.
    Hence $\Mpre$ is a regular branch group in each of those cases.
\end{proof}

\subsection{Characterizing the level one stabilizer}

We previously observed that $\st_1 \Mper \subseteq \Mper^d$ and $\st_1\Mpre \subseteq \Mpre^d$ (see the remark before Lemma~\ref{lemma: stabilizer projections are surjective}). 
In this section we characterize $\st_1 \Mper$ and $\st_1 \Mpre$ as subgroups of $\Mper^d$ and $\Mpre^d$, respectively. In the pre-periodic case, the branching property is essential.

\begin{lemma}
\label{lemma: involutions}
    There exists an involution $\tau : \Mpre/\br \to \Mpre/\br$ uniquely determined by 
    \[
        \tau(b_m) \equiv b_n \bmod \br.
    \]
\end{lemma}

\begin{proof}
    The quotient $\Mpre/\br$ is generated by $b_m$ and $b_n$, so any involution exchanging them is unique.
    
    In case $(A)$, observe that $\Mpre/\br \cong (\ZZ/d\ZZ)^2$ where the generators $b_m$ and $b_n$ are identified with $(1,0)$ and $(0,1)$, respectively, and exchanging them is an involution.
    
    If $(B)$ or $(C)$, then $\Mpre/\br \cong \H_d$ by Lemma~\ref{lemma: N quotients}, under which $b_m$ and $b_n$ are sent to the $g_1$ and $g_2$ from the presentation of Lemma~\ref{lemma: H_d presentation}. 
    In Lemma~\ref{lemma: H_d involution}, we observed that $\H_d$ has an involution exchanging $g_1$ and $g_2$, and therefore $\Mpre/\br$ has an involution exchanging $b_m$ and $b_n$.

    If $(D)$, then $\br = 1$ and Lemma~\ref{lemma: N quotients} implies that $\Mpre$ is isomorphic to the pro-2 dihedral group.
    There is a unique involution $\tau$ of $\Mpre$ satisfying $\tau(b_1) = b_2$.
\end{proof}

\begin{lemma}
\label{lemma: small truncation B}
    If $\ell \leq n$, then $\Mpre =_\ell [C_d]^\ell$.
\end{lemma}

\begin{proof}
    We proceed by induction on $\ell$.
    The base case $\ell = 0$ is trivial, so suppose that $0 < \ell \leq n$ and the claim is true for $\ell - 1$. 
    Since $b_i =_\ell 1$ if and only if $\ell < i$, it follows that each $b_i$ is supported in at most one component in $\rho_\ell(\Mpre)$ and that
    $
        \st_1 \Mpre =_\ell \Mpre^d.
    $
    Our inductive hypothesis implies that $\Mpre =_{\ell-1} [C_d]^{\ell-1}$, hence 
    \[
        \st_1 \Mpre =_\ell ([C_d]^{\ell-1})^d = \st_1 [C_d]^\ell.
    \]
    Since $b_1 = \s$, we conclude that $\Mpre =_\ell [C_d]^\ell$.
\end{proof}

\begin{proposition}
\label{prop: stabilizer criteria}
The level one stabilizer in the periodic case is
\[
    \st_1 \Mper = \{(g_1,\ldots,g_d) \in \Mper^d : \eta_n(g_i) = \eta_n(g_j) \text{ for all $1 \leq i, j \leq n$}\}.
\]

The level one stabilizer in the pre-periodic cases is
\[
    \st_1 \Mpre = 
    \begin{cases}
        \{(g_1,\ldots,g_d) \in \Mpre^d : \chi_m(g_i) = \chi_n(g_{i+\omega}) \text{ for all $1 \leq i \leq n$}\} & \text{if $(A)$,}\\
        \{(g_1,\ldots,g_d) \in \Mpre^d : \tau(g_i) \equiv g_{i+\omega} \bmod \br \text{ for all $1 \leq i \leq n$}\} & \text{if $(B)$, $(C)$, or $(D)$.}
    \end{cases}
\]
\end{proposition}

\begin{proof}
    First consider the periodic case.
    Let $\SS \subseteq \Mper^d$ be the subgroup defined by
    \[
        \SS := \{(g_1,\ldots,g_d) \in \Mper^d : \eta_n(g_i) = \eta_n(g_j) \text{ for all $1 \leq i, j \leq n$}\}.
    \]
    Lemma \ref{lemma: A subgroup properties}(3) implies that $\st_1\Mper = \lprof a_1^d\rprof \wh\Mper_n^d$.
    Since $a_1^d = (a_n,\ldots,a_n)$ and $\wh\Mper_n = \ker(\eta_n)$, we have $\st_1\Mper \subseteq \SS$.
    On the other hand, if $g := (g_1,\ldots,g_d) \in \SS$, then Lemma \ref{lemma: A subgroup properties}(2) implies that for each $i$ we can write $g_i = a_n^{k_i} h_i$ for some $k_i \in \ZZ_d$ and some $h_i \in \wh\Mper_n$.
    The definition of $\SS$ implies that $k_i = \eta_n(g_i) = \eta_n(g_j) = k_j$ for all $i,j$.
    Hence $g = a_1^{dk_1}(h_1,\ldots,h_d) \in \st_1\Mper$.
    Therefore $\st_1\Mper = \SS$.

    Next consider the preperiodic case.
    Suppose $(A)$.
    Let $\SS \subseteq \Mpre^d$ be the subgroup defined by 
    \[
        \SS := \{(g_1,\ldots, g_d) \in \Mpre^d : \chi_m(g_i) = \chi_n(g_{i+\omega}) \text{ for all $i$}\}.
    \]
    If $i \neq 1, m + 1$, then 
    \[
        b_i \in \wh{\Mpre}_{1,m+1} = \wh{\Mpre}_{m,n}^d \subseteq \SS,
    \]
    where the equality follows from Lemma~\ref{lemma: pre-branch}.
    Recall that $b_{m+1} = (1,\ldots, 1, b_n, 1, \ldots, 1,b_m)$ where $b_n$ is in the $\omega$th component.
    Hence $b_{m+1} \in \SS$.
    The group $\SS$ is closed under conjugation by elements of $\Mpre$, thus $\st_1\Mpre = \wh{\Mpre}_1 \subseteq \SS$.

    Now suppose that $g = (g_1,\ldots, g_d) \in \SS$.
    Consider the element
    \[
        h := (h_1, \ldots, h_d) = b_{m+1,1}^{\chi_m(g_1)}\cdots b_{m+1,d}^{\chi_m(g_d)} \in \st_1\Mpre.
    \]
    Proposition~\ref{prop: branch} implies that $\Mpre/\br$ is abelian, hence
    \[
        h_{i+\omega} \equiv b_m^{\chi_m(g_{i+\omega})}b_n^{\chi_m(g_{i})}
	\equiv b_m^{\chi_m(g_{i+\omega})}b_n^{\chi_n(g_{i+\omega})}
	\equiv g_{i+\omega} \bmod \br,
    \]
    where the second congruence follows from $g \in \SS$.
    Thus $g \equiv h \bmod \br^d$ and Proposition~\ref{prop: branch} implies $g \in \st_1\Mpre$.
    Therefore $\st_1\Mpre = \SS$ in case $(A)$.

    Next suppose either $(B)$ or $(C)$.
    Let $\SS \subseteq \Mpre^d$ be the subgroup defined by
    \[
        \SS := \{(g_1,\ldots, g_d) \in \Mpre^d : \tau(g_i) \equiv g_{i+\omega} \bmod \br \text{ for all $i$}\}.
    \]
    Since $\wh{\Mpre}_{m,n} \subseteq \br$, for any $i \neq 1, m + 1$,
    \[
        b_i \in \wh{\Mpre}_{1,m+1} = \wh{\Mpre}_{m,n}^d \subseteq \SS.
    \]
    The definition of $\tau$ implies that
    \[
    	b_{m+1}  \equiv (1,\ldots, 1, \tau(b_m), 1, \ldots, 1,b_m) \bmod \br^d,
    \]
    hence that $b_{m+1} \in \SS$.
    The group $\SS$ is closed under conjugation by $\Mpre$ and $\st_1\Mpre = \wh{\Mpre}_1$, hence $\st_1\Mpre \subseteq \SS$.

    For the reverse inclusion, suppose that $g = (g_1,\ldots,g_d) \in \SS$.
    Lemma~\ref{lemma: N quotients} and the presentation Lemma~\ref{lemma: H_d presentation} implies that for each $g_i$ we have
    \[
        g_i \equiv b_m^{r_i}b_n^{s_i}[b_m,b_n]^{t_i} \bmod \br,
    \]
    for some unique $r_i, s_i, t_i \in \ZZ/d\ZZ$.
    Since $g \in \SS$, it follows that
    \[
        g_{i+\omega} \equiv \tau(g_i) \equiv b_n^{r_i}b_m^{s_i}[b_n,b_m]^{t_i} \bmod \br.
    \]
    Recall that $\omega = d/2$ in cases $(B)$ and $(C)$, hence $\pi_i(b_{m+1,i}) = b_m$ and $\pi_i(b_{m+1,i+\omega}) = b_n$.
    For $1 \leq i \leq \omega$, let
    \[
        h_i' := b_{m+1,i}^{r_i}b_{m+1,i+\omega}^{s_i}[b_{m+1,i},b_{m+1,i+\omega}]^{t_i} \in \Mpre.
    \]
    Then $h_i'$ is supported in the $i$th and $(i + \omega)$th coordinates and
    \begin{align*}
        \pi_i(h_i') &= b_m^{r_i}b_n^{s_i}[b_m,b_n]^{t_i} \equiv g_i \bmod \br\\
        \pi_{i+\omega}(h_i') &= b_n^{r_i}b_m^{s_i}[b_n,b_m]^{t_i} \equiv g_{i+\omega} \bmod \br.
    \end{align*}
    Let $h = (h_1,\ldots,h_d) := h_1'\cdots h_\omega' \in \Mpre$.
    The above congruences imply that $g \equiv h \bmod \br^d$, hence Proposition~\ref{prop: branch} implies that $g \in \st_1\Mpre$.
    Therefore $\st_1\Mpre = \SS$.

    Finally, suppose $(D)$.
    Lemma~\ref{lemma: N quotients} implies that $\br = 1$ and $\Mpre$ is pro-2 dihedral.
    Recall that $b_\infty := b_2b_1$.
    Hence every element of $\Mpre$ may be uniquely expressed as $b_\infty^\ve \s^i$ for some $\ve \in \ZZ_2$ and $i \in \ZZ/2\ZZ$.
    Let 
    \[
        \SS := \{(g,\tau(g)) \in \Mpre^2\}.
    \]
    The stabilizer $\st_1\Mpre$ is topologically generated by $b_\infty^2 = (b_\infty,b_\infty^{-1})$ and $b_\infty\s = b_2 = (b_2,b_1)$, both of which belong to $\SS$.
    Hence $\st_1\Mpre \subseteq \SS$.
    On the other hand,
    \[
        (b_\infty^{\ve},\tau(b_\infty^\ve)) = (b_\infty^\ve,b_\infty^{-\ve}) = b_\infty^{2\ve} \in \Mpre,
    \]
    and
    \[
        (b_\infty^{\ve}\s, \tau(b_\infty^{\ve}\s)) = (b_\infty^\ve\s, b_\infty^{-\ve}b_2) = b_\infty^{2\ve-1} \s \in \Mpre.
    \]
    Thus $\st_1\Mpre = \SS$.
\end{proof}

In Proposition \ref{prop: branch} we showed that $\br^d \subseteq \Mpre$.
Since $\br$ is a normal subgroup of $\Mpre$, it follows that $\br^d$ is also normal in $\Mpre$.
We identify the quotients $\Mpre/\br^d$ in Lemma \ref{lemma: N^d quotient}; the index of $\br^d$ is essential for our calculation of the orders of finite level truncations of $\Mpre$.

\begin{lemma}
\label{lemma: N^d quotient}
The quotients $\Mpre/\br^d$ are as follows
    \begin{enumerate}
        \item If $(A)$, then $\Mpre/\br^d \cong [C_d]^2$.
        \item If $(B)$ or $(C)$, then $\Mpre/\br^d \cong  \langle \s \rangle \ltimes (\Mpre/\br)^\omega$ where $\s$ acts on $(\Mpre/\br)^\omega$ by
        \[
            \s\inv (g_1, \ldots, g_\omega) \s := (g_2,\ldots, g_{\omega},\tau(g_1)).
        \]
    \end{enumerate}
    Hence
    \[
        [\Mpre : \br^d] =
        \begin{cases}
            d^{d+1} & \text{if $(A)$},\\
            d^{3d/2 + 1} & \text{if $(B)$ or $(C)$}.
        \end{cases}
    \]
\end{lemma}

\begin{proof}
    Lemma~\ref{lemma: pre-branch} says $\wh{\Mpre}_{1,m+1}  = \wh \Mpre_{m,n}$ while
    Proposition~\ref{prop: branch} implies that $(\wh \Mpre_{m,n})^d\subseteq \br^d$. 
    Hence $\Mpre/\br^d$ factors through $\Mpre/\wh \Mpre_{1,m+1}$ and hence is generated by $b_1=\s$ and $b_{m+1}$.
    Furthermore,
    \[
        \st_1 \Mpre/\br^d \subseteq \Mpre^d/\br^d \cong (\Mpre/\br)^d.
    \]
    
    (1)
    Suppose $(A)$.
    Since $\Mpre/\br$ is abelian, it follows that $b_{m+1,i}$ and $b_{m+1,j}$ commute modulo $\Mpre/\br^d$ for all $i$ and $j$.
    Thus $\Mpre/\br^d$ is a quotient of $[C_d]^2$.
    On the other hand, suppose that
    \[
        \s^j(g_1, \ldots, g_d) := \s^jb_{m+1,1}^{e_1}\cdots b_{m+1,d}^{e_d}.
    \]
    Then $\chi_m(g_i) = e_i$, which by Lemma~\ref{lemma: N quotients}(1) implies that these $d^{d+1}$ elements are all distinct in $\Mpre/\br^d$.
    Therefore $\Mpre/\br^d \cong [C_d]^2$ and $[\Mpre : \br^d] = d^{d+1}$ in case $(A)$.

    (2) Suppose $(B)$ or $(C)$.
    Proposition~\ref{prop: stabilizer criteria}(2) implies that $\s^j(g_1, \ldots, g_d) \in \Mpre$ if and only if $\tau(g_i) \equiv g_{i+\omega} \bmod \br$ for all $i$.
    Hence the map $\s^j(g_1,\ldots, g_d) \mapsto \s^j(g_1,\ldots, g_\omega)$ defines a surjective homomorphism from $\Mpre$ onto 
    $\langle \s \rangle \ltimes (\Mpre/\br)^\omega$ with kernel $\br^d$.
    Since $|\Mpre/\br| = d^3$ by Lemma~\ref{lemma: N quotients}(2) and $\omega = d/2$ in cases $(B)$ and $(C)$, we conclude that $[\Mpre : \br^d] = d^{3d/2 + 1}$.
\end{proof}

\subsection{Semirigidity}

The following lemma shows that the generators of $\Mper$ and $\Mpre$ satisfy a weak form of rigidity in the sense of the inverse Galois problem: if our distinguished generators are replaced by arbitrary conjugates, they still generate the group.
\begin{lemma}
\label{lemma: generator rigidity}
If $u_i \in \Mper$ and $v_i \in \Mpre$ for $1 \leq i \leq n$ are arbitrary elements, then
\begin{align*}
    \Mper &= \lprof u_1a_1u_1^{-1}, \ldots, u_na_nu_n^{-1} \rprof,\\
    \Mpre &= \lprof v_1b_1v_1^{-1}, \ldots, v_nb_nv_n^{-1} \rprof.
\end{align*}
\end{lemma}

\begin{proof}
    First consider $\Mper$.
    Let $a_i' := u_i a_i u_i^{-1}$ and define $\Mper' := \lprof a_1', \ldots, a_n' \rprof \subseteq \Mper$.
    We claim that $\Mper' = \Mper$.
    The groups $\Mper$ and $\Mper'$ are both closed, hence it suffices by Proposition~\ref{prop: inductive criteria} to prove that $\Mper' =_\ell \Mper$ for all $\ell \geq 0$; we proceed by induction on $\ell$.

    The base case $\ell = 0$ is immediate since $\Mper$ and $\Mper'$ are both trivial at level $\ell = 0$.
    Now suppose that $\ell \geq 1$ and $\Mper' =_{\ell -1} \Mper$.
    Since $u\Mper u^{-1} = \Mper$ for all $u \in \Mper$, we may conjugate everything by $u_1^{-1}$ to assume without loss of generality that $a_1 = a_1' \in \Mper'$ and $u_1 = 1$.
    Thus it suffices to show that $a_i \in_\ell \Mper'$ for all $2 \leq i \leq n$.

    If $1 \leq i < n$, then after potentially replacing $a_{i+1}'$ with a conjugate by a power of $a_1'$, there is some element $w_{i} \in \Mper$ such that
    \[
        a_{i+1}' = u_{i+1}a_{i+1}u_{i+1}^{-1} = (1,\ldots, 1, w_{i}a_{i}w_{i}^{-1}).
    \]
    We also have
    \[
        a_1'^d = a_1^d = (a_n,\ldots,a_n).
    \]
    Let $\pi_d : \st_1\Mper' \to \Mper$ denote projection onto the $d$th coordinate.
    The above observations imply that $\pi_d(\st_1\Mper')$ is a subgroup of $\Mper$ containing $\Mper$-conjugates of each of the $a_i$, hence $\pi_d(\st_1\Mper') =_{\ell-1} \Mper$ by our inductive hypothesis.
    Thus for $2 \leq i \leq n$ there exists elements $w_{i}' \in \st_1\Mper'$ such that $\pi_d(w_{i}') =_{\ell-1} w_{i-1}^{-1}$, hence
    \[
        a_i = (1,\ldots,1,a_{i-1}) =_\ell w_{i}'a_{i}'w_{i}'^{-1} \in \Mper'.
    \]
    That completes the induction, hence proves that
    $
        \Mper = \lprof u_1a_1u_1^{-1}, \ldots, u_na_nu_n^{-1} \rprof.
    $
    
    Next consider $\Mpre$.
    Let $b_i' := v_i b_i v_i^{-1}$ and define $\Mpre' := \lprof b_1', \ldots, b_n' \rprof \subseteq \Mpre$.
    We claim that $\Mpre' = \Mpre$.
    The groups $\Mpre$ and $\Mpre'$ are both closed, hence it suffices by Proposition~\ref{prop: inductive criteria} to prove that $\Mpre' =_\ell \Mpre$ for all $\ell \geq 0$; we proceed by induction on $\ell$.

    The base case $\ell = 0$ is again immediate; suppose that $\ell \geq 1$ and $\Mpre' =_{\ell -1} \Mpre$.
    Since $v\Mpre v^{-1} = \Mpre$ for all $v \in \V$, we may conjugate everything by $v_1^{-1}$ to assume without loss of generality that $\s = b_1 = b_1' \in \Mpre'$ and $v_1 = 1$.
    Thus it suffices to show that $b_i \in_\ell \Mpre'$ for all $2 \leq i \leq n$.

    If $1 \leq i < n$ and $i \neq m$, then after potentially replacing $b_{i+1}'$ with a conjugate by a power of $\s$, there is some element $w_i \in \Mpre$ such that
    \[
        b_{i+1}' = v_{i+1}b_{i+1}v_{i+1}^{-1} = (1,\ldots, 1, w_ib_iw_i^{-1}).
    \]
    Similarly after potentially replacing $b_{m+1}'$ with a conjugate by a power of $\s$, there are elements $w_m, w_n \in \Mpre$ such that
    \[
        b_{m+1}' = v_{m+1}b_{m+1}v_{m+1}^{-1} = (1,\ldots,1,w_nb_nw_n^{-1},1,\ldots,1,w_mb_mw_m^{-1}).
    \]
    Conjugating by powers of $\s$ gives us all cyclic shifts of these elements in $\Mpre'$.
    It follows that $\pi_d(\st_1\Mpre')$ contains $\Mpre$-conjugates of $b_i$ for each $1\leq i \leq n$.
    Thus $\pi_d(\st_1\Mpre') =_{\ell-1} \Mpre$ by our inductive hypothesis.

    If $i \neq 1, m+1$ and $v \in \Mpre$, then $\pi_d(\st_1\Mpre') =_{\ell-1} \Mpre$ implies that there exists a $w_i' \in \st_1\Mpre'$ such that $\pi_d(w_i') =_{\ell-1} vw_{i-1}^{-1}$.
    Hence
    \[
        w_i'b_i'w_i'^{-1} =_\ell (1,\ldots,1, vb_{i-1}v^{-1}) \in_\ell \Mpre'
    \]
    and we conclude that $\wh{\Mpre}_{1,m+1} \subseteq_\ell \Mpre'$.
    Observe that
    \begin{align*}
        b_{m+1}' 
        &= v_{m+1}b_{m+1}v_{m+1}^{-1}\\
        &= (1,\ldots,1,w_nb_nw_n^{-1},1,\ldots,1,w_mb_mw_m^{-1})\\
        &= (1\ldots,1,[w_n,b_n]b_n,1,\ldots,1,[w_m,b_m]b_m)\\
        &= (1,\ldots,1,[w_n,b_n],1,\ldots,1,[w_m,b_m])b_{m+1}
    \end{align*}
    Without loss of generality we may assume that $v_{m+1} \in \wh{\Mpre}_{m+1}$.
    Since $\wh{\Mpre}_{m+1} = (\wh{\Mpre}_{m,n})^d$, it follows that $w_m, w_n \in \wh{\Mpre}_{m,n}$.
    Note that $[w_m, b_m]$ and $[w_n, b_n]$ belong to $\wh{\Mpre}_{m,n}$ since $\wh{\Mpre}_{m,n}$ is a normal subgroup of $\Mpre$.
    Therefore $b_{m+1} \in_\ell \Mpre'$.
    This completes the induction and proves that
    $
        \Mpre = \lprof v_1b_1v_1^{-1}, \ldots, v_nb_nv_n^{-1} \rprof.
    $
\end{proof}

\begin{remark}
    When $d = p$ is prime, the groups $\Mper, \Mper', \Mpre, \Mpre'$ from the proof of Lemma~\ref{lemma: generator rigidity} are pro-$p$ groups.
    In that case, $\Mper' = \Mper$ and $\Mpre' = \Mpre$ follow from observing that $\Mper'$ and $\Mpre'$ surject onto $\Mper^{ab}$ and $\Mpre^{ab}$, combined with the observation that the Frattini subgroup of a $p$-group contains the commutator; semirigidity holds for all pro-$p$ groups.
    This is the approach taken by Pink when analyzing the $d = 2$ case (see \cite[Lem. 1.3.2]{PinkPCF}).
\end{remark}

We now prove the first main result of this section.
This result generalizes the main assertion of \cite[Thm. 2.4.1]{PinkPCF} in the case of a periodic critical point.

\begin{theorem}[Semirigidity of $\Mper$]
    \label{theorem: A semi rigidity}
    If $a_i' \in [C_d]^\infty$ are elements such that $a_i' \sim a_i'$ in $[C_d]^\infty$ for $1 \leq i \leq n$, then there exists an element $w \in [C_d]^\infty$ and elements $u_i \in \Mper$ such that for each $i$,
    \[
        wa_i'w^{-1} = u_ia_iu_i^{-1}.
    \]
    In particular,
    $w\lprof a_1',\ldots, a_n'\rprof w^{-1} = \Mper.$
\end{theorem}

\begin{proof}
    For $\ell \geq 0$ let $(*_\ell)$ be the claim,

    \begin{quote}
        $(*_\ell)$: If $a_i' \in [C_d]^\infty$ are elements such that $a_i' \sim_\ell a_i$ in $[C_d]^\infty$ for $1 \leq i \leq n$, then there exists an element $w \in [C_d]^\infty$ and elements $u_i \in \Mper$ such that $wa_i'w^{-1} =_\ell u_ia_iu_i^{-1}$ for each $i$.
    \end{quote}

    If $(*_\ell)$, then Proposition~\ref{prop: inductive criteria} implies that proving $(*_\ell)$ for all $\ell \geq 0$ furnishes $u_i$ and $w$ as in the statement of the theorem, while Lemma~\ref{lemma: generator rigidity} implies that $\Mper =_\ell \lprof u_1a_1u_1^{-1}, \ldots, u_na_nu_n^{-1}\rprof$. So it suffices to verify \((*_\ell)\), which we do by induction on $\ell$. The base case is immediate.
    
    Suppose that $\ell \geq 1$ and that $(*_{\ell-1})$ holds.
    
    Simultaneously conjugating all of the $a_i'$ we may assume without loss of generality that $a_1 = a_1'$. 
    If $2 \leq i \leq n$, then any conjugate of $a_{i} = (1,\ldots,1,a_{i-1})$ has a single non-trivial component.
    Hence there is an integer $j_{i}$ such that
    \[
        a_1^{j_{i}}a_{i}'a_1^{-j_{i}} =_\ell (1,\ldots,1,a_{i-1}'')
    \]
    for some element $a_{i-1}'' \in [C_d]^\infty$.
    Proposition~\ref{prop: general wreath conj crit} implies that $a_{i-1} \sim_{\ell-1} a_{i-1}''$ in $[C_d]^\infty$.
    Let $a_n'' := a_n$, so that $a_i \sim_{\ell - 1} a_i''$ for all $1\leq i \leq n$.
    Then $(*_{\ell-1})$ implies there is an element $w' \in [C_d]^\infty$ and elements $u_i' \in \Mper$ such that for each $i$,
    \[
        w'a_i''w'^{-1} =_{\ell-1} u_i'a_iu_i'^{-1}
    \]
    Define $w := (u_n'^{-1}w', \ldots, u_n'^{-1}w') \in [C_d]^\infty$;
    note that $w$ commutes with $\s$.
    Let $u_1 = 1$, then
    \begin{align*}
        wa_1'w^{-1} &= wa_1w^{-1}\\
        &= \s w(1,\ldots,1,a_n)w^{-1}\\
        &= \s(1,\ldots,1,u_n'^{-1}(w'a_n''w'^{-1})u_n')\\
        &=_\ell \s (1,\ldots,1,a_n)\\
        &= u_1a_1u_1^{-1}.
    \end{align*}
    In particular, $w$ commutes with $a_1$.
    
    Next suppose that $2 \leq i \leq n$.
    Since $u_n'^{-1}u_{i-1}' \in \Mper$, Lemma~\ref{lemma: stabilizer projections are surjective} implies there is a $v_i' \in \st_1\Mper$ such that $\pi_d(v_i) = u_n'^{-1}u_{i-1}'$.
    Therefore,
    \begin{align*}
        w(a_1^{j_i}a_i'a_1^{-j_i})w^{-1} &=_\ell w(1,\ldots,1,a_{i-1}'')w^{-1}\\
        &= (1,\ldots,1,u_n'^{-1}w'a_{i-1}''w'^{-1}u_n')\\
        &=_\ell (1,\ldots,1,u_n'^{-1}u_{i-1}'a_{i-1}u_{i-1}'^{-1}u_n')\\
        &= v_i(1,\ldots,1,a_{i-1})v_i^{-1}\\
        &=v_ia_iv_i^{-1}.
    \end{align*}
    Let $u_i := a_1^{-j_i}v_i$.
    Since $a_1$ commutes with $w$, the above identity is equivalent to
    \[
        wa_i'w^{-1} =_\ell u_ia_iu_i^{-1}.
    \]
    Then $u_i \in \Mper$ for all $i$ and
    $
        wa_i'w^{-1} =_\ell u_ia_iu_i^{-1}
    $
    which completes our induction.
\end{proof}

Next we prove the analog of Theorem~\ref{theorem: A semi rigidity} for the group $\Mpre$.
First, a technical lemma.

\begin{lemma}
\label{lemma: stabilizer technical}
    If $u_1, u_2 \in \Mpre$, then there exists an element $v \in \Mpre$ and $i, j \in \ZZ/d\ZZ$ such that 
    $
        (1,\ldots,1,vu_1b_n^i,1,\ldots,1,vu_2b_m^j) \in \Mpre,
    $
    where the $vu_1b_n^i$ term is in the $\omega$th component.
\end{lemma}

\begin{proof}
    We proceed by cases.
    First suppose $(A)$.
    Then $(1,\ldots,1,vu_1b_n^i,1,\ldots,1,vu_2b_m^j) \in \Mpre$ is equivalent, by Proposition~\ref{prop: stabilizer criteria}, to
    \[
        \chi_n(v) + \chi_n(u_1) + i = \chi_m(v) + \chi_m(u_2) + j,
    \]
    and either
    \[
    \chi_m(v) + \chi_m(u_1) =
    \begin{cases}
         \chi_n(v) + \chi_n(u_2) & \text{if }\omega = d/2,\\
        0 & \text{if }\omega \neq d/2.
    \end{cases}
    \]
    Either system of equations has a solution for $i, j \in \ZZ/d\ZZ$ and $v \in \Mpre$.

    Next suppose $(B)$, $(C)$, or $(D)$.
    Then $(1,\ldots,1,vu_1b_n^i,1,\ldots,1,vu_2b_m^j) \in \Mpre$ is equivalent, by Proposition~\ref{prop: stabilizer criteria}, to 
    \begin{equation}
    \label{eqn: tau constraint}
        \tau(vu_1b_n^i) \equiv vu_2b_m^j \bmod \br.
    \end{equation}
    Let $w_1 := vu_1$ and $w_2 := u_1\inv u_2$, so that \eqref{eqn: tau constraint} becomes
    \[
        \tau(w_1)b_m^i \equiv \tau(w_1b_n^i) \equiv w_1w_2b_m^j   \bmod \br,
    \]
    or equivalently,
    \[
        w_1\inv \tau(w_1) \equiv w_2 b_m^{j-i} \bmod \br.
    \]
    Writing $w_1 = b_m^a b_n^b [b_m,b_n]^c$ and $w_2 = b_m^rb_n^s[b_m,b_n]^t$ this expands to
    \[
        b_n^{a-b}b_m^{b - a}[b_m,b_n]^{-a^2 - 2c} \equiv b_n^s b_m^{r+j-i}[b_m,b_n]^{t + s(j-i)} \bmod \br.
    \]
    Comparing exponents, we then see that \eqref{eqn: tau constraint} is equivalent to the following system of congruences having a solution for $a, b, c, i, j$ in terms of $r, s, t$
    \begin{align*}
        a - b &\equiv s \bmod d\\
        b - a &\equiv r + j - i \bmod d\\
        -a^2 - 2c &\equiv t + s(j - i) \bmod d.
    \end{align*}
    The first two congruences determine $b$ and are equivalent to $i - j \equiv r + s \bmod d$.
    Hence the system reduces to solving
    \[
        a^2 + 2c \equiv (r + s)s - t \bmod d.
    \]
    This congruence may be solved by setting $a \equiv 0, 1 \bmod d$ depending on the parity of the right hand side.
    In particular, a solution always exists.
\end{proof}

\begin{theorem}[Semirigidity of $\Mpre$]
    \label{theorem: B semi rigidity}
    If $b_i' \in [C_d]^\infty$ are elements such that $b_i' \sim b_i$ in $[C_d]^\infty$ for $1 \leq i \leq n$, then there exists an element $w \in [C_d]^\infty$ and elements $u_i \in \Mpre$ such that for each $i$,
    \[
        wb_i'w^{-1} = u_ib_iu_i^{-1}.
    \]
    In particular, $\lprof b_1', \ldots, b_n'\rprof = w \Mpre w^{-1}$.
\end{theorem}

\begin{proof}
    As in the proof of Theorem~\ref{theorem: A semi rigidity}, by combining Proposition~\ref{prop: inductive criteria} and Lemma~\ref{lemma: generator rigidity}, it suffices to prove the following proposition for all $\ell \geq 0$,
    
    \begin{quote}
    \noindent $(*_\ell)$: If $b_i' \in [C_d]^\infty$ are elements such that $b_i' \sim_\ell b_i$ in $[C_d]^\infty$ for $1 \leq i \leq n$, then there exists an element $w \in [C_d]^\infty$ and elements $u_i \in \Mpre$ such that 
    $wb_i'w^{-1} =_\ell u_ib_iu_i^{-1}$.
    \end{quote}

    We proceed by induction on $\ell$.
    The base case is trivial; suppose $\ell \geq 1$ and that $(*_{\ell-1})$ is true.
    Replacing all the $b_i$ by a simultaneous conjugation we may assume that $b_1' = b_1 = \s$.
    If $i \neq 1, m + 1$, then $b_i' \sim_\ell b_i$ implies there is some integer $j_i$ such that
    \[
        b_i' =_\ell \s^{j_i}(1,\ldots,1,b_{i-1}'')\s^{-j_i},
    \]
    where $b_{i-1}'' \in [C_d]^\infty$.
    Similarly, $b_{m+1}' \sim_\ell b_{m+1}$ implies there is some $j_{m+1}$ such that
    \[
        b_{m+1}' =_\ell \s^{j_{m+1}}(1,\ldots,1,b_n'',1,\ldots,1,b_m'')\s^{-j_{m+1}},
    \]
    where $b_m'', b_n'' \in [C_d]^\infty$ and the $b_n''$ is in the $\omega$th coordinate.
    Proposition~\ref{prop: general wreath conj crit} implies that $b_i \sim_{\ell-1} b_i''$ for each $1\leq i \leq n$.
    Thus by $(*_{\ell-1})$ we have an element $w' \in [C_d]^\infty$ and elements $u_i' \in \Mpre$ such that for each $1 \leq i \leq n$,
    \[
        w'b_i''w'^{-1} =_{\ell-1} u_i'b_iu_i'^{-1}.
    \]

    Let $v \in \Mpre$ and $i, j \in \ZZ/d\ZZ$ be the elements provided by Lemma~\ref{lemma: stabilizer technical} associated to $u_n', u_m' \in \Mpre$ such that
    \[
        v_{n+1} := (1,\ldots,1,vu_n'b_n^i,1\ldots,1,vu_m'b_m^j) \in \Mpre.
    \]
    Define $w := (vw',\ldots,vw') \in [C_d]^\infty$; note that $w$ commutes with $\s$.
    Lemma~\ref{lemma: stabilizer projections are surjective} implies that for $i\neq 1, m +1$ there exists an element $v_i \in \st_1\Mpre$ such that $\pi_d(v_i) = vu_{i-1}'$.
    Let $u_1 = 1$ and for $i \neq 1$ let $u_i = \s^{j_i}v_i$.
    Hence each $u_i \in \Mpre$

    Since $w$ commutes with $\s$ and $u_1 = 1$, we have
    \[
        wb_1'w^{-1} = \s = u_1b_1u_1^{-1}.
    \]
    If $i = m + 1$, then
    \begin{align*}
        wb_{m+1}'w^{-1} &=_\ell w\s^{j_{m+1}}(1,\ldots,1,b_n'',1\ldots,1,b_m'')\s^{-j_{m+1}}w^{-1}\\
        &=\s^{j_{m+1}}(1,\ldots,1,vw'b_n''w'^{-1}v^{-1},1,\ldots,1,vw'b_m''w'^{-1}v^{-1})\s^{-j_{m+1}}\\
        &=_\ell\s^{j_{m+1}}(1,\ldots,1,vu_n'b_nu_n'^{-1}v^{-1},1,\ldots,1,vu_m'b_mu_m'^{-1}v^{-1})\s^{-j_{m+1}}\\
        &= \s^{j_{m+1}}v_{m+1}b_{m+1}v_{m+1}^{-1}\s^{-j_{m+1}}\\
        &= u_{m+1}b_{m+1}u_{m+1}^{-1}.
    \end{align*}
    Finally, if $i\neq 1, m+1$, then
    \begin{align*}
        wb_i'w^{-1} &=_\ell w\s^{j_i}(1,\ldots,1,b_{i-1}'')\s^{-j_i}w^{-1}\\
        &= \s^{j_i}(1,\ldots,1,vw'b_{i-1}''w'^{-1}v^{-1})\s^{-j_i}\\
        &=_\ell \s^{j_i}(1,\ldots,1,vu_{i-1}'b_{i-1}u_{i-1}'^{-1}v^{-1})\s^{-j_i}\\
        &= \s^{j_i}v_ib_iv_i^{-1}\s^{-j_i}\\
        &= u_ib_iu_i^{-1}.
    \end{align*}
    That completes our induction, hence our proof.
\end{proof}

Semirigidity allows us to identify $\barb f$ with $\Mper$ or $\Mpre$ in the periodic and preperiodic cases, respectively

\begin{corollary}
\label{corollary: final IMG calculation}
    Let $f \in K[x]$ be a unicritical PCF polynomial with degree $d$ coprime to $\cchar K$.
    Suppose that $f$ has $n$ distinct finite post-critical points.
    Let $\barb f = \lprof c_1, c_2, \ldots, c_n\rprof$ be as defined in Section~\ref{section: unicritical PCF}.
    \begin{enumerate}
        \item In the periodic case, there exists a $w \in [C_d]^\infty$ and $u_i \in \Mper = \Mper(d,n)$ such that $w c_i w\inv = u_i a_i u_i\inv$ for each $1 \leq i \leq n$.
        In particular, $w\barb f w\inv = \Mper$.

        \item In the preperiodic case, there exists a $w \in [C_d]^\infty$ and $u_i \in \Mpre = \Mpre(d,m,n,\omega)$ such that $w c_i w\inv = u_i b_i u_i\inv$ for each $1 \leq i \leq n$.
        In particular, $w \barb f w\inv = \Mpre$.
    \end{enumerate}
\end{corollary}

\begin{proof}
    First consider the periodic case.
    Proposition~\ref{prop: geometric generators} implies that the generators $c_i$ satisfy the same cyclic conjugate recurrences as the $a_i$.
    Hence Proposition~\ref{prop: weak recursion} implies $c_i \sim a_i$ in $[C_d]^\infty$ for each $i$.
    Therefore Theorem~\ref{theorem: A semi rigidity} yields the conclusion.
    The preperiodic case is the same, \emph{mutatis mutandis}.
\end{proof}

\begin{remark}
    Recall that $c_i \in \barb f$ is the image of an inertia generator over the point $p_i = f^i(0)$.
    If $v \in \barb f$ is any element, then $vc_iv\inv$ is also an inertia generator over $p_i$.
    Thus after replacing $c_i$ by the appropriate conjugates, Corollary \ref{corollary: final IMG calculation} implies that $\barb f = \lprof c_1, c_2,\ldots, c_n\rprof$ where the $c_i$ are inertia generators and satisfy the recursive identities defining the model groups.
    Note, however, that with this choice of generators we typically will not have $c_1c_2\cdots c_n$ equal to the standard odometer.
\end{remark}

\section{Finite level truncations, Hausdorff dimension, and normalizers}\label{section: further properties}

The semirigidity result of the preceding section establishes that $\barb f$ is isomorphic to a model group $\Mper(d,n)$ or $\Mpre(d,m,n,\omega)$ according to the combinatorics of its critical orbit. 
In this section we establish additional properties of the model groups---hence $\barb f$---including determining the order of their finite level truncations, calculating their Hausdorff dimension, and analyzing the structure of their normalizers in $[C_d]^\infty$.
The latter is essential to our determination of the constant field extensions in Section \ref{section: constant fields}.

\subsection{Finite level truncations}

If $K \subseteq H \subseteq [C_d]^\infty$ are subgroups and $\ell \geq 0$, then we define 
\[
    [H : K]_\ell := [\rho_\ell(H) : \rho_\ell(K)] = [H \St_\ell [C_d]^\infty : K\st_\ell [C_d]^\infty].
\]
Note that $[H : K]_\ell$ is a weakly increasing function of $\ell$ and that $[H : K] \geq [H : K]_\ell$ for all $\ell \geq 0$. 
For $K$ open, $[H:K] = \lim_{\ell\to\infty} [H:K]_\ell$, so the sequence necessarily stabilizes. 
We say the index of $K$ \emph{stabilizes at $k$} if $[H : K]_\ell = [H : K]$ for all $\ell \geq k$.
In Lemma \ref{lemma: branch index stabilization} we determine when the index stabilizes for the subgroups $\br$ and $\br^d$ in $\Mpre$.

\begin{lemma}
\label{lemma: branch index stabilization}
The table below shows the index of stabilization for $\br$ and $\br^d$ in cases $(A), (B), (C)$.
\begin{center}
    \vspace{1em}
    \begin{tabular}{|c|c|c|}
    \hline
            \rule[-.5em]{0pt}{1.5em} & $\br$ & $\br^d$ \\
    \hline
       $(A)$\rule[-.5em]{0pt}{1.5em} & $n$ & $m+1$ \\
    \hline
       $(B)$\rule[-.5em]{0pt}{1.5em} & $n$ & $n+1$ \\
    \hline
       $(C)$\rule[-.5em]{0pt}{1.5em} & $n+1$ & $n+2$ \\
    \hline
    \end{tabular}
    \vspace{1em}
\end{center}

Furthermore, $[\Mpre : \br]_3 = 4$ and $[\Mpre : \br^d]_4 = 8$.
\end{lemma}

\begin{proof}
    Suppose $(A)$.
    Lemma~\ref{lemma: N quotients} implies that 
    \[
        \br= \Mpre \cap \ker(\chi_m) \cap \ker(\chi_n).
    \]
    From this and $\st_\ell [C_d]^\infty \subseteq \ker(\chi_\ell)$, it follows that $\st_\ell\Mpre \subseteq \br$ for all $\ell \geq n$.
    Therefore the index of $\br$ stabilizes at $n$.
    The proof of Lemma \ref{lemma: N^d quotient} shows that $\Mpre/\br^d$ is generated by $\s$ and $b_{m+1}$ in case $(A)$.
    These elements generate a subgroup of order $d^{d+1} = [\Mpre : \br^d]$ at each level $\ell \geq m + 1$.
    Therefore $\br^d$ stabilizes at $m + 1$.
    
    Suppose $(B)$.
    Lemma~\ref{lemma: N quotients} implies that $[\Mpre : \br] = d^3$ and $\br = \ker\psi_B \cap \Mpre$.
    The formulas defining $\psi_B$ only depend on the level $n$ truncation of $\Mpre$, hence $\st_n\Mpre \subseteq \br$.
    Therefore $\br$ stabilizes at level $n$ and $\br^d$ stabilizes at level $n + 1$.

    Finally, suppose $(C)$.
    Lemma~\ref{lemma: N quotients} implies that $[\Mpre : \br] = 8$.
    In this case we do not have a description of $\br$ as the kernel of a homomorphism depending on the $\chi_i$ characters.
    We check via computer calculation (see Appendix~\ref{calculation: gap code}) that $[\Mpre : \br]_3 = 4$ and $[\Mpre : \br^d]_4 = 8$, while $[\Mpre : \br]_4 = 8 = [\Mpre : \br]$ and $[\Mpre : \br^d]_5 = 16 = [\Mpre : \br^d]$.
    Therefore the index of $\br$ stabilizes at level 4 and the index of $\br^d$ stabilizes at level 5.
\end{proof}

Recall that for a subgroup $H \subseteq [C_d]^\infty$ we write $\ord_\ell(H)$ to denote the order of $\rho_\ell(H) \subseteq [C_d]^\ell$.
Note that $\log_d|[C_d]^\ell| = {[\ell]_d}$ where
\[
    [\ell]_d := \frac{d^\ell-1}{d - 1} = 1 + d + \ldots + d^{\ell-1}.
\]
We will make several uses of the identity $[\ell + k]_d = d^k[\ell]_d + [k]_d$ in the proof of the next proposition.

\begin{proposition}
\label{prop: finite level order}
    Let $\ell \geq 0$,
    \begin{enumerate}
        \item In the periodic case, let $q_\ell$ and $r_\ell$ be the unique integers such that $\ell = q_\ell n + r_\ell$ and $0 \leq r_\ell < n$. 
        Then,
        \[
            \log_d\ord_\ell(\Mper) = [\ell]_d - d^{r_\ell}[q_\ell]_{d^n} + q_\ell.
        \]
    
        \item In the preperiodic case except $(D)$, let $\delta := \log_d[\Mpre : \br]$ and $\ve := \log_d[\Mpre : \br^d]$.
        \begin{enumerate}
            \item If $(A)$ or $(B)$ and $\ell \geq n$, then,
            \[
                \log_d\ord_\ell(\Mpre) =
                \begin{cases}
                    [\ell]_d & \text{if $\ell \leq n$}\\
                    [\ell]_d + (\ve - 1)[\ell-n]_d - \delta[\ell-n+1]_d + \delta & \text{if $\ell > n$.}
                \end{cases}
            \]
            \item If $(C)$, then
            \[
                \log_2\ord_\ell(\Mpre) =
                \begin{cases}
                    2^\ell - 1 & \text{if $\ell \leq 3$,}\\
                    13 & \text{if $\ell = 4$,}\\
                    11\cdot 2^{\ell-4} + 2 & \text{if $\ell > 4$.}
                \end{cases}
            \]
            \item If $(D)$, then
            \[
                \log_d\ord_\ell(\Mpre) = \ell + 1.
            \]
        \end{enumerate}
    \end{enumerate}
\end{proposition}
\begin{proof}

    (1) First consider the periodic case.
    Let $\alpha_\ell := \log_d \ord_\ell(\Mper)$ and let $\gamma_\ell := \log_d\ord_\ell(\wh{\Mper}_n)$.
    Lemma~\ref{lemma: A subgroup properties} implies that $\Mper = \lprof a_1 \rprof \wh{\Mper}_n^d$ and $\Mper = \lprof a_n \rprof \wh{\Mper}_n$.
    Thus by Proposition~\ref{prop: order of generators} we have
    \begin{align*}
        \alpha_\ell &= \gamma_{\ell} + \lfloor\tfrac{\ell-n}{n}\rfloor + 1
        = \gamma_\ell + \lfloor\tfrac{\ell}{n}\rfloor,\\
        \alpha_\ell &= d\gamma_{\ell-1} + \lfloor\tfrac{\ell-1}{n}\rfloor + 1.
    \end{align*}
    Eliminating $\alpha_\ell$ gives a recurrence for $\gamma_\ell$ with initial value $\gamma_0 = 0$ and
    \[
        \gamma_\ell = d\gamma_{\ell-1} + \lfloor\tfrac{\ell-1}{n}\rfloor - \lfloor\tfrac{\ell-n}{n}\rfloor  = 
        \begin{cases}
            d\gamma_{\ell-1} + 1 & \text{if }n \nmid \ell,\\
            d\gamma_{\ell-1} & \text{if }n \mid \ell.\\
        \end{cases}
    \]
    Define
    $
        \gamma_\ell' := [\ell]_d - d^{r_\ell}[q_\ell]_{d^n}.
    $
    We show that $\gamma_\ell'$ satisfies the same recurrence as $\gamma_\ell$.
    Note that $\gamma_0' = 0$, so their initial conditions agree.
    Now suppose $\ell > 0$.
    If $n\nmid \ell$, then $\ell - 1 = q_\ell n + r_\ell - 1$ where $0 \leq r_\ell - 1 < n$.
    Hence $q_{\ell-1} = q_\ell$ and $r_{\ell-1} = r_\ell - 1$.
    Thus
    \begin{align*}
        \gamma_\ell' &= [\ell]_d - d^{r_\ell}[q_\ell]_{d^n}\\
        &= d[\ell-1]_d - d^{r_\ell}[q_\ell]_{d^n} + 1\\
        &= d([\ell-1]_d - d^{r_\ell-1}[q_\ell]_{d^n}) + 1\\
        &= d([\ell-1]_d - d^{r_{\ell-1}}[q_{\ell-1}]_{d^n}) + 1\\
        &= d\gamma_{\ell-1}' + 1.
    \end{align*}
    \indent Otherwise, $n \mid \ell$ so $r_\ell=0$ and $\ell - 1 = (q_\ell-1)n + (n-1)$.
    Hence $q_{\ell-1} = q_\ell - 1$ and $r_{\ell-1} = n - 1$.
    Therefore,
    \begin{align*}
        \gamma_\ell' &= [\ell]_d - [q_\ell]_{d^n}\\
        &= (d[\ell-1]_d + 1) - (d^n[q_\ell-1]_{d^n} + 1)\\
        &= d([\ell-1]_d - d^{n-1}[q_\ell-1]_{d^n})\\
        &= d([\ell-1]_d - d^{r_{\ell-1}}[q_{\ell-1}]_{d^n})\\
        &= d\gamma_{\ell-1}'.
    \end{align*}
    Therefore, $\gamma_\ell = \gamma_\ell'$ for all $\ell \geq 0$.
    Note that $q_\ell = \lfloor\tfrac{\ell}{n}\rfloor$.
    Hence for all $\ell \geq 0$,
    \[
        \alpha_\ell =  [\ell]_d - d^{r_\ell}[q_\ell]_{d^n} + q_\ell.
    \]

    (2) Next consider the preperiodic case.
    If $\ell \leq n$, then Lemma \ref{lemma: small truncation B} implies that $\Mpre =_\ell [C_d]^\ell$.
    Hence 
    \[
        \log_d\ord_\ell(\Mpre) = \log_d\ord_\ell([C_d]^\ell) = [\ell]_d.
    \]
    Now suppose that $\ell > n$.
    First suppose $(A)$, $(B)$, or $(C)$.
    Let $\beta_\ell := \log_d \ord_\ell(\Mpre)$.
    Proposition~\ref{prop: branch} implies that $\br^d \subseteq \Mpre$ and $[\Mpre : \br] < \infty$.
    Define
    \begin{align*}
        \gamma_\ell
        &:= \log_d\ord_\ell(\br),\\
        \delta_\ell
        &:= \log_d [\Mpre : \br]_\ell,\\
        \ve_\ell
        &:= \log_d[\Mpre : \br^d]_\ell.
    \end{align*}
    Then for all $\ell \geq 1$,
    \begin{align*}
        \beta_\ell &= \gamma_\ell + \delta_\ell\\
        \beta_\ell &= d\gamma_{\ell-1} + \ve_\ell.
    \end{align*}
    Eliminating $\beta_\ell$ gives the recursion
    \[
        \gamma_\ell = d \gamma_{\ell-1} + \ve_\ell - \delta_\ell.
    \]
    
    Suppose either $(A)$ or $(B)$.
    Lemma~\ref{lemma: branch index stabilization} implies that $\delta_\ell = \delta := [\Mpre :\br]$ for all $\ell \geq n$ and $\ve_{\ell} = \ve := [\Mpre : \br^d]$ for all $\ell \geq n + 1$.
    Hence the recursion simplifies, for all $\ell \geq n + 1$, to
    \[
        \gamma_{\ell} = d\gamma_{\ell-1} + \ve - \delta.
    \]
    From this, a straightforward induction implies for all $k \geq 1$
    \begin{equation}\label{eqn: first pass order over n'}
        \gamma_{n+k} = d^k\gamma_{n} + (\ve-\delta)[k]_d
    \end{equation}
    Since $\gamma_{n} = [n]_d - \delta$ we can further simplify this to
    \begin{align*}
        \gamma_{n+k} 
        &= d^k\gamma_{n} + (\ve-\delta)[k]_d\\
        &= d^k[n]_d - d^k\delta + (\ve - \delta)[k]_d\\
        &= [n+k]_d + (\ve - 1)[k]_d - \delta[k+1]_d.
    \end{align*}
    Setting $\ell=n+k$ and regrouping terms we have for $\ell \geq n + 1$,
    \[
        \gamma_\ell = [\ell]_d + (\ve - 1)[\ell-n]_d - \delta [\ell-n+1]_d.
    \]
    From $\beta_\ell = \gamma_\ell+\delta_\ell$ and $\delta_\ell = \delta$ when $\ell \geq n$, we conclude that
    \[
        \beta_\ell = \gamma_\ell + \delta = [\ell]_d + (\ve - 1)[\ell-n]_d - \delta [\ell-n+1]_d + \delta.
    \]

    Now suppose $(C)$; hence $d = 2$ and $n = 3$.
    Then Lemma~\ref{lemma: branch index stabilization} implies that $\delta_\ell = \delta = 3$ for $\ell \geq 4$; $\ve_\ell = \ve = 4$ for $\ell \geq 5$; $\delta_3 = 2$; and $\ve_4 = 3$.
    Note that
    \[
        \gamma_4 = 2\gamma_3 + \ve_4 - \delta_4
        = 2(\beta_3 - \delta_3) + \ve_4 - \delta_4
        = 2([3]_2 - 2) + 3 - 3
        = 10.
    \]
    Hence
    \[
        \beta_4 = \gamma_4 + \delta_4 = 10 + 3 = 13.
    \]
    If $\ell > 4$, then 
    \[
        \gamma_\ell = d\gamma_{\ell-1} + \ve - \delta,
    \]
    and a simple induction implies that for all $k \geq 1$,
    \[
        \gamma_{k + 4} = 2^k\gamma_{4} + [k]_2 = 11\cdot2^{k} - 1.
    \]
    Thus for $\ell > 4$,
    \[
        \beta_\ell = \gamma_\ell + \delta = 11\cdot 2^{\ell-4} + 2.
    \]
    
    Suppose $(D)$.
    The identity $(b_1b_2)^2 = (b_1b_2, b_2b_1)$ implies that
    \[
        \ord_\ell (b_1b_2) = 2\ord_{\ell-1}(b_1b_2).
    \]
    Since $\ord_0(b_1b_2) = 1$, it follows by induction that $\ord_\ell(b_1b_2) = 2^\ell$.
    Hence $\rho_\ell(\Mpre)$ is isomorphic to the dihedral group of order $2^{\ell+1} = d^{\ell+1}$.
    Thus $\log_d\ord_\ell(\Mpre) = \ell + 1$ in case $(D)$.
\end{proof}

Given a subgroup $H \subseteq [C_d]^\infty$, we define the \emph{Hausdorff dimension} of $H$ to be
\[
    \haus(H) := \lim_{\ell\to\infty}\frac{\log_d\ord_\ell(H)}{\log_d\ord_\ell([C_d]^\infty)} = \lim_{\ell\to\infty}\frac{\log_d\ord_\ell(H)}{[\ell]_d},
\]
provided the limit exists.

\begin{corollary}\,
\label{cor: haus}
    \begin{enumerate}
        \item In the periodic case,
        \[
            \haus(\Mper) = 1 - \frac{d - 1}{d^{n} - 1}.
        \]
        \item In the preperiodic case,
        \begin{enumerate}
            \item If $(A)$ or $(B)$ then, using the notation of Proposition~\ref{prop: finite level order} we have,
            \[
                \haus(\Mpre) = 1 + \frac{\ve -1}{d^{n}} - \frac{\delta}{d^{n-1}}.  
            \]

            \item If $(C)$, then
            \[
                \haus(\Mpre) = \frac{11}{16}.
            \]

            \item If $(D)$, then $\haus(\Mpre) = 0$.
        \end{enumerate}
    \end{enumerate}
\end{corollary}

\begin{proof}
    (1)
    Consider the periodic case.
    Using the notation from Proposition~\ref{prop: finite level order}(1), observe that
    \[
        d^{r_\ell}[q_\ell]_{d^n} = d^{\ell -n} + d^{\ell-2n} + \ldots + d^{r_\ell}.
    \]
    Hence
    \[
        \lim_{\ell\to\infty}\frac{d^{r_\ell}[q_\ell]_{d^n}}{[\ell]_d} 
        = \frac{\sum_{i=1}^\infty d^{-ni}}{\sum_{i=1}^\infty d^{-i}}
        = \frac{d^{-n}(1 - d^{-n})^{-1}}{d^{-1}(1 - d^{-1})^{-1}}
        = \frac{d - 1}{d^n - 1}.
    \]
    Therefore Proposition~\ref{prop: finite level order}(1) implies
    \[
        \haus(\Mper) = \lim_{\ell\to\infty}\frac{[\ell]_d - d^{r_\ell}[q_\ell]_{d^n} + q_\ell}{[\ell]_d}
        = 1 - \frac{d-1}{d^n-1}.
    \]

    (2) 
    Now consider the preperiodic case.\\
    (2a)
    Suppose $(A)$ or $(B)$.
    Note that for any $k \geq 0$,
    \[
        \lim_{\ell\to\infty}\frac{[\ell-k]_d}{[\ell]_d} = \frac{\sum_{i=1}^\infty d^{-k-i}}{\sum_{i=1}^\infty d^{-i}} = \frac{1}{d^k}.
    \]
    Thus Proposition~\ref{prop: finite level order}(2a) implies
    \[
        \haus(\Mpre) = \lim_{\ell\to\infty}\frac{[\ell]_d + (\ve -1)[\ell-n]_d - \delta[\ell-n+1]_d}{[\ell]_d} = 1 + \frac{\ve-1}{d^{n}} - \frac{\delta}{d^{n-1}}.
    \]

    (2b)
    Suppose $(C)$.
    Then
    \[
        \haus(\Mpre) = \lim_{\ell\to\infty}\frac{11\cdot 2^{\ell-4} + 2}{[\ell]_2} =
        \lim_{\ell\to\infty}\frac{11\cdot 2^{\ell-4} + 2}{2^\ell-1} = \frac{11}{16}.
    \]

    (2c)
    Suppose $(D)$.
    Proposition~\ref{prop: finite level order}(2b) implies that 
    \[
        \haus(\Mpre) = \lim_{\ell\to\infty}\frac{\ell+1}{[\ell]_d} = 0.\qedhere
    \]    
\end{proof}

The following table combines the results of Proposition~\ref{prop: finite level order}, Corollary~\ref{cor: haus}, and the index calculations of Lemma~\ref{lemma: N quotients} and Lemma~\ref{lemma: N^d quotient}.

\begin{center}
    \vspace{1em}
    \begin{tabular}{|c|l|l|}
    \hline
            \rule[-.5em]{0pt}{1.5em} & $\log_d\ord_\ell(\Mpre)$ & $\haus(\Mpre)$ \\
    \hline
       $(A)$\rule[-.5em]{0pt}{1.5em} & $[\ell]_d + d[\ell-n]_d - 2[\ell-n+1]_d + 2$ & $1 - \frac{1}{d^{n-1}}$ \\
    \hline
       $(B)$\rule[-.5em]{0pt}{1.5em} & $[\ell]_d + \frac{3d}{2}[\ell-n]_d - 3[\ell-n+1]_d + 3$ & $1 - \frac{3}{2d^{n-1}}$ \\
    \hline
       $(C)$\rule[-.5em]{0pt}{1.5em} & $11\cdot 2^{\ell-4} + 2$ & $\frac{11}{16}$ \\
    \hline
       $(D)$\rule[-.5em]{0pt}{1.5em} & $\ell + 1$ &  $0$\\
    \hline
    \end{tabular}
    \vspace{1em}
\end{center}

\subsection{Normalizers}
\label{section: normalizer of B}

Let $N(\Mper)$ and $N(\Mpre)$ denote the normalizers of $\Mper$ and $\Mpre$, respectively, in $[C_d]^\infty$. 
Some understanding of these normalizers is required for our analysis of the constant field extensions $\wh K_f/K$.
For example, in Proposition \ref{prop: N/B finite exponent} we show that $N(\Mpre)/\Mpre$ has finite exponent except in case $(D)$, which translates into $\wh K_f/K$ being finite in those cases.

\begin{lemma}\,
\label{lemma: N(B) properties}
    \begin{enumerate}
        \item $\st_1N(\Mper) \subseteq N(\Mper)^d$ and $\st_1N(\Mpre) \subseteq N(\Mpre)^d$,
        \item If $v \in N(\Mpre)$ and $b \in \Mpre$, then $\chi_\ell(v\inv bv) = \chi_\ell(b)$ for each $\ell \geq 1$,
        \item In cases $(B)$, $(C)$, or $(D)$, for every $v \in N(\Mpre)$, then there exists an element $u \in \Mpre$ such that $v\inv bv \equiv u\inv bu \bmod \br$ for all $b \in \Mpre$.
    \end{enumerate}
\end{lemma}

\begin{proof}    
    (1)
    Since $N(\Mper)$ is contained in $[C_d]^\infty$, it preserves the level structure on $\Mper$. 
    Therefore, $N(\Mper)$ and its subgroup $\st_1 N(\Mper)$ normalize $\st_1\Mper$. Lemma~\ref{lemma: stabilizer projections are surjective} implies that $\pi_i(\st_1\Mper) = \Mper$ for each $1 \leq i \leq d$, which in turn implies that $\st_1 N(\Mpre) \subseteq N(\Mpre)^d$.
    The same argument applies to $\Mpre$.

    (2)
    The functions $\chi_\ell$ are defined on all of $[C_d]^\infty$ and $N(\Mpre) \subseteq [C_d]^\infty$ by definition.
    Hence $\chi_\ell(v\inv b v) = \chi_\ell(b)$ for all $v \in N(\Mpre)$ and all $b \in \Mpre$.
    
    (3)
    Suppose either $(B)$, $(C)$, or $(D)$.
    Let $v \in N(\Mpre)$.
    Part (2) implies that
    \begin{align*}
        v\inv b_mv &\equiv b_m[b_m,b_n]^i \bmod \br\\
        v\inv b_nv &\equiv b_n[b_m,b_n]^j \bmod \br
    \end{align*}
    for some $i, j \in \ZZ/d\ZZ$.
    Let $u = b_m^{-j}b_n^{i} \in \Mpre$.
    A straightforward calculation using the identity
    \[
        b_mb_n \equiv b_nb_m[b_m,b_n] \bmod \br
    \]
    implies that $v\inv b_mv \equiv u\inv b_mu \bmod \br$ and $v\inv b_nv \equiv u\inv b_nu \bmod \br$, hence that
    \[
        v\inv bv \equiv u\inv bu \bmod \br
    \]
    for all $b \in \Mpre$.
\end{proof}

Lemma \ref{lemma: A normalizer image} characterizes the action of the normalizer $N(\Mper)$ on the generators $a_i$.

\begin{lemma}
\label{lemma: A normalizer image}
If $w \in N(\Mper)$, then for each $1 \leq i \leq n$ there exists an $\ve_i \in 1 + d\ZZ_d$ such that $w\inv a_i w \sim a_i^{\ve_i}$ in $\Mper$.
\end{lemma}

\begin{proof}
    By Proposition \ref{prop: inductive criteria}, it suffices to prove the following proposition for all $1 \leq i \leq n$ and $\ell\geq 0$,

    \begin{quote}
        $(*_{i,\ell})$: If $w \in N(\Mper)$, then there exists an $\ve_i \in 1 + d\ZZ_d$ and a $v_i \in \Mper$ such that
        \(
            w\inv a_iw =_\ell v_i\inv a_i^{\ve_i}v_i.
        \)
    \end{quote}

    Note that if $w \in N(\Mper)$, we may replace $w$ by $wa_1^k$ in order to assume without loss of generality that $w := (w_1,\ldots,w_d) \in \st_1\Mper \subseteq \Mper^d$, where the last containment follows from Lemma \ref{lemma: N(B) properties}.

    We proceed by induction to show that $(*_{i-1,\ell-1})$ implies $(*_{i,\ell})$ for $1 \leq i \leq n$ and $\ell\geq 1$, where the $i$ subscripts are interpreted modulo $n$.
    First observe that $(*_{i,0})$ holds trivially for each $1 \leq i \leq n$ since $a_i =_0 1$.
    Now suppose that $2 \leq i \leq n$, $\ell \geq 1$, and that $(*_{i-1,\ell-1})$ is true.
    Since $a_i = (1,\ldots,1,a_{i-1})$, we have
    \[
        w\inv a_iw = (1,\ldots,1,w_d\inv a_{i-1}w_d),
    \]
    where $w_d\inv a_{i-1}w_d \in \Mper$.
    The inductive hypothesis $(*_{i-1,\ell-1})$ implies there is some $\ve_{i} \in 1 + d\ZZ_d$ and some $v_d \in \Mper$ such that
    \[
        w_d\inv a_{i-1}w_d =_{\ell-1} v_d\inv a_{i-1}^{\ve_{i}} v_d.
    \]
    Lemma \ref{lemma: stabilizer projections are surjective} provides an element $v_i \in \st_1(\Mper)$ such that $\pi_d(v_i) = v_d$.
    Thus
    \[
        w\inv a_i w = (1,\ldots,1,w_d\inv a_{i-1}w_d) 
        =_\ell (1,\ldots,1, v_d\inv a_{i-1}^{\ve_{i}}v_d)
        = v_i\inv a_i^{\ve_i} v_i,
    \]
    hence $(*_{i,\ell})$ is true.
    
    Next suppose that $(*_{n,\ell-1})$ is true for some $\ell\geq 1$.
    Then
    \[
        w\inv a_1^dw = (w_1\inv a_nw_1,\ldots,w_d\inv a_nw_d) \in \st_1(\Mper) \subseteq \Mper^d.
    \]
    In particular, $w_1\inv a_nw_1 \in \Mper$ and $(*_{n,\ell-1})$ implies that there exists an $\ve_1 \in 1 + d\ZZ_d$ and a $v_1 \in \Mper$ such that
    \begin{equation}
    \label{eqn u1 w1}
        w_1\inv a_nw_1 =_{\ell-1} v_1\inv a_n^{\ve_1}v_1.
    \end{equation}
    On the other hand, if $w\inv a_1w := \sigma(u_1, \ldots, u_d)$, then
    \[
        w\inv a_1^dw = (u_du_{d-1}\cdots u_2u_1, u_1u_du_{d-1}\cdots u_2,\ldots, u_{d-1}\cdots u_2u_1u_d),
    \]
    which implies that 
    \begin{equation}
    \label{eqn: w1 u prod}
        w_1\inv a_nw_1 = u_du_{d-1}\cdots u_2u_1.
    \end{equation}
    Recall that $\st_1(\Mper) = \lprof a_1^d \rprof\wh{\Mper}_n^d$ by
    Lemma \ref{lemma: A subgroup properties}(3).
    Hence there are $c_i \in \wh{\Mper}_n$ for $1 \leq i \leq d$ and some $k \in \ZZ_d$ such that
    \[
        a_1^{-1}w\inv a_1w = (u_1, u_2,\ldots,u_{d-1}, a_n^{-1}u_d) = (c_1a_n^{m}, \ldots, c_da_n^{m}).
    \]
    Therefore
    \begin{equation}
        \label{eqn ell_n u_i}
        \eta_n(u_j) = 
        \begin{cases}
            m & \text{if } 1 \leq j < d,\\
            m + 1 & \text{if } j = d.,
        \end{cases}
    \end{equation}
    which implies that
    \(
        u_du_{d-1}\cdots u_2 u_1 \equiv a_n^{1 + dm} \bmod \wh\Mper_n.
    \)
    Using \eqref{eqn u1 w1} we see that
    \[
        v_1\inv a_n^{\ve_1}v_1
        =_{\ell-1} w_1\inv a_nw_1
        = u_du_{d-1}\cdots u_2u_1
        \equiv a_n^{1 + dm} \bmod \wh\Mper_n.
    \]
    Thus $a_n^{\ve_1} =_{\ell-1} a_n^{1 + dm}$.
    
    Now we let $v := y_1y_2y_3$ where
    \begin{align*}
        y_1 &:= (1, a_n^{m}, a_n^{2m},\ldots, a_n^{(d-1)m})\\
        y_2 &:= (v_1, \ldots, v_1)\\
        y_3 &:= (1, u_1\inv, (u_2u_1)\inv, \ldots, (u_{d-1}\cdots u_2u_1)\inv)
    \end{align*}
    Proposition \ref{prop: stabilizer criteria} implies that $v \in \st_1\Mper$.
    
    We finish the proof by showing $w\inv a_1 w =_\ell v\inv a_1^{1 + dm} v$ in several steps, starting by conjugating
    \[
        a_1^{1+dm} = \sigma(a_n^{m},\ldots, a_n^{m}, a_n^{1 + m}).
    \]
    by $y_1$ to get
    \begin{align*}
        y_1\inv a_1^{1+dm}y_1 &= (1, a_n^{-m}, a_n^{-2m},\ldots, a_n^{-(d-1)m})\s(a_n^{m}, \ldots, a_n^{m}, a_n^{1 + m})(1, a_n^{m}, a_n^{2m},\ldots, a_n^{(d-1)m})\\
        &= \s(a_n^{-m}, a_n^{-2m},\ldots, a_n^{-(d-1)m},1)(a_n^{m}, \ldots, a_n^{m}, a_n^{1 + m})(1, a_n^{m}, a_n^{2m},\ldots, a_n^{(d-1)m})\\
        &= \s(1,\ldots,1,a_n^{1+dm}).
    \end{align*}
    Next, conjugating by $y_2$ gives us
    \[
        (y_1y_2)\inv a_1^{1 + dm_2}(y_1y_2) = \sigma (1, \ldots, 1, v_1\inv a_n^{1 + dm}v_1)
        =_\ell \sigma(1,\ldots,1, w_1\inv a_n w_1).
    \]
    Finally conjugating by $y_3$ and using \eqref{eqn: w1 u prod} we have
    \begin{align*}
        v\inv a_1^{1 + dm}v &= (y_1y_2y_3)\inv a_1^{1 + dm}(y_1y_2y_3)\\
        &=_\ell (1, u_1, \ldots, u_{d-1}\cdots u_2u_1)\sigma(1,\ldots,1, w_1\inv a_nw_1) (1, u_1^{-1}, \ldots, (u_{d-1}\cdots u_2u_1)^{-1})\\
        &= \sigma (u_1, \ldots, u_{d-1}\cdots u_2u_1,1)(1,\ldots,1, w_1\inv a_nw_1) (1, u_1^{-1}, \ldots, (u_{d-1}\cdots u_2u_1)^{-1})\\
        &= \sigma (u_1, u_2,\ldots, u_{d-1}, (w_1\inv a_nw_1)(u_{d-1}\cdots u_2u_1)^{-1})\\
        &= \sigma (u_1, u_2,\ldots, u_{d-1}, u_d)\\
        &= w\inv a_1w.
    \end{align*}
    This completes the proof of  $(*_{1,\ell})$ and hence our induction.
\end{proof}

Lemma \ref{lemma: A normalizer image} allows us to characterize $N(\Mper)$ as a subset of $N(\Mper)^d$.

\begin{lemma}
\label{lemma: A normalizer membership criteria}
    Let $v \in N(\Mper)^d$.
    Then $v \in N(\Mper)$ if and only if $v\inv a_\infty v \in \Mper$.
\end{lemma}

\begin{proof}
    Let $v := (v_1, \ldots, v_d) \in N(\Mper)^d$.
    If $2 \leq i \leq n$, then Lemma \ref{lemma: A normalizer image} implies there is some $u_d \in \Mper$ and some $\ve_i \in 1 + d\ZZ_d$ such that
    \[
        v\inv a_i v = (1,\ldots,1,v_d\inv a_{i-1}v_d)
        = (1,\ldots,1,u_d\inv a_{i-1}^{\ve_i} u_d).
    \]
    Lemma \ref{lemma: stabilizer projections are surjective} implies there is an element $u \in \st_1\Mper$ such that $\pi_d(u) = u_d$, hence $v\inv a_i v = u\inv a_i^{\ve_i} u \in \Mper$.
    Thus the identity $a_\infty = a_1a_2\cdots a_n$ implies that $v \in N(\Mper)$ if and only if $v\inv a_\infty v \in \Mper$.
\end{proof}

Next we prove an analog of Lemma \ref{lemma: A normalizer membership criteria} for $\Mpre$.

\begin{lemma}
\label{lemma: normalizer membership criteria}
    Let $v \in N(\Mpre)^d$.
    Then $v\inv b_iv \in \Mpre$ for $i \neq 1, m + 1$.
    Furthermore, $v \in N(\Mpre)$ if and only if $v\inv b_iv \in \Mpre$ for $i = 1, m + 1$.
\end{lemma}

\begin{proof}
    In cases $(A)$ or $(B)$, Lemma~\ref{lemma: psi properties} implies that $\br = \ker\psi_A \cap \Mpre$ or $\ker\psi_B \cap \Mpre$, respectively.
    Hence $N(\Mpre) \subseteq [C_d]^\infty$ normalizes $\br$.
    In case $(C)$, Lemma~\ref{lemma: branch index stabilization} implies that $\st_4(\Mpre) \subseteq \N$.
    The group $N(\Mpre)$ normalizes $\st_4(\Mpre) = \st_4[C_d]^\infty \cap \Mpre$.
    We then check via a computer calculation (see Appendix~\ref{calculation: gap code}) that $N(\Mpre)$ normalizes $\br/\st_4\Mpre$.
    Hence $N(\Mpre)_4$ normalizes $\br$ in case $(C)$ as well.
    In case $(D)$ we have $\br = 1$, hence $N(\Mpre)$ trivially normalizes $\br$.
    Therefore in all cases $N(\Mpre)^d$ normalizes $\br^d$.
    Since $\wh{\Mpre}_{1,m+1} \subseteq \br^d$, we have $v\inv b_iv \in \br^d \subseteq \Mpre$ for all $i \neq 1,m+1$.
    Thus $v \in N(\Mpre)$ if and only if $v\inv b_iv \in \Mpre$ for $i=1,m+1$.
\end{proof}

We give a detailed analysis of the constant field extensions for unicritical PCF polynomials in Section~\ref{section: constant fields}.
In the preperiodic cases, except for case $(D)$, the constant field extensions are always finite.
Mere finiteness can be deduced from Proposition~\ref{prop: N/B finite exponent}, though we obtain far more refined control over these extensions in Section~\ref{section: power conjugators}.

\begin{proposition}
\label{prop: N/B finite exponent}
    If $(A)$, $(B)$, or $(C)$, then $N(\Mpre)/\Mpre$ has a finite exponent.
\end{proposition}

\begin{proof}
    Lemma~\ref{lemma: branch index stabilization} implies that $\st_{n'}\Mpre \subseteq \br$.
    Let $\ve$ be the exponent of $[C_d]^{n'} \cong [C_d]^\infty/\st_{n'}[C_d]^\infty$.
    We will show that $v^\ve \in \Mpre$ for all $v \in N(\Mpre)$.
    By Proposition~\ref{prop: inductive criteria}, it suffices to show that $v^{\ve} \in_\ell \Mpre$ for all $\ell \geq 0$ and all $v \in N(\Mpre)$.
    We proceed by induction.
    The case $\ell = 0$ is immediate.
    Suppose that $\ell \geq 1$ and that $v^{\ve} \in_{\ell-1} \Mpre$ for all $v \in N(\Mpre)$.
    Let $v \in N(\Mpre)$. 
    Note that $[\s, v] \in \Mpre$, hence for each $i$
    \[
        (\s^i v)^{\ve} \equiv \s^{i\ve}v^\ve \equiv v^{\ve} \bmod \Mpre.
    \]
    Thus we may suppose without loss of generality that $v = (v_1, \ldots, v_d) \in \st_1 N(\Mpre) \subseteq N(\Mpre)^d$, where the last containment is Lemma~\ref{lemma: N(B) properties}(1).
    Our inductive hypothesis implies that
    \[
        v^{\ve} = (v_1^{\ve},\ldots, v_d^{\ve}) \in_\ell \Mpre^d.
    \]
    Thus for each $i$ we have $v_i^\ve \in_{\ell-1} \Mpre$ and by the definition of $\ve$ we have $v_i^\ve \in \st_{n'}[C_d]^\infty$.
    Hence $v_i^\ve \in_{\ell-1} \st_{n'}\Mpre \subseteq \N$.
    Therefore $v^\ve \in_\ell \N^d \subseteq \Mpre$, where the last containment is by Proposition~\ref{prop: branch}.
    This completes our induction.
\end{proof}

\subsection{Odometers}

Recall the elements $a_\infty \in \Mper$ and $b_\infty \in \Mpre$ defined by
\begin{align*}
    a_\infty &= a_1a_2\cdots a_n = \s(1,\ldots,1,a_na_1\cdots a_{n-1}) = \s(1,\ldots,1,a_na_\infty a_n\inv),\\
    b_\infty &= b_1b_2\cdots b_n = \s(1,\ldots,1,b_n,1,\ldots,1,b_1\cdots b_{n-1}) = \s (1,\ldots,1, b_n,1,\ldots,1, b_\infty b_n\inv).
\end{align*}

In Section~\ref{section: constant fields} we reduce the analysis of constant field extensions $\wh K_{f,\ell}$ in $\barb f$ to the problem of deciding for which $\ve \in \ZZ_d^\times$ we have $a_\infty \sim_\ell a_\infty^\ve$ in $\Mper$ and $b_\infty \sim_\ell b_\infty^\ve$ in $\Mpre$.
The following proposition is a crucial step in this reduction.

\begin{proposition}
\label{prop: N(B) odometer criteria}
    Let $v \in [C_d]^\infty$.
    Then
    \begin{enumerate}
        \item $v\inv a_\infty v \in \Mper$ if and only if $v \in N(\Mper)$,
        \item $v\inv b_\infty v \in \Mpre$ if and only if $v \in N(\Mpre)$.
    \end{enumerate}
\end{proposition}

\begin{proof}
    (1)
    If $v \in N(\Mper)$, then clearly $v\inv a_\infty v \in \Mper$.
    For the other direction we prove by induction on $\ell \geq 0$ that if $v \in [C_d]^\infty$ is an element such that $v\inv a_\infty v \in_\ell \Mper$, then $v \in_\ell N(\Mper)$.
    The base case is trivial, so suppose $\ell \geq 1$ and that the assertion holds for $\ell - 1$.

    Let $v \in [C_d]^\infty$ be an element such that $v\inv a_\infty v \in_\ell \Mper$.
    Replacing $v$ by $va_\infty^k$ we may assume that $v =_1 1$, hence $v = (v_1, \ldots, v_d)$.
    Then $v\inv a_\infty^d v \in \Mper$ and we calculate
    \[
        v\inv a_\infty^d v = (v_1\inv a_na_\infty a_n\inv v_1,\ldots,v_d\inv a_n a_\infty a_n\inv v_d) \in_\ell \st_1\Mper \subseteq \Mper^d.
    \]
    Therefore $(v_i\inv a_n)a_\infty(v_i\inv a_n)\inv \in_{\ell-1} \Mper$ for each $i$, which by our inductive hypothesis implies that $v_i \in_{\ell-1} N(\Mper)$.
    Hence $v \in_\ell N(\Mper)^d$.
    Then Lemma \ref{lemma: A normalizer membership criteria} implies that $v \in_\ell N(\Mper)$, which completes our induction.

    (2)
    The forward direction is immediate.
    For the reverse direction it suffices by Proposition~\ref{prop: inductive criteria} to prove by induction on $\ell \geq 0$ that if $v \in [C_d]^\infty$ is an element such that $v\inv b_\infty v \in_\ell \Mpre$, then $v \in_\ell N(\Mpre)$.
    If $\ell \leq n$, then Lemma~\ref{lemma: small truncation B} implies that $\Mpre =_\ell [C_d]^\infty$ which renders both assertions immediate.
    We now assume that $\ell > n$ and that the claim holds for $\ell - 1$.

    Now suppose that $v \in [C_d]^\infty$ is an element such that $v\inv b_\infty v \in_\ell \Mpre$.
    Replacing $v$ by $b_\infty^iv$ for some $i$, we may assume without loss of generality that $v = (v_1, \ldots, v_d)$ for some $v_i \in [C_d]^\infty$.

    Observe that $b_\infty^d = (b_\infty,\ldots,b_\infty,b_nb_\infty b_n\inv,\ldots, b_nb_\infty b_n\inv)$.
    Since $v\inv b_\infty^d v \in_\ell \st_1\Mpre \subseteq \Mpre^d$, we have $v_i\inv b_\infty v_i \in_{\ell-1} \Mpre$ for $1 \leq i \leq \omega$ and $v_i\inv b_n\inv b_\infty b_n v_i \in_{\ell-1} \Mpre$ for $\omega < i \leq d$.
    Hence the inductive hypothesis implies that $v_i \in_{\ell-1} N(\Mpre)$ for each $i$.
    Thus $v \in_\ell N(\Mpre)^d$.

    Lemma~\ref{lemma: normalizer membership criteria} implies that $v\inv b_iv \in_\ell \Mpre$ for all $i \neq 1, m + 1$ and that $v \in_\ell N(\Mpre)$ if and only if $v\inv b_iv \in_\ell \Mpre$ for $i = 1, m + 1$.
    Since $b_\infty = b_1\cdots b_n$ and $v\inv b_\infty v \in_\ell \Mpre$, to prove that $v \in_\ell N(\Mpre)$, it suffices to show that $v\inv b_1v = v\inv \sigma v \in_\ell \Mpre$, or equivalently that 
    \[
        u = (u_1,\ldots, u_d) := [\s,v] \in_\ell \Mpre.
    \]
    We proceed by comparing $u$ to $u' = (u_1',\ldots, u_d') := \s\inv v\inv b_\infty v$, which is in $\st_1\Mpre \subseteq \Mpre^d$. Note that 
    \begin{equation}
    \label{eqn: u}
        u = \s\inv v\inv \s v = (\s\inv v\inv \s b_2\cdots b_n v)(v\inv (b_2\cdots b_n)\inv v)
        = u' v\inv b_\infty\inv \s v.
    \end{equation}
    Since $b_\infty\inv \s \in \st_1\Mpre \subseteq \Mpre^d$ and $v \in_\ell N(\Mpre)^d$, we conclude that $u \in_\ell \Mpre^d$.

    Observe that $b_\infty\inv \s = (1,\ldots,1,b_n\inv,1,\ldots,1,b_nb_\infty\inv)$, hence comparing coordinates in \eqref{eqn: u} we have
    \[
        u_i =
        \begin{cases}
            u_i' & i \neq \omega, d,\\
            u_\omega' v_\omega\inv b_n\inv v_\omega & i = \omega,\\
            u_d' v_d\inv b_nb_\infty\inv v_d & i = d.
        \end{cases}
    \]
    
    First consider case $(A)$.
    Since $u' = \sigma \inv (v\inv b_\infty v) \in \Mpre$, Proposition~\ref{prop: stabilizer criteria}(1) implies that $\chi_m(u_i') = \chi_n(u_{i+\omega}')$ for each $i$.
    Thus $\chi_m(u_i) = \chi_n(u_{i+\omega})$ for $i \neq d$.
    Furthermore,
    \[
        \chi_m(u_d) = \chi_m(u_d') + \chi_m(v_d\inv b_n b_\infty\inv v_d)
        = \chi_n(u_\omega') - 1
        = \chi_n(u_\omega') + \chi_n(v_\omega\inv b_n\inv v_\omega) = \chi_n(u_\omega).
    \]
    Therefore Proposition~\ref{prop: stabilizer criteria}(1) implies that $u \in_\ell \Mpre$
    This completes our induction in case $(A)$.
    
    Next suppose either $(B)$, $(C)$, or $(D)$.
    It is still the case that $u'\in\Mpre$, so we apply Proposition~\ref{prop: stabilizer criteria}(2) which says that $\tau(u_i') \equiv u_{i+\omega}' \bmod \br$ for each $i$.
    Thus $\tau(u_i) \equiv u_{i+\omega} \bmod \br$ for $i \neq \omega, d$.
    We calculate that
    \begin{align*}
        u_{\omega-1}' \cdots u_1'u_d' &= v_{\omega}\inv b_\infty b_n\inv v_d\\
        u_{d-1}'\cdots u_{\omega}' &= v_d\inv b_n v_{\omega}.
    \end{align*}
    Then $\tau(u_{\omega-1}' \cdots u_1'u_d') \equiv u_{d-1}'\cdots u_{\omega}' \bmod \br$ implies that $\tau(v_{\omega}\inv b_\infty b_n\inv v_d) \equiv v_d\inv b_n v_{\omega} \bmod \br$.
    In particular, 
    \begin{equation}
    \label{eqn: chi identity}
        \chi_m(v_\omega\inv v_d) 
        = \chi_m(v_{\omega}\inv b_\infty b_n\inv v_d) - 1
        = \chi_n(v_d\inv b_n v_\omega) - 1
        = \chi_n(v_d\inv v_\omega).
    \end{equation}
    To show that $[\s,v] \in_\ell \Mpre$, it only remains to show that $\tau(u_\omega) \equiv u_d \bmod \br$.
    Observe that
    \begin{align*}
        u_{\omega-1} \cdots u_1u_d &= u_{\omega-1}' \cdots u_1'u_d'v_d\inv b_n b_\infty\inv v_d =  v_{\omega}\inv v_d\\
        u_{d-1}\cdots u_{\omega} &= u_{d-1}'\cdots u_{\omega}' v_\omega\inv b_n\inv v_\omega =  v_d\inv v_{\omega}.
    \end{align*}
    Hence it suffices to show that 
    \[
        \tau(v_\omega\inv v_d) \equiv v_d\inv v_\omega \equiv (v_\omega\inv v_d)^{-1} \bmod \br,
    \]
    which, because $v_\omega\inv v_d \in_\ell \Mpre$, is equivalent to \eqref{eqn: chi identity}.
    Therefore $[\s,v] \in_\ell \Mpre$. 
    Hence in any case we have $v \in_\ell N(\Mpre)$.
\end{proof}

\subsection{Power conjugators}
\label{section: power conjugators}

In Lemma \ref{lemma: d-adic unit conjugate} we showed that $c_\infty \sim c_\infty^\ve$ for any $\ve \in 1 + d\ZZ_d$.
Here we introduce a subgroup of elements $\cc\ve\in [C_d]^\infty$ which realize these conjugations.
We use this subgroup and its conjugates to determine the structure of the constant field extension $\wh K_{f,\ell}$.

Let $v := (1,\ldots,1,b_n,\ldots,b_n)$ where the first $b_n$ occurs in the $(\omega+1)$th component.
Given $\ve = 1 + de \in 1 + d\ZZ_d$ where $e \in \ZZ_d$, let $\ca\ve, \cb\ve \in [C_d]^\infty$ be the elements defined recursively by
\begin{align*}
    \ca\ve &:= a_1^d(1,a_\infty^e,a_\infty^{2e},\ldots, a_\infty^{(d-1)e})(\ca\ve,\ldots,\ca\ve)a_1^{-d},\\
    \cb\ve &:= v(1,b_\infty^e, b_\infty^{2e}, \ldots, b_\infty^{(d-1)e})(\cb\ve, \ldots, \cb\ve)v\inv.
\end{align*}
Recall that
\begin{align*}
    a_\infty &= a_1a_2\cdots a_n = \s(1,\ldots,1,a_na_1\cdots a_{n-1}) = \s(1,\ldots,1,a_na_\infty a_n\inv),\\
    b_\infty &= b_1b_2\cdots b_n = \s(1,\ldots,1,b_n,1,\ldots,1,b_1\cdots b_{n-1}) = \s (1,\ldots,1, b_n,1,\ldots,1, b_\infty b_n\inv),\\
    c_\infty &= \s(1, \ldots, 1, c_\infty).
\end{align*}

\begin{lemma}
\label{lemma: c_ve}
    Let $\ve, \ve_1, \ve_2 \in 1 + d\ZZ_d$, then
    \begin{enumerate}
        \item $\ca\ve a_\infty \ca\ve\inv = a_\infty^{\ve}$ and $\cb\ve b_\infty c_{\ve}^{-1} = b_\infty^{\ve}$,
        \item $\ca{\ve_1}\ca{\ve_2} = \ca{\ve_1\ve_2}$ and $\cb{\ve_1}\cb{\ve_2} = \cb{\ve_1\ve_2}$,
        \item If $0 \leq \ell \leq \infty$, then $\ca\ve \in_\ell \Mper$ if and only if $a_\infty \sim_\ell a_\infty^\ve$ in $\Mper$, and $\cb\ve \in_\ell \Mpre$ if and only if $b_\infty \sim_\ell b_\infty^\ve$ in $\Mpre$,
        \item $\chi_1(\cb\ve) = 0$ and $\chi_\ell(\cb\ve) = \binom{d}{2}e$ for $\ell > 1$.
    \end{enumerate}
\end{lemma}

\begin{proof}
    (1)
    Given $\ve = 1 + de \in 1 + d\ZZ_d$, let $\cc\ve \in [C_d]^\infty$ be the element defined recursively by
    \[
        \cc\ve := (1,c_\infty^e, c_\infty^{2e},\ldots, c_\infty^{(d-1)e})(\cc\ve,\ldots,\cc\ve).
    \]
    We calculate
    \begin{align*}
        \cc\ve c_\infty \cc\ve\inv 
        &= (1,c_\infty^e, c_\infty^{2e},\ldots, c_\infty^{(d-1)e})\s(1,\ldots,1,\cc\ve c_\infty \cc\ve\inv)(1,c_\infty^{-e}, c_\infty^{-2e},\ldots, c_\infty^{-(d-1)e})\\
        &= \s(c_\infty^e, \ldots, c_\infty^e, (\cc\ve c_\infty \cc\ve\inv)c_\infty^{-(d-1)e}).
    \end{align*}
    On the other hand,
    \[
        c_\infty^\ve = c_\infty c_\infty^{de}  
        = \s(c_\infty^e, \ldots, c_\infty^e, c_\infty^{1 + e})
        = \s(c_\infty^e, \ldots, c_\infty^e, c_\infty^\ve c_\infty^{-(d-1)e}).
    \]
    Hence Lemma~\ref{lemma: recursive solution} implies $\cc\ve c_\infty \cc\ve\inv = c_\infty^\ve$.

    First, the periodic case.
    Let $w \in [C_d]^\infty$ be the element defined recursively by
    \[
        w = a_1^d(w,\ldots,w) = (a_nw, \ldots, a_nw).
    \]
    Note that $w$ commutes with $\s$.
    Then
    \[
        wc_\infty w\inv = \s (1,\ldots,1,a_nwc_\infty w\inv a_n\inv).
    \]
    Since $wc_\infty w\inv$ satisfies the same recurrence as $a_\infty$, we conclude that $wc_\infty w\inv = a_\infty$ by Lemma~\ref{lemma: recursive solution}.
    Observe that
    \begin{align*}
        w\cc\ve w\inv &= a_1^d(1,wc_\infty^ew\inv,\ldots, wc_\infty^{(d-1)e}w\inv)(w\cc\ve w\inv,\ldots, w\cc\ve w\inv) a_1^{-d}\\
        &= a_1^d(1,a_\infty^e,\ldots,a_\infty^{(d-1)e})(w\cc\ve w\inv,\ldots, w\cc\ve w\inv) a_1^{-d},
    \end{align*}
    hence $w\cc\ve w\inv$ satisfies the same recurrence as $\ca\ve$, implying that $w \cc\ve w\inv = \ca\ve$ by Lemma~\ref{lemma: recursive solution}.
    Therefore
    \[
        \ca\ve a_\infty \ca\ve\inv = (w\cc\ve w\inv)(w c_\infty w\inv)(w \cc\ve\inv w\inv)
        = w \cc\ve c_\infty \cc\ve\inv w\inv 
        = w c_\infty^\ve w\inv
        = a_\infty^\ve.
    \]

    For the preperiodic case, let $w' \in [C_d]^\infty$ be the element defined recursively by
    \[
        w' = v(w',\ldots, w'),
    \]
    where $v = (1,\ldots,1,b_n,\ldots,b_n)$ and the first $b_n$ occurs in the $(\omega + 1)$th component.
    Observe that
    \begin{align*}
        w'c_\infty w'^{-1} &= (1,\ldots, 1, b_n,\ldots, b_n)\s(1,\ldots, 1, w'c_\infty w'^{-1})(1,\ldots,1,b_n^{-1},\ldots,b_n^{-1}) \\
        &= \s(1, \ldots, 1, b_n, 1,\ldots, 1, w'c_\infty w'^{-1}b_n^{-1}).
    \end{align*}
    Hence $w'c_\infty w'^{-1}$ satisfies the same recurrence as $b_\infty$, implying $w'c_\infty w'^{-1} = b_\infty$.

    Observe that
    \begin{align*}
        w'\cc\ve w'^{-1} &= v(1,w'c_\infty^ew'^{-1},\ldots, w'c_\infty^{(d-1)e}w'^{-1})(w'\cc\ve w'^{-1},\ldots, w'\cc\ve w'^{-1}) v\inv\\
        &= v(1,b_\infty^e,\ldots,b_\infty^{(d-1)e})(w'\cc\ve w'^{-1},\ldots, w'\cc\ve w'^{-1}) v\inv.
    \end{align*}
    Therefore $w'\cc\ve w'^{-1} = \cb\ve$, again by Lemma~\ref{lemma: recursive solution}.
    We then calculate
    \[
        \cb\ve b_\infty \cb\ve\inv = (w'\cc\ve w'^{-1})(w' c_\infty w'^{-1})(w' \cc\ve\inv w'^{-1})
        = w' \cc\ve c_\infty \cc\ve\inv w'^{-1} 
        = w' c_\infty^\ve w'^{-1}
        = b_\infty^\ve.
    \]

    (2)
    Let $\ve_1 = 1 + de_1$ and $\ve_2 = 1 + de_2$.
    We calculate
    \begin{align*}
        \cc{\ve_1}\cc{\ve_2} &= (1,c_\infty^{e_1},\ldots,c_\infty^{(d-1)e_1})(\cc{\ve_1},\ldots,\cc{\ve_1})(1,c_\infty^{e_2},\ldots,c_\infty^{(d-1)e_2})(\cc{\ve_2},\ldots,\cc{\ve_2})\\
        &= (1,c_\infty^{e_1 + \ve_1e_2},\ldots,c_\infty^{(d-1)({e_1 + \ve_1e_2})})(\cc{\ve_1}\cc{\ve_2},\ldots,\cc{\ve_1}\cc{\ve_2}).
    \end{align*}
    Note that
    \[
        \ve_1\ve_2 = 1 + de_1 + de_2 + d^2e_1e_2 = 1 + d(e_1 + \ve_1 e_2).
    \]
    Hence $\cc{\ve_1}\cc{\ve_2}$ satisfies the same recurrence as $\cc{\ve_1\ve_2}$ and we conclude that $\cc{\ve_1}\cc{\ve_2} = \cc{\ve_1\ve_2}$.
    Conjugating by $w$ and $w'$ gives us $\ca{\ve_1}\ca{\ve_2} = \ca{\ve_1\ve_2}$ and $\cb{\ve_1}\cb{\ve_2} = \cb{\ve_1\ve_2}$.

    (3)
    The forward direction is an immediate consequence of (2).
    For the converse direction, suppose that there exists some $g \in \Mper$ such that $ga_\infty g\inv =_\ell a_\infty^\ve = \ca\ve a_\infty \ca\ve\inv$.
    Then $[g\inv \ca\ve, a_\infty] =_\ell 1$ and Proposition~\ref{prop: odometers are self-centralizing} implies that $\ca\ve \in_\ell g\lprof a_\infty \rprof \subseteq \Mper$.
    The same argument, \emph{mutatis mutandis}, implies that $\cb{\ve} \in_\ell \Mpre$ if and only if $b_\infty \sim_\ell b_\infty^\ve$ in $\Mpre$.

    (4)
    By construction we have $\chi_1(\cb\ve) = 0$.
    If $\ell > 1$, we calculate
    \[
        \chi_\ell(\cb\ve) = \sum_{i=0}^{d-1}\chi_{\ell-1}(b_\infty^{ie}\cb\ve)
        = \sum_{i=0}^{d-1} \big(ie + \chi_{\ell-1}(\cb\ve)\big)
        = \binom{d}{2}e + d\chi_{\ell-1}(\cb\ve)
        \equiv \binom{d}{2}e \bmod d.\qedhere
    \]
\end{proof}

In Theorem~\ref{theorem: outer action determined by infinity} we prove that for polynomials $f$ in the preperiodic case, there is a direct correspondence between elements in the Galois group of the constant field extension $\wh{K}_f/K$ and $\ve \in \ZZ_d^\times$ such that $\cb\ve \in \Mpre$. 
Under this correspondence, the following proposition allows us to precisely determine $\gal(\wh{K}_f/K)$. This is a substantial refinement of Proposition~\ref{prop: N/B finite exponent}.

\begin{proposition}
\label{prop: c_ve in B}
    Let $\ve \in 1 + d\ZZ_d$, let $0 \leq \ell \leq \infty$, and let
    \[
        \kappa := 
        \begin{cases}
            d^2/\gcd(d,\omega) & \text{if $(A)$ and either $m > 1$ or $d$ odd,}\\
            d^2/\gcd(d,\omega + d/2) & \text{if $(A)$ and $m = 1$ and $d$ even,}\\
            4d & \text{if $(B)$ or $(C)$}.
        \end{cases}
    \]
    Then
    \begin{enumerate}
        \item $a_\infty \sim_\ell a_\infty^\ve$ in $\Mper$ if and only if $\ve \equiv 1 \bmod d^{\lfloor (\ell-1)/n\rfloor + 1}$,
        \item $b_\infty \sim_\ell b_\infty^\ve$ in $\Mpre$ if and only if $\ell \leq n$ or $\ell > n$ and $\ve \equiv 1 \bmod \kappa$.
    \end{enumerate}
\end{proposition}

\begin{proof}
    (1)
    Lemma \ref{lemma: c_ve} implies that $a_\infty \sim_\ell a_\infty^\ve$ if and only if $\ca\ve \in_\ell \Mper$.
    A simple induction implies that $\ca\ve \in \Mper^d$; let $\ca\ve := (h_1,\ldots,h_d)$.
    Recall that Proposition \ref{prop: stabilizer criteria} implies that
    \[
        \st_1\Mper = \{(g_1,\ldots,g_d) \in \Mper^d : \eta_n(g_i) = \eta_n(g_j) \text{ for all $i,j$}\}.
    \]
    Proposition \ref{prop: order of generators} implies that $\ord_\ell(a_n) = d^{\lfloor \ell/n\rfloor}$, hence $\ca\ve \in_\ell \Mper$ if and only if $\eta_n(h_i) \equiv \eta_n(h_j) \bmod d^{\lfloor (\ell-1)/n\rfloor}$ for each $i$ and $j$.
    The definition of $\ca\ve$ implies that $\eta_n(h_i) = (i-1)e$.
    Hence in particular we must have
    \[
        e \equiv \eta_n(h_2) \equiv \eta_n(h_1) \equiv 0 \bmod d^{\lfloor (\ell-1)/n\rfloor}.
    \]
    This is also clearly sufficient.
    Therefore $a_\infty \sim_\ell a_\infty^\ve$ in $\Mper$ if and only if $\ve \equiv 1 \bmod d^{\lfloor (\ell-1)/n\rfloor + 1}$.

    (2)
    If $\ell \leq n$, then Lemma \ref{lemma: small truncation B} implies that $\Mpre =_\ell [C_d]^\infty$ and Lemma \ref{lemma: d-adic unit conjugate} implies that $b_\infty \sim_\ell b_\infty^\ve$ for all $\ve \in 1 + d\ZZ_d$.

    Suppose that $\ell > n$.
    Lemma \ref{lemma: c_ve} implies that $b_\infty \sim_\ell b_\infty^\ve$ in $\Mpre$ if and only if $\cb\ve \in_\ell \Mpre$.
    Our argument above implies that $\cb\ve \in_{\ell-1} \Mpre$, hence that $\cb\ve \in_\ell \Mpre^d$.
    Thus we may use the criteria provided by Proposition \ref{prop: stabilizer criteria} to characterize when $\cb\ve \in_\ell \st_1 \Mpre$.
    Note that these criteria are vacuous for $\ell \leq n$ but impose constraints for all $\ell > n$.
    Let us write $\cb\ve := (g_1, \ldots, g_d)$.
    The definition of $\cb\ve$ implies that for all $k \geq 1$,
    \begin{equation}
    \label{eqn: chi c_ve}
        \chi_k(g_i) = \chi_k(b_\infty^{ie}\cb\ve)
        = ie + \chi_k(\cb\ve)
        = 
        \begin{cases}
            ie & k = 1,\\
            ie + \binom{d}{2}e & k > 1.
        \end{cases}
    \end{equation}

    Suppose $(A)$.
    Proposition~\ref{prop: stabilizer criteria}(1) implies that $\cb\ve \in_\ell \Mpre$ if and only if $\chi_m(g_i) = \chi_n(g_{i+\omega})$ for each $i$.
    By \eqref{eqn: chi c_ve}, this is equivalent to
    \[
        ie + \chi_m(\cb\ve) \equiv (i + \omega)e + \chi_n(\cb\ve) \bmod d,
    \]
    which reduces to
    \[
        \omega e  \equiv \chi_m(\cb\ve) - \chi_n(\cb\ve) \bmod d,
    \]
    or equivalently
    \[
        d \text{ divides } 
        \begin{cases}
            \omega e & m >1,\\
            \omega e + \binom{d}{2}e & m = 1.
        \end{cases}
    \]
    If $m > 1$, this is equivalent to $d/\gcd(d,\omega)$ dividing $e$, hence $\ve \equiv 1 \bmod d^2/\gcd(d,\omega)$.
    If $m = 1$, then $\cb\ve \in_\ell \Mpre$ is equivalent to $d/\gcd\big(d,\omega + \binom{d}{2}\big)$ dividing $e$.
    Note that
    \[
        \binom{d}{2} \equiv
        \begin{cases}
            0 \bmod d & \text{if $d$ odd,}\\
            \frac{d}{2} \bmod d & \text{if $d$ even.}
        \end{cases}
    \]
    Hence 
    \[
        \gcd\Big(d, \omega + \binom{d}{2}\Big) = 
        \begin{cases}
            \gcd(d,\omega) & \text{if $d$ odd},\\
            \gcd(d,\omega + d/2) & \text{if $d$ even}.
        \end{cases}
    \]
    Therefore if $m = 1$ and $d$ odd, then $\cb\ve \in_\ell \Mpre$ if and only if $\ve \equiv 1 \bmod d^2/\gcd(d,\omega)$.
    If $m = 1$ and $d$ is even, then $\cb\ve \in_\ell \Mpre$ if and only if $\ve \equiv 1 \bmod d^2/\gcd(d,\omega + d/2)$.

    Suppose $(B)$ or $(C)$.
    Proposition~\ref{prop: stabilizer criteria}(2) implies that $\cb\ve \in_\ell \Mpre$ if and only if 
    \[
        \tau(g_i) \equiv g_{i+\omega} \bmod \br
    \]
    for each $0 \leq i < \omega$, or equivalently
    \begin{equation}
    \label{eqn: tau c_ve}
        \tau(b_\infty^{ie}\cb\ve) \equiv b_nb_\infty^{(i+\omega)e}\cb\ve b_n^{-1} \bmod \br.
    \end{equation}
    Since $b_\infty \equiv b_mb_n \bmod \br$, it follows by a simple induction that for all $k \in \ZZ$, 
    \[
        b_\infty^k \equiv b_m^kb_n^k[b_n,b_m]^{\binom{k}{2}} \bmod \br.
    \]
    Suppose that $a, b, c\in \ZZ/d\ZZ$ satisfy
    \[
        \cb\ve \equiv b_m^ab_n^b [b_n,b_m]^c \bmod \br.
    \]
    We calculate
    \begin{align*}
        \tau(b_\infty^{ie}\cb\ve) &\equiv (b_n^{ie}b_m^{ie}[b_n,b_m]^{-\binom{ie}{2}})(b_n^ab_m^b[b_n,b_m]^{-c}) \bmod \br\\
        &\equiv b_m^{ie + b}b_n^{ie + a}[b_n,b_m]^{-\binom{ie}{2}-c + (ie)^2 + bie + ab} \bmod \br,\\
        b_nb_\infty^{(i+\omega)e}\cb\ve b_n^{-1} &\equiv b_n(b_m^{(i+\omega)e}b_n^{(i+\omega)e}[b_n,b_m]^{\binom{(i+\omega)e}{2}})(b_m^{a}b_n^{b}[b_n,b_m]^{c})b_n^{-1} \bmod \br\\
        &\equiv b_m^{(i + \omega)e + a}b_n^{(i + \omega)e + b}[b_n,b_m]^{\binom{(i+\omega)e}{2} + c + a(i+\omega)e + (i + \omega)e + a} \bmod \br.
    \end{align*}
    Comparing exponents of $b_n$ and $b_m$ we conclude that $\omega e \equiv 0 \bmod d$.
    Since $\omega = d/2$ in cases $(B)$ and $(C)$, this is equivalent to $e$ being even.
    Comparing exponents of $[b_n,b_m]$ we have
    \[
        -\binom{ie}{2}-c + (ie)^2 + bie + ab 
        \equiv \binom{(i+\omega)e}{2} + c + a(i+\omega)e + (i + \omega)e + a \bmod d,
    \]
    which simplifies to
    \begin{equation}
    \label{eqn: simplified commutator exponent}
        \frac{de(de - 2)}{8} + 2c - ab + a +(a - b)ie \equiv 0 \bmod d.
    \end{equation}
    In case $(C)$ we have $d = 2$, $m = 2$, and $i = 0$.
    Thus Lemma~\ref{lemma: c_ve}(4) implies $a = \chi_m(\cb\ve) = e = \chi_n(\cb\ve) = b$.
    Hence \eqref{eqn: simplified commutator exponent} reduces to
    \[
        \frac{e(e-1)}{2} \equiv 0 \bmod 2.
    \]
    Since $e$ is even, this is equivalent to $4$ dividing $e$.

    Now suppose we are in case $(B)$.
    Since $m = 1$, Lemma~\ref{lemma: c_ve}(4) implies $a = \chi_m(\cb\ve) = 0$ and $b = \chi_n(\cb\ve) = \binom{d}{2}e$.
    Furthermore, Lemma~\ref{lemma: psi properties} implies that
    \[
        2c = 2\chi'(\cb\ve) 
        = 2\sum_{i=1}^d i \chi_{n-1}(g_i)
        = 2\sum_{i=1}^d i^2e + i\chi_{n-1}(\cb\ve)
        \equiv \frac{d(d+1)(2d+1)e}{3}
        \equiv \frac{(2d^2 + 1)de}{3} \bmod d.
    \]
    Since \eqref{eqn: simplified commutator exponent} holds for all $i$, we have
    \begin{align*}
        (a - b)e &\equiv 0 \bmod d,\\
        \frac{de(de - 2)}{8} + 2c - ab + a &\equiv 0 \bmod d.
    \end{align*}
    Substituting for $a, b, c$ gives us
    \begin{align}
        \binom{d}{2}e^2 &\equiv 0 \bmod d\nonumber\\
        \frac{de(de - 2)}{8} + \frac{(2d^2+1)de}{3} &\equiv 0 \bmod d.\label{eqn: second line}
    \end{align}
    Since $e$ is even, $\binom{d}{2}e \equiv 0 \bmod d$.
    If 3 does not divide $d$, then $\frac{(2d^2+1)de}{3} \equiv 0 \bmod d$.
    Let $e' = e/2$.
    Then
    \[
        0 \equiv \frac{de(de-2)}{8} \equiv \frac{de'(de'-1)}{2} \bmod d.
    \]
    Since $d$ is even, this reduces to $\frac{de'}{2} \equiv 0 \bmod d$, which is equivalent to 4 dividing $e$.
    If 3 divides $d$, then \eqref{eqn: second line} is equivalent to
    \[
        6k = 3e'(de'-1) + 2e
    \]
    for some integer $k$.
    This implies 3 divides $e$ and 4 divides $e$.
    Hence in either case, \eqref{eqn: second line} is equivalent to 4 dividing $e$ in case $(B)$.
    Therefore $\cb\ve \in_\ell \Mpre$ if and only if 4 divides $e$, or equivalently, if and only if $\ve \equiv 1 \bmod 4d$.
\end{proof}

\section{Constant field extensions and the outer action}
\label{section: constant fields}

In this section we study the constant field extensions $\wh{K}_{f,\ell}/K$ in iterated pre-image extensions $K(f^{-\infty}(t))/K(t)$.
We show that for any polynomial $f$ with degree $d$ coprime to $\cchar K$, the constant field extension is contained in the pro-$d$ cyclotomic extension $K(\zeta_{d^\infty})/K$ and is completely encoded within the structure of $\barb f$.
Then we turn to the case of unicritical polynomials where we can leverage our analysis of $\barb f$ to precisely determine $\wh{K}_{f,\ell}/K$ for all $\ell \geq 1$. 
In particular, show that $\Gal(\wh K_f/K)$ has a faithful outer action on $\bArb f$ which factors through the cyclotomic character.
Applying the results of Section~\ref{section: further properties} tells us how much farther the action factors.

\subsection{Preliminary bounds}

First we establish a general upper bound on the constant field extensions of a polynomial.

\begin{proposition}\label{prop: cyclotomic bound}
    Let $f(x) \in K[x]$ be a polynomial of degree $d \geq 2$, and assume that $d$ is coprime to $\cchar K$. Then the constant field extension $\wh K_{f,\ell}$ is contained in $K(\zeta_{d^\ell})$, and hence $\wh K_f \subseteq K(\zeta_{d^\infty})$.
\end{proposition}

\begin{proof}
    Our assumption that $f$ is a polynomial and $d$ is coprime to $\cchar K$ implies that $f$ has a totally tamely ramified fixed point at $\infty$.
    Therefore
    \[
        K\lp 1/t\rp(f^{-\ell}(t)) = K(\zeta_{d^\ell})\lp 1/t^{1/d^\ell}\rp.
    \]
    Hence
    \[
        \wh{K}_{f,\ell} = K\sep \cap K(f^{-\ell}(t)) \subseteq K\sep \cap K\lp 1/t\rp(f^{-\ell}(t)) = K(\zeta_{d^\ell}).\qedhere
    \]
\end{proof}

This upper bound on the constant field extension need not be sharp but a result of Hamblen and Jones~\cite{hamblen_jones} shows that it is when $f(x) = ax^d + b$ has a periodic critical point.

\begin{corollary}
\label{cor: periodic const field}
    Let $f(x) = ax^d + b \in K[x]$ where $d$ is coprime to $\cchar K$ and $0$ is periodic under $f$, then $\wh K_f = K(\zeta_{d^\infty})$.
\end{corollary}
    
\begin{proof}
    Hamblen and Jones \cite[Thm. 2.1]{hamblen_jones} prove, in this situation, that $K(\zeta_{d^\infty}) \subseteq \wh{K}_f$.
    Therefore Proposition~\ref{prop: cyclotomic bound} implies that $K(\zeta_{d^\infty}) = \wh{K}_f$.
\end{proof}

The situation when 0 is strictly preperiodic under $f(x) = ax^d + b$ is quite different.
Note that $K(\zeta_d) \subseteq \wh{K}_f$ for any unicritical polynomial $f(x) = ax^d + b$ since
\begin{equation}
\label{eqn: constant field simple lower bound}
    K(f^{-1}(t)) = K(\zeta_d, \big(\tfrac{t - b}{a}\big)^{1/d}).
\end{equation}

In the preperiodic case we can prove the following coarse finiteness result.
\begin{proposition}
\label{prop: preper finite const}
    If $f(x) = ax^d + b \in K[x]$ where $d$ is coprime to the characteristic of $K$ and 0 is strictly preperiodic under $f$, then $\wh K_f$ is a finite extension of $K(\zeta_d)$, hence of $K$, unless $d = 2$ and $f$ is conjugate over $K$ to $x^2 - 2$, the degree 2 Chebyshev polynomial.
\end{proposition}

\begin{proof}
    Recall that
    \(
        \gal(\wh{K}_f/K) \cong \arb f/\barb f.
    \)
    Corollary~\ref{corollary: final IMG calculation} implies that there is a $w \in [C_d]^\infty$ such that $w\barb f w\inv = \Mpre$ for the appropriate values of $d, m, n, \omega$.
    Since $\barb f$ is a normal subgroup of $\arb f$, it follows that $w \arb f w\inv \subseteq N(\Mpre)$ where $N(\Mpre)$ is the normalizer of $\Mpre$ in $[C_d]^\infty$.
    Thus we have an injective homomorphism $\gal(\wh{K}_f/K) \hookrightarrow N(\Mpre)/\Mpre$.

    On the other hand, Proposition~\ref{prop: cyclotomic bound} implies that $\gal(\wh{K}_f/K)$ is a quotient of $\gal(K(\zeta_{d^\infty})/K)$, which is isomorphic to a finitely generated subgroup of $\ZZ_d^\times$.
    Proposition~\ref{prop: N/B finite exponent} implies that $N(\Mpre)/\Mpre$ has a finite exponent unless $(d,m,n,\omega) = (2,1,2,1)$.
    A finitely generated abelian group with a finite exponent is finite, therefore $\gal(\wh K_f/K)$ is finite in all but one case.

    If $(d,m,n,\omega) = (2,1,2,1)$, then $f$ is conjugate over $K$ to a polynomial of the form $f(x) = x^2 + c$ such that $f^2(c) = f(c)$ and $c \neq 0$.
    The only such polynomial is $x^2 - 2$.
\end{proof}

One may extract an explicit degree bound on $\wh K_f/K$ from the argument in Proposition \ref{prop: preper finite const}, however in Theorem~\ref{theorem: intro constant field} we refine this result to an exact determination of $\wh K_{f,\ell}/K$ for all $\ell \geq 1$.

The finiteness of the constant field extensions shown in Proposition~\ref{prop: preper finite const} stems from the fact that $\barb f$ is a regular branch group in each preperiodic case except case $(D)$ (see Proposition~\ref{prop: branch}).
In case $(D)$ and the periodic case, the group $\barb f$ is not a regular branch group and the constant field extension is all or nearly all of the cyclotomic extension $K(\zeta_{d^\infty})/K$.
This suggests the following question:

\begin{question}
    Does $\barb f$ being a regular branch group imply that $\wh{K}_f/K$ is a finite extension for any polynomial $f(x)$? For any rational function?
\end{question}

As the proof of Proposition~\ref{prop: preper finite const} shows, for an affirmative answer it would suffice to prove that $\barb f$ regular branch implies $N(\barb f)/\barb f$ has finite exponent.

\subsection{Outer actions and sharper bounds}
\label{sec: explicit constant fields}

In this section, we show that the structure of the constant field extension can be interpreted group-theoretically, in terms of an outer Galois action. This allows us to refine the results of the preceding subsection and explicitly calculate the constant field extensions.

In the preperiodic cases except for $(D)$, Proposition~\ref{prop: preper finite const} implies that $K(\zeta_d) \subseteq \wh{K}_{f,\ell} \subseteq K(\zeta_{d^\ell})$.
Combining Proposition~\ref{prop: c_ve in B} with the branch cycle lemma (Lemma~\ref{lemma: BCL}) we give a precise determination of $\wh{K}_{f,\ell}$.

Let $P \subseteq \PP_{K\sep}^1$ be a Galois-stable set of points and let $K(t)_P/K(t)$ denote the maximal tamely ramified extension of $K\sep(t)$ unramified outside of $P$.
As discussed in the proof of Lemma~\ref{lemma: inertia generators}, the Galois group of $K\sep(t)_P/K\sep(t)$ is topologically generated by inertia generators over the points in $P$.
The branch cycle lemma (first appearing in Fried~\cite{FriedBCL}) constrains how $\gal(\wh K_{f,\ell}/K)$ acts on these inertia generators.
Our proof of the branch cycle lemma generalizes one appearing in Malle and Matzat~\cite[Thm. 2.6]{MalleMatzat} for $\QQ$. 
In characteristic zero, their proof requires no changes; tameness allows their argument to be extended into positive characteristic.

\begin{lemma}[Branch Cycle Lemma]
\label{lemma: BCL}
    Let $K$ be a field and let $P \subseteq \PP_{K\sep}^1$ be a Galois-stable set.
    Let $\chi_\cyc : \gal(K\sep/K) \to \wh \ZZ^\times$ be the cyclotomic character defined by
    \[
        \gamma(\zeta) = \zeta^{\chi_\cyc(\gamma)}
    \]
    for all $\gamma \in \gal(K\sep/K)$ and all roots of unity $\zeta \in K\sep$.
    Let $b \in P$ and let $\gamma_b$ be a topological generator for an inertia group over $b$ in $\gal(K\sep(t)_P/K\sep(t))$.
    If $\tau \in \gal(K\sep/K)$ and $\tilde\tau \in \gal(K(t)_P/K(t))$ is any lift, then
    \[
        \tilde\tau \gamma_b \tilde\tau\inv \sim \gamma_{\tau(b)}^{\chi_\cyc(\tau)},
    \]
    where the conjugacy takes place in $\gal(K\sep(t)_P/K\sep(t))$.
\end{lemma}

\begin{proof}
    Let $\tau \in \gal(K\sep/K)$ and let $b \in P$.
    After a change of coordinates over $K$, we may assume for simplicity that $b$ and $\tau(b)$ are both finite.
    Embedding $K\sep(t)_P$ into a separable closure of $K\sep\lp t - b\rp$ we have (by our tameness hypothesis) that $K\sep(t)_PK\sep\lp t - b\rp/K\sep \lp t - b\rp$ is generated by elements $(t - b)^{1/n}$ for all $n$ coprime to $\cchar K$ and that the Galois group of this extension is topologically generated by a lift of $\gamma_b$.
    We may assume that these elements are compatible in the sense that $((t - b)^{1/mn})^m = (t - b)^{1/n}$.
    Let $(\zeta_n)$ be primitive $n$th roots of unity in $K\sep$ such that
    \[
        \gamma_b (t - b)^{1/n} = \zeta_n (t - b)^{1/n}.
    \]
    Note that $\gamma_b$ is determined by how it acts on the elements $(t - b)^{1/n}$.
    
    The element $\tilde\tau \gamma_b \tilde\tau\inv \in \gal(K\sep(t)_P/K\sep(t))$ generates an inertia group over $\tau(b)$.
    The element $\tilde\tau$ may be extended to a $K\sep$-isomorphism of completions such that
    \[
        (t - \tau(b))^{1/n} := \tilde\tau((t - b)^{1/n}).
    \]
    Replacing $\gamma_{\tau(b)}$ by a conjugate in $\gal(K\sep(t)_P/K\sep(t))$, we may assume that
    \[
        \gamma_{\tau(b)}(t - \tau(b))^{1/n} = \zeta_n (t - \tau(b))
    \]
    for all $n$. 
    We calculate
    \begin{align*}
        \tilde\tau \gamma_b \tilde \tau\inv(t - \tau(b))^{1/n} &= \tilde \tau (t - b)^{1/n}\\
        &= \tilde\tau \zeta_n (t - b)^{1/n}\\
        &= \zeta_n^{\chi_\cyc(\tau)} (t - \tau(b))^{1/n}\\
        &= \gamma_{\tau(b)}^{\chi_\cyc(\tau)} (t - \tau(b))^{1/n}.
    \end{align*}
    Since $\tilde\tau \gamma_b \tilde\tau\inv$ and $\gamma_{\tau(b)}^{\chi_\cyc(\tau)}$ are determined by how they acts on $(t - \tau(b))^{1/n}$ for all $n$ coprime to $\cchar K$, we conclude that $\tilde \tau \gamma_b \tilde \tau\inv \sim \gamma_{\tau(b)}^{\chi_\cyc(\tau)}$ in $\gal(K\sep(t)_P/K\sep(t))$.
\end{proof}

\begin{theorem}
\label{theorem: outer action determined by infinity}
    Let $K$ be a field, let $0 \leq \ell \leq \infty$, and let $f(x) \in K[x]$ be a polynomial with degree $d$ coprime to $\cchar K$.
    If $\tau \in \gal(K\sep/K)$, then $\tau$ fixes $\wh K_{f,\ell}$ if and only if $\gamma_\infty \sim_\ell \gamma_\infty^{\chi_{\cyc}(\tau)}$ in $\bpimg f$.
\end{theorem}

\begin{proof}
    Let $P_f \subseteq \PP_{K\sep}^1$ be the post-critical set of $f$.
    Then by assumption $P_f$ is a Galois-stable set of points.
    Hence we have the following commutative diagram with exact rows:
    \begin{equation}
            \label{diagram: two row outer}
            \begin{tikzcd}
                0 & \gal(K\sep(t)_{P_f}/K\sep(t)) & \gal(K(t)_{P_f}/K(t)) & \gal(K\sep/K) & 0 \\
                0 & \rho_\ell({\bpfImg f}) & \rho_\ell({\pfImg f}) & {\gal(\wh K_{f,\ell}/K)} & 0
                \arrow[from=1-1, to=1-2]
                \arrow[from=1-2, to=1-3]
                \arrow["{\bar\rho}", two heads, from=1-2, to=2-2]
                \arrow[from=1-3, to=1-4]
                \arrow["\rho", two heads, from=1-3, to=2-3]
                \arrow[from=1-4, to=1-5]
                \arrow["{\wh\rho}", two heads, from=1-4, to=2-4]
                \arrow[from=2-1, to=2-2]
                \arrow[from=2-2, to=2-3]
                \arrow[from=2-3, to=2-4]
                \arrow[from=2-4, to=2-5]
            \end{tikzcd}
        \end{equation}

    Let $\tau \in \gal(K\sep/K)$.
    First suppose that $\tau$ fixes $\wh{K}_{f,\ell}$, hence belongs to the kernel of $\wh\rho$.
    Thus any lift $\tilde\tau$ of $\tau$ acts trivially on $\rho_\ell(\bpimg f)$, which means that $\tilde\tau\gamma_\infty\tilde\tau\inv \sim_\ell \gamma_\infty$ in $\bpimg f$.
    On the other hand, the branch cycle lemma implies that $\tilde\tau\gamma_\infty\tilde\tau\inv \sim \gamma_\infty^{\chi_\cyc(\tau)}$ in $\bpimg f$.
    Therefore $\gamma_\infty \sim_\ell \gamma_\infty^{\chi_\cyc(\tau)}$ in $\bpimg f$.

    Next suppose that $\gamma_\infty \sim_\ell \gamma_\infty^{\chi_\cyc(\tau)}$ in $\bpimg f$.
    If $\tilde\tau$ is any lift of $\tau$, then the branch cycle lemma implies that $\tilde\tau\gamma_\infty\tilde\tau\inv \sim \gamma_\infty^{\chi_\cyc(\tau)}\sim_\ell \gamma_\infty$ in $\bpimg f$.
    Let $\delta \in \bpimg f$ be an element such that $\tilde\tau\gamma_\infty\tilde\tau\inv =_\ell \delta \gamma_\infty \delta\inv$.
    Then $[\delta\inv \tilde\tau, \gamma_\infty] =_\ell 1$.
    Since $f$ is a polynomial, Proposition~\ref{prop: odometer labeling} and Proposition \ref{prop: odometers are self-centralizing} together imply that $\tilde\tau \in_\ell \delta\lprof\gamma_\infty\rprof \subseteq \bpimg f$.
    The exactness of the bottom row of \eqref{diagram: two row outer} implies that $\wh\rho(\tau) = 1$. 
    Hence $\tau$ fixes $\wh{K}_{f,\ell}$.
\end{proof}

Theorem~\ref{theorem: outer action determined by infinity} shows that for polynomials, the constant field extension is entirely encoded in the structure of the geometric profinite iterated monodromy group.
As a first illustration of this result we determine the constant field extension for post-critically infinite polynomials.

\begin{proposition} [PCI Constant Field]
\label{prop: PCI constant field}
    Let $f(x) = ax^d + b$ be post-critically infinite.
    Then $\wh{K}_{f,\ell} = K(\zeta_d)$ for $1 \leq \ell \leq \infty$.
\end{proposition}

\begin{proof}
    In Proposition~\ref{prop: PCI arb} we showed that $\barb f = [C_d]^\infty$ when $f$ is post-critically infinite.
    Proposition~\ref{prop: odometer labeling} implies that $\gamma_\infty$ is a strict odometer and Lemma~\ref{lemma: d-adic unit conjugate} implies that $\gamma_\infty \sim_\ell \gamma_\infty^\ve$ in $[C_d]^\infty$ if and only if $\ve \equiv 1 \bmod d$.
    Therefore Theorem~\ref{theorem: outer action determined by infinity} implies that $\wh{K}_{f,\ell}$ is the fixed field of $K(\zeta_{d^\infty})$ corresponding to the subgroup of $\ZZ_d^\times$ generated by all $\ve \equiv 1 \bmod d$, which is precisely $K(\zeta_d)$.
\end{proof}

Next we use Theorem~\ref{theorem: outer action determined by infinity} to determine $\wh K_{f,\ell}$ for all $1 \leq \ell \leq \infty$ in the periodic case; this also provides an alternative proof of Corollary~\ref{cor: periodic const field}.
Let $\chi_{\cyc,d} : \gal(K(\zeta_{d^\infty})/K) \to \ZZ_d^\times$ denote the pro-$d$ cyclotomic character of $K$.
Observe that
\[
    \gal(K(\zeta_{d^\infty})/K(\zeta_{\kappa})) = \{\tau \in \gal(K(\zeta_{d^\infty})/K) : \chi_{\cyc,d}(\tau) \equiv 1 \bmod \kappa\}.
\]

\begin{proposition}[Periodic Constant Field]
\label{prop: periodic case constant field}
    Let $f(x) = ax^d + b$ and suppose that $0$ is periodic under $f$ with period $n$.
    Then $\wh K_{f,\ell} = K(\zeta_{d^{\lfloor \ell/n\rfloor + 1}})$ for $\ell \geq 1$.
    In particular, $\wh K_f = K(\zeta_{d^\infty})$.
\end{proposition}

\begin{proof}
    Recall that $\barb f = \lprof c_1, \ldots, c_n\rprof$ and that $c_\infty = c_1c_2\cdots c_n$ is the image of $\gamma_\infty$ in $\barb f$.
    Let $\Mper = \Mper(d,n)$.    
    Corollary~\ref{corollary: final IMG calculation} implies that there exists a $w \in [C_d]^\infty$ and elements $u_i \in \Mper$ such that
    \[
        wc_iw\inv = u_i a_i u_i\inv
    \]
    for each $1 \leq i \leq n$.
    Therefore
    \[
        wc_\infty w\inv = (u_1a_1u_1\inv)\cdots (u_n a_n u_n\inv).
    \]
    Note that $wc_\infty w\inv$ is a strict odometer in $\Mper$, hence there exists a $v \in [C_d]^\infty$ such that $v\inv a_\infty v = w c_\infty w\inv \in \Mper$.
    Thus Lemma \ref{prop: N(B) odometer criteria} implies that $v \in N(\Mper)$.
    Therefore $c_\infty \sim_\ell c_\infty^\ve$ in $\barb f$ if and only if $a_\infty \sim_\ell a_\infty^\ve$ in $\Mper$.
    Proposition \ref{prop: c_ve in B} implies that $a_\infty \sim_\ell a_\infty^\ve$ in $\Mper$ if and only if $\ve \equiv 1 \bmod d^{\lfloor (\ell-1)/n\rfloor + 1}$.
    Hence Theorem \ref{theorem: outer action determined by infinity} implies that $\wh K_{f,\ell} = K(\zeta_{d^{\lfloor (\ell-1)/n\rfloor + 1}})$.
\end{proof}

Finally we determine $\wh K_{f,\ell}$ in the preperiodic case.

\begin{theorem}[Preperiodic Constant Field]
\label{theorem: preper const field}
    Let $f(x) = ax^d + b \in K[x]$ where $d$ is coprime to the characteristic of $K$ and 0 is strictly preperiodic.
    Let $m < n$ be the smallest integers such that $f^m(b) = f^n(b)$ and let $1 < \omega < d$ be such that $f^n(0) = \zeta_d^\omega f^m(0)$.
    If $(d,m,n,\omega) \neq (2,1,2,1)$ (case $(D)$), then
    \[
        K(\zeta_d) \subseteq \wh K_f \subseteq K(\zeta_{2d^2}).
    \]
    More precisely, if $1 \leq \ell \leq \infty$, then for $\ell \leq n$ we have $\wh K_{f,\ell} = K(\zeta_d)$ and for $\ell > n$ we have
    \[
        \wh K_{f,\ell} =
        \begin{cases}
            K(\zeta_{d^2/\gcd(d,\omega)}) & \text{if $(A)$ and either $m > 1$ or $d$ odd,}\\
            K(\zeta_{d^2/\gcd(d,\omega+d/2)}) & \text{if $(A)$ and $m = 1$ and $d$ even,}\\
            K(\zeta_{4d}) & \text{if $(B)$ or $(C)$,}\\
            K(\zeta_{2^\ell} + \zeta_{2^\ell}\inv) & \text{if $(D)$.}
        \end{cases}
    \]
\end{theorem}

\begin{proof}
    Suppose we are not in case $(D)$.
    The lower bound on $\wh K_f$ comes from \eqref{eqn: constant field simple lower bound}.
    Recall that $\barb f = \lprof c_1, \ldots, c_n\rprof$ and that $c_\infty = c_1c_2\cdots c_n$ is the image of $\gamma_\infty$ in $\barb f$.
    Let $\Mpre = \Mpre(d,m,n,\omega)$.
    Corollary~\ref{corollary: final IMG calculation} implies that there exists a $w \in [C_d]^\infty$ and $u_i \in \Mpre$ such that
    \(
        wc_iw\inv = u_i b_i u_i\inv
    \)
    for each $1 \leq i \leq n$.
    Therefore
    \[
        wc_\infty w\inv = (u_1b_1u_1\inv)\cdots (u_n b_n u_n\inv).
    \]
    Note that $wc_\infty w\inv$ is a strict odometer in $\Mpre$, hence is conjugate to $b_\infty := b_1b_2\cdots b_n$ in $[C_d]^\infty$.
    Proposition~\ref{prop: N(B) odometer criteria} implies that $w c_\infty w\inv \sim b_\infty$ in $N(\Mpre)$.
    Altogether this implies that $c_\infty \sim_\ell c_\infty^\ve$ in $\barb f$ if and only if $b_\infty \sim_\ell b_\infty^\ve$ in $\Mpre$.

    Proposition~\ref{prop: c_ve in B} implies that $b_\infty \sim_\ell b_\infty^\ve$ in $\Mpre$ for all $\ve \in 1 + d\ZZ_d$ when $\ell \leq n$, and $b_\infty \sim_\ell b_\infty^\ve$ in $\Mpre$ for $\ell > n$ if and only if $\ve \equiv 1 \bmod \kappa$ where
    \[
        \kappa := 
        \begin{cases}
            d^2/\gcd(d,\omega) & \text{if $(A)$ and either $m > 1$ or $d$ odd,}\\
            d^2/\gcd(d,\omega + d/2) & \text{if $(A)$ and $m = 1$ and $d$ even,}\\
            4d & \text{if $(B)$ or $(C)$}.
        \end{cases}
    \]
    This translates directly via Theorem \ref{theorem: outer action determined by infinity} to the calculations for $\wh K_f$ in cases $(A)$, $(B)$, and $(C)$.

    In case $(D)$, $\Mpre$ is isomorphic to the pro-2 dihedral group (Lemma~\ref{lemma: N quotients}).
    We have $b_\infty \sim_\ell b_\infty^\ve$ in $\Mpre$ if and only if $\ve \equiv -1 \bmod 2^\ell$.
    Hence in this case $\wh K_{f,\ell} = K(\zeta_{2^\ell} + \zeta_{2^\ell}\inv)$.
\end{proof}

\appendix

\section{Computer calculations}
\label{calculation: gap code}

We use Laurent Bartholdi's FR GAP package~\cite{FR2.4.13} to verify calculations in case $(d,m,n,\omega) = (2,2,3,1)$ (case $(C)$).
The code below calculates the indices $[\Mpre : \br]_3$, $[\Mpre : \br]_4$, and the index $[\Mpre : \br'\cap \Mpre]_4$ where $\br'$ is the normal closure of $\br$ in $[C_2]^\infty$.

\vspace{1em}

\begin{lstlisting}[language=GAP]
LoadPackage("fr");

# Construct the (discrete) IMG
B := FRGroup("b1=(1,2)", "b2=<1,b1>", "b3=<b3,b2>");
AssignGeneratorVariables(B);

# Branching subgroup
N := NormalClosure(B, [b1, Comm(b2,Comm(b2,b3)), Comm(b3,Comm(b2,b3))]);

# Level 3 and 4 truncations
B3 := PermGroup(B,3);
N3 := PermGroup(N,3);
B4 := PermGroup(B,4);
N4 := PermGroup(N,4);

# Truncated $[C_2]^4$ via FRGroups to ensure compatibility
C4 := PermGroup(FRGroup("c1=(1,2)", "c2=<1,c1>", "c3=<1,c2>", "c4=<1,c3>"),4);

# Normalizer of $B_4$ in $[C_2]^4$
NB4 := Normalizer(C4, B4);

# Normal closure of the branching subgroup in $[C_2]^4$
Np := NormalClosure(C4,N4);

# Results
Print("[B:N]_4 = ", Index(B4,N4),"\n",
      "[B:N]_3 = ", Index(B3,N3), "\n",
      "[B:N'\\cap B]_4 = ", Index(B4,Intersection(Np,B4)),"\n",
      "N(B)_4 normalizes N_4? ",IsNormal(NB4,N4));
\end{lstlisting}

\newpage
\printbibliography

\end{document}